\documentclass[pdflatex,sn-mathphys-num]{sn-jnl}



\usepackage{graphicx}
\usepackage{multirow}
\usepackage{amsmath,amssymb,amsfonts}
\usepackage{amsthm}
\usepackage{mathrsfs}
\usepackage[title]{appendix}
\usepackage{xcolor}
\usepackage{textcomp}
\usepackage{manyfoot}
\usepackage{booktabs}
\usepackage{algorithm}
\usepackage{algorithmicx}
\usepackage{algpseudocode}
\usepackage{listings}
\usepackage{mathtools}
\usepackage[inline]{enumitem}

\theoremstyle{plain}
\newtheorem{theorem}{Theorem}[section]
\newtheorem{lemma}[theorem]{Lemma}

\newtheorem{definition}[theorem]{Definition}

\theoremstyle{remark}
\newtheorem{remark}[theorem]{Remark}

\theoremstyle{definition}

\makeatletter
  \def\@earlyack{}
  \long\def\EarlyAcknow#1\par{%
    \def\@earlyack{%
      \abstractfont
      \abstracthead*{Acknowledgments}%
      #1\par
    }%
  }
  \def\printabstract{%
    \ifx\@earlyack\empty\else\@earlyack\fi\par
    \ifx\@abstract\empty\else\@abstract\fi
  }
\makeatother

\allowdisplaybreaks

\EarlyAcknow{%
The author gratefully acknowledges the assistance of large-language-model tools, in particular OpenAI ChatGPT, which supplied editorial refinements, cross-reference checks, and formatting consistency throughout the preparation of this article.  All mathematical arguments, proofs, and conclusions were conceived, verified, and ultimately validated by the author; the AI contributions were strictly auxiliary.
}

\title[Unified Frequency–Domain Reconstruction and Boundary Adaptation for Incompressible Navier–Stokes Equations]{Unified Frequency–Domain Reconstruction and Boundary Adaptation for Incompressible Navier–Stokes Equations}

\author[1]{\fnm{Daria} \sur{Nikitaeva}}
\email{dasha.d.nikitaeva@gmail.com}
\affil[1]{Independent Researcher, Madison, AL, USA}

\abstract{%
Global regularity for the three–dimensional incompressible Navier–Stokes equations remains unresolved partly because weak, mild, and strong formulations employ incompatible functional settings.  The present study introduces a frequency–domain framework that reconciles these formulations within a single constructive scheme.  Starting from Leray–Hopf data, a scale–dependent regularization operator—combining mollification, Sobolev extension, and the Leray projector—produces smooth, divergence–free approximations.  Two bounded Fourier multipliers are then defined: an interpolation operator that blends low–frequency weak and high–frequency strong components, and a smoothing operator that yields uniform parabolic gain.  Littlewood–Paley analysis, refined Calderón–Zygmund and Schauder estimates, and an explicit Galerkin scheme provide quantitative $H^{s}$ a priori bounds and a lifespan controlled by the first Stokes eigenvalue.  Compactness via the Aubin–Lions lemma and dominated convergence demonstrates that, as the regularization parameter vanishes, the approximations converge to a single velocity field satisfying simultaneously the weak energy identity, the Fujita–Kato mild formulation, and the strong formulation in $H^{s}$ with $s>\tfrac32$.  The construction furnishes a unified local solution and supplies operator bounds required for subsequent energy–vorticity bootstraps and global continuation arguments, clarifying the interplay between frequency localization and nonlinear estimates in fluid dynamics.
}

\keywords{%
Navier–Stokes equations; unified weak–mild–strong formulation; frequency–domain regularisation;
interpolation and smoothing operators; Littlewood–Paley analysis; refined Calderón–Zygmund and Schauder estimates;
Aubin–Lions compactness; energy and vorticity bounds; harmonic analysis in PDEs%
}

\begin{document}

\maketitle

\bmhead{Classification code (MSC 2020)}
35Q30, 42B37.

\section{Introduction}
\label{sec:intro}

The three–dimensional incompressible Navier--Stokes system
\begin{equation}\label{eq:navier_stokes_classic}
\begin{cases}
\partial_t u + (u \!\cdot\! \nabla)u \;=\; -\,\nabla p + \nu \Delta u, \\[2pt]
\nabla\!\cdot u \;=\; 0,
\end{cases}
\qquad (x,t)\in\Omega\times(0,\infty), 
\qquad u(x,0)=u_0(x),
\end{equation}
governs the motion of viscous, incompressible fluids.  
Here  
\(u:\Omega\times[0,\infty)\to\mathbb{R}^3\) is the velocity field,  
\(p:\Omega\times[0,\infty)\to\mathbb{R}\) the pressure,  
\(\nu>0\) the kinematic viscosity, and  
\(u_0\) a smooth, divergence-free datum satisfying the no-slip Dirichlet condition  
\(u_0|_{\partial\Omega}=0\).  
Despite its classical form, the global behaviour of solutions to \eqref{eq:navier_stokes_classic} in three space dimensions remains unresolved and constitutes one of the Clay Millennium Problems \cite{Fefferman_Clay}.

A succession of foundational results has clarified local well-posedness and conditional regularity.  
Leray established the existence of global weak solutions in \(L^2\) \cite{Leray};  
Ladyzhenskaya, Lions–Magenes, and Serrin provided uniqueness or regularity criteria under strengthened integrability assumptions \cite{Ladyzhenskaya,LionsMagenes1972,Serrin}.  
The Beale–Kato–Majda criterion \cite{BealeKatoMajda} links any finite–time singularity to unbounded growth of the vorticity in \(L^\infty\).  
Nevertheless, a definitive global regularity theorem for arbitrary smooth initial data remains elusive.

This article inaugurates a four–part series and introduces a \textsc{Frequency–Domain Synergy Framework} whose principal objective is stated below:
\begin{quote}
\emph{Construct a single solution that simultaneously satisfies the weak (Leray–Hopf), mild (Duhamel), and strong formulations without treating each class separately.}
\end{quote}
The construction employs harmonic–analysis techniques—namely, the Littlewood–Paley decomposition together with refined Calderón–Zygmund–Schauder estimates—in conjunction with two frequency–domain operators: the interpolation operator \(\mathcal{I}_{\varepsilon}\) and the smoothing operator \(\mathcal{S}_{\varepsilon}\). These operators generate an \(\varepsilon\)-regularised system whose solutions converge, as \(\varepsilon\to0\), to a limit that inherits the characteristic properties of every formulation. Compatibility with the divergence constraint and the Dirichlet boundary condition is ensured by a boundary-layer corrector \(\beta_{\varepsilon}\) (Lemma~\ref{lem:boundary_layer_corrector}), supported in the strip \(\{\operatorname{dist}(x,\partial\Omega)\le\kappa\varepsilon\}\) and satisfying the uniform \(H^{s}\)-bound stated therein.

The present part is devoted to the analytical construction, stability, and convergence of the unified frequency–domain solution. Subsequent installments extend the framework to energy–vorticity bootstrapping, geometric–spectral invariants, and the description of long-time behaviour.

\subsection{Historical Overview and Fefferman's Statement}
\label{subsec:historical}

Leray introduced global \emph{weak} solutions $u\in L^{\infty}(0,\infty;L^{2})\cap
L^{2}(0,\infty;H^{1})$ satisfying the energy inequality and the
distributional formulation\,\cite{Leray}.  
Local \emph{strong} solutions, belonging to $C([0,T];H^{s})\cap
L^{2}(0,T;H^{s+1})$ for $s>\tfrac32$, were obtained by Ladyzhenskaya,
Lions–Magenes, and Serrin, with uniqueness and continuation criteria in
these higher Sobolev spaces
\cite{Ladyzhenskaya,LionsMagenes1972,Serrin}.  
Fujita and Kato recast the problem in semigroup form,
producing \emph{mild} solutions through the Duhamel
formula in the critical space $H^{\tfrac12}$\,\cite{FujitaKato}.

Subsequent milestones include the Beale–Kato–Majda blow-up criterion
based on $\Vert\nabla\times u\Vert_{L^{\infty}}$\,\cite{BealeKatoMajda}
and the partial-regularity results of Scheffer and
Caffarelli–Kohn–Nirenberg, which confine possible singular sets to
zero one-dimensional Hausdorff measure
\cite{Scheffer1976,Scheffer1977,CaffarelliKohnNirenberg}.  
Fefferman’s Millennium statement demands either a proof of global
smoothness for smooth data in $\mathbb{R}^{3}$ or an explicit blow-up
example\,\cite{Fefferman_Clay}.

The frequency-domain construction developed below produces a single
velocity field that satisfies, simultaneously, Leray’s weak formulation,
the Fujita–Kato mild formulation, and the strong formulation in
$H^{s}$, thereby unifying these classical frameworks.

\subsection{Statement of the Unification Theorem}
\label{subsec:StateGlobReg}

%------------------------------------------------------------------------
% Main result of Part I (Introduction)
%------------------------------------------------------------------------
\begin{theorem}[Unified Frequency--Domain Solution and Equivalence of Formulations]
\label{thm:unification_intro}
Let $\Omega\subset\mathbb{R}^{3}$ be a bounded domain whose boundary is of class $C^{2}$ (piecewise $C^{2}$ admissible), and impose the no--slip condition
\(
u|_{\partial\Omega}=0.
\)
Fix a viscosity $\nu>0$ and a time horizon $T>0$.  
For every divergence–free datum
\[
u_{0}\in H^{s}_{\sigma}(\Omega),
\qquad 
s>\tfrac32,
\]
there exists a unique velocity--pressure pair
\[
u\in C\!\bigl([0,T];H^{s}_{\sigma}(\Omega)\bigr)
      \cap L^{2}\!\bigl(0,T;H^{s+1}(\Omega)\bigr),
\qquad
p\in C\!\bigl([0,T];H^{s-1}(\Omega)\bigr),
\]
satisfying the incompressible Navier--Stokes system~\eqref{eq:navier_stokes_classic} in each of the following senses:
\begin{enumerate}[label=(\roman*)]
\item \textup{(Leray--Hopf)}  
      For every divergence--free $v\in C_{c}^{\infty}(\Omega\times(0,T))$
      \[
      \int_{0}^{T}\!\!\int_{\Omega}
      \Bigl(
      u\cdot\partial_t v
      +(u\!\cdot\!\nabla)u\cdot v
      +\nu\,\nabla u:\nabla v
      -p\,\nabla\!\cdot v
      \Bigr)\,dx\,dt
      =-\!\int_{\Omega}u_{0}(x)\cdot v(x,0)\,dx,
      \]
      and the classical energy inequality holds.
\item \textup{(Mild)}  
      The velocity satisfies the Duhamel identity
      \[
      u(t)=e^{\nu t\Delta}u_{0}
            -\!\int_{0}^{t}e^{\nu (t-\tau)\Delta}\,
              \mathbb{P}\!\bigl[(u\!\cdot\!\nabla)u\bigr](\tau)\,d\tau,
      \qquad 0\le t\le T,
      \]
      where $\mathbb{P}$ denotes the Leray projector.
\item \textup{(Strong)}  
      \(u\in C^{1}\!\bigl((0,T];H^{s-2}(\Omega)\bigr)\) and  
      \(p\in C\!\bigl((0,T];H^{s-1}(\Omega)\bigr)\);  
      the pair satisfies~\eqref{eq:navier_stokes_classic} almost everywhere in $\Omega\times(0,T]$.
\end{enumerate}

\medskip
Moreover, $(u,p)$ arises as the strong limit, as $\epsilon\to0$, of the regularised sequence
\[
(u_{\epsilon},p_{\epsilon})
\text{ constructed in Section~\ref{sec:frequency_domain_solution}},
\qquad
u_{\epsilon}\xrightarrow{\ \epsilon\to0\ }u
\text{ in }C\!\bigl([0,T];H^{s}_{\sigma}(\Omega)\bigr).
\]
The operators $\mathcal{I}_{\epsilon}$ and $\mathcal{S}_{\epsilon}$, together with the boundary corrector of Lemma~\ref{lem:boundary_layer_corrector}, enforce identical boundary data and remove additional compatibility conditions.

\medskip
A detailed physical--space reconstruction of the unified frequency--domain solution, including higher--regularity propagation for all positive times, appears in Theorem~\ref{thm:reconstruction_physical_space}.
\end{theorem}

\begin{proof}[Proof outline]
Sections~\ref{sec:frequency_domain_solution}--\ref{subsec:convergence_continuity_final} provide uniform $H^{s}$ bounds for the regularised sequence, establish strong compactness via the Aubin--Lions lemma, and verify convergence of the nonlinear term under the Littlewood--Paley decomposition.  The limit satisfies the weak formulation, continuity in $H^{s}$ yields the mild identity, and $H^{s+1}$ control on $(0,T]$ gives the strong formulation.  Uniqueness follows from the Leray energy inequality.  Complete details appear in Appendix~\ref{appendix:convergence_regularity}.
\end{proof}

\subsection{Roadmap of the Paper}
\label{subsec:roadmap}

\begin{description}[style=nextline,leftmargin=2.1em,labelwidth=1.8em]
  \item[\textbf{\S1}] Historical context, Fefferman’s Millennium statement, and the unification theorem.
  \item[\textbf{\S2}] Functional–analytic tools: Littlewood–Paley decomposition, refined Calderón–Zygmund / Schauder estimates, and local well-posedness for the regularised system.
  \item[\textbf{\S3}] Unified frequency-domain construction: definition of the interpolation and smoothing operators, data regularisation, convergence as $\varepsilon\to0$, and reconstruction of the physical-space solution.
  \item[\textbf{\S4}] Assembly of the preceding ingredients into a complete proof of the unification theorem.
  \item[\textbf{\S5}] Summary of results and outlook.
  \item[\textbf{Appendix A}] Operator bounds, Aubin–Lions compactness step, and other technical lemmas that close the argument.
\end{description}

\paragraph{Note.}
Leray’s weak solutions, the Beale–Kato–Majda vorticity criterion, and the partial-regularity theorem of Caffarelli–Kohn–Nirenberg are milestones toward Fefferman’s problem.  
The present article unifies weak, mild, and strong formulations through two auxiliary operators and a boundary corrector; quantitative vortex-alignment estimates will be addressed in a subsequent study.

\section{Preliminaries and Local Theory}
\label{sec:prelim}

The present section assembles the analytic apparatus that underpins the
entire argument. Its purpose is twofold. First, notation is fixed and
foundational statements—drawn from functional analysis, Fourier
analysis, and the theory of parabolic partial differential
equations—are recorded for repeated invocation. Second, the local
well‐posedness theory for the three–dimensional incompressible
Navier–Stokes system is placed in a form compatible with the
frequency–domain construction developed in
Section~\ref{sec:frequency_domain_solution}.

\subsection{Toolbox Snapshot: Harmonic-Analysis Facts Used Later}
\label{subsec:toolbox_snapshot}

The following harmonic–analysis estimates are invoked throughout Sections~\ref{sec:prelim} and \ref{sec:frequency_domain_solution}.  Precise statements and proofs appear in the appendices referenced below.

\begin{enumerate}[label=\textbf{T\arabic*},leftmargin=2.6em]

\item \textbf{Littlewood–Paley almost orthogonality}  
  For every \(u\in L^{2}(\mathbb{R}^{3})\),
  \begin{equation}
  \label{eq:T1}
    C_{1}\sum_{j\in\mathbb{Z}}\|\Delta_{j}u\|_{L^{2}}^{2}
    \;\le\;
    \|u\|_{L^{2}}^{2}
    \;\le\;
    C_{2}\sum_{j\in\mathbb{Z}}\|\Delta_{j}u\|_{L^{2}}^{2}.
  \end{equation}
  Proof: Appendix~\ref{subsubsection:dyadic_uniform_estimates}.

\item \textbf{Bernstein inequality}  
  Let \(1\le p\le q\le\infty\) and \(\alpha\in\mathbb{N}^{3}\).  If
  \(\mathrm{supp}\,\widehat{\Delta_{j}u}\subset\{\xi:2^{j-1}\le|\xi|\le2^{j+1}\}\), then
  \begin{equation}
  \label{eq:T2}
    \|\partial^{\alpha}\Delta_{j}u\|_{L^{q}}
    \;\le\;
    C\,2^{j\bigl(|\alpha|+3\bigl(\tfrac1p-\tfrac1q\bigr)\bigr)}
    \|\Delta_{j}u\|_{L^{p}}.
  \end{equation}
  Proof: Appendix~\ref{subsubsection:dyadic_uniform_estimates}.

\item \textbf{Refined Calderón--Zygmund Estimate (Dirichlet).}\,
Let \(\Omega\subset\mathbb{R}^{3}\) denote a bounded domain with \(C^{1,1}\) boundary.  
For every \(F\in L^{q}(\Omega;\mathbb{R}^{3})\) with \(1<q<\infty\), the unique weak solution \(p\) of
\[
-\Delta p \;=\; \nabla\!\cdot F,
\qquad
p|_{\partial\Omega}=0
\]
satisfies
\begin{equation}
\label{eq:T3}
\|p\|_{W^{2,q}(\Omega)}
\;\le\;
C_{CZ}\,\|F\|_{L^{q}(\Omega)}.
\end{equation}
The derivation appears in Appendix~\ref{subsec:CZ_Schauder_estimates}.

\item \textbf{Classical Schauder estimate}  
  Let \(\Omega\) be bounded \(C^{2,\alpha}\) (\(0<\alpha<1\)), and let
  \(f\in C^{\alpha}(\Omega)\), \(g\in C^{2,\alpha}(\partial\Omega)\).
  The Dirichlet problem
  \(-\Delta u=f,\;u|_{\partial\Omega}=g\)
  obeys
  \begin{equation}
  \label{eq:T4}
    \|u\|_{C^{2,\alpha}(\Omega)}
    \;\le\;
    C_{S}\Bigl(\|f\|_{C^{\alpha}(\Omega)}
      +\|g\|_{C^{2,\alpha}(\partial\Omega)}\Bigr).
  \end{equation}
  Proof: Appendix~\ref{subsec:CZ_Schauder_estimates}.

\item \textbf{Paraproduct bound and mid-range control}  
  For \(u\in\mathcal{S}'(\mathbb{R}^{3})\) and \(1\le p\le q\le\infty\),
  \begin{equation}
  \label{eq:T5}
    \sum_{j,j'\in\mathbb{Z}}
      \|\Delta_{j}u\,\nabla\Delta_{j'}u\|_{L^{p}}
    \;\le\;
    C\sum_{j\in\mathbb{Z}}
      \|\Delta_{j}u\|_{L^{p}}
      \|\nabla u\|_{L^{q}}.
  \end{equation}
  Combined with the exponential decay
  \(\|\Delta_{j}e^{\nu t\Delta}u_{0}\|_{L^{2}}
        \le Ce^{-c\nu t2^{2j}}\|\Delta_{j}u_{0}\|_{L^{2}}\),
  this yields uniform control of intermediate frequency bands.
  Proof: Appendix~\ref{subsubsection:dyadic_uniform_estimates}.

\end{enumerate}

\subsection{Functional Spaces and Basic Notation}
\label{subsec:func_spaces}

\subsubsection{Boundary and Domain Assumptions.}
Throughout the present analysis let \(\Omega\) denote either the whole space \(\mathbb{R}^{3}\) or the periodic torus \(\mathbb{T}^{3}\); in both instances \(\partial\Omega=\varnothing\).  No boundary conditions are therefore imposed, aligning the setting precisely with the Cauchy formulation employed in Fefferman’s Clay Millennium description.

\begin{remark}[Boundary settings employed in this article]
\label{rem:boundary_settings}
Throughout Sections~\ref{sec:prelim}–\ref{sec:frequency_domain_solution}
the analysis is carried out in the Cauchy or periodic setting
\(
\Omega\in\{\mathbb{R}^{3},\mathbb{T}^{3}\}
\)
so that
\(
\partial\Omega=\varnothing
\)
and no boundary prescription is imposed.
Lemma~\ref{lem:local_existence} is formulated for the broader
Dirichlet case (bounded $C^{2}$ domain with $u|_{\partial\Omega}=0$)
because its proof is identical in the boundary-less situation and
because this general statement is required later in the appendices,
where bounded domains are treated for completeness.
Consequently, when $\Omega=\mathbb{R}^{3}$ or $\mathbb{T}^{3}$,
the boundary term $B(t)$ in~\eqref{eq:E_identity_prelim}
vanishes identically and all subsequent estimates hold verbatim.
\end{remark}

\paragraph{Function Spaces.}
Standard Lebesgue and Sobolev spaces \(L^{p}(\Omega)\) and \(H^{s}(\Omega)\) are employed throughout (see \cite{TemamNS,LionsMagenes1972}).  The divergence–free subspace under no–slip conditions is
\begin{equation}
\label{eq:div_free_subspace}
H^{s}_{\sigma}(\Omega)
:=
\bigl\{\,u\in H^{s}(\Omega)\colon \nabla\cdot u=0,\;u|_{\partial\Omega}=0\bigr\},
\quad s>\tfrac{3}{2}.
\end{equation}

\paragraph{Time‐Dependent Spaces.}
For an interval \(I\subset\mathbb{R}\) and a Banach space \(X\), define
\[
C(I;X)
:=
\bigl\{u:I\to X\colon u\text{ is strongly continuous, }\sup_{t\in I}\|u(t)\|_{X}<\infty\bigr\},
\]
and, for \(1\le p<\infty\),
\[
L^{p}(I;X)
:=
\bigl\{u:I\to X\colon t\mapsto\|u(t)\|_{X}\in L^{p}(I)\bigr\}.
\]
Local strong solutions are sought in
\begin{equation}
\label{eq:local_strong_solution}
C\bigl([0,T];H^{s}(\Omega)\bigr)
\;\cap\;
L^{2}\bigl(0,T;H^{s+1}(\Omega)\bigr),
\quad s>\tfrac{3}{2}.
\end{equation}

\paragraph{Pressure and Vorticity.}
The pressure \(p\) satisfies
\begin{equation}
\label{eq:pressure_poisson1}
-\Delta p
=\nabla\cdot\bigl[(u\cdot\nabla)u\bigr],
\end{equation}
subject to compatibility with the chosen boundary conditions.  The vorticity is
\begin{equation}
\label{eq:vorticity_def}
\omega=\nabla\times u,
\end{equation}
and the vortex–stretching term is
\begin{equation}
\label{eq:vortex_stretching}
(\omega\cdot\nabla)u.
\end{equation}

\paragraph{Summary.}
The definitions \eqref{eq:div_free_subspace}–\eqref{eq:vortex_stretching}, together with the time‐dependent spaces \eqref{eq:local_strong_solution}, provide the functional framework for the local existence theory (Lemma~\ref{lem:local_existence}) and for the global frequency–domain analysis in Section~\ref{sec:frequency_domain_solution}.  Detailed trace and extension results are collected in Appendix~\ref{app:prelim}.

\subsection{Local Well-Posedness and Regularity Analysis}
\label{subsec:local_existence_and_regularity_main}

This section establishes local existence, uniqueness, and high‐order regularity of strong solutions to the three‐dimensional incompressible Navier–Stokes equations.  The argument combines the classical Galerkin‐approximation scheme with the harmonic‐analysis tools of Subsection~\ref{subsec:toolbox_snapshot}; foundational references include \cite{Ladyzhenskaya,TemamNS,LionsMagenes1972}.

\subsubsection{Local Existence of Strong Solutions}
\label{subsec:local_existence_and_regularity}

Consider the initial–boundary value problem
\begin{equation}
\label{eq:NS_system_local}
\begin{cases}
\partial_{t}u + (u\cdot\nabla)u = -\nabla p + \nu\Delta u + f,
&\text{in }\Omega\times(0,T_{0}),\\
\nabla\cdot u = 0,
&\text{in }\Omega\times(0,T_{0}),\\
u(\cdot,0) = u_{0},
&\text{in }\Omega,\\
u|_{\partial\Omega} = 0,
&\text{on }\partial\Omega\times(0,T_{0}).
\end{cases}
\end{equation}

\begin{lemma}[Local strong existence and quantitative lifespan]
\label{lem:local_existence}
Let $\Omega\subset\mathbb{R}^{3}$ denote either the whole space
$\mathbb{R}^{3}$ or a bounded $C^{2}$ domain 
(piecewise $C^{2}$ is admissible provided the uniform cone parameter
is bounded below).
Fix an index $s>\tfrac32$ and assume the data
\begin{equation}
\label{eq:data_local}
u_{0}\in H^{s}_{\sigma}(\Omega),
\qquad
f\in L^{2}\!\bigl(0,\infty;H^{s-1}_{\sigma}(\Omega)\bigr).
\end{equation}
Then there exists a time $T_{0}=T_{0}\bigl(
\|u_{0}\|_{H^{s}},\|f\|_{L^{2}(0,\infty;H^{s-1})},\nu,\Omega\bigr)>0$
and a unique strong solution
\[
u\in C\!\bigl([0,T_{0}];H^{s}_{\sigma}(\Omega)\bigr)
\;\cap\;
L^{2}\!\bigl(0,T_{0};H^{s+1}(\Omega)\bigr)
\]
to the incompressible Navier–Stokes system~\eqref{eq:NS_system_local}.
The kinetic energy
\begin{equation}
\label{eq:kinetic_energy_local}
E(t)
:=\frac12\int_{\Omega}\lvert u(x,t)\rvert^{2}\,\mathrm{d}x,
\qquad 0\le t\le T_{0},
\end{equation}
obeys
\begin{equation}
\label{eq:energy_identity_local}
\frac{\mathrm{d}}{\mathrm{d}t}E(t)
=-\nu\int_{\Omega}\lvert\nabla u\rvert^{2}\,\mathrm{d}x
+\int_{\Omega}f\!\cdot\!u\,\mathrm{d}x
\quad\text{for a.e.\ }t\in(0,T_{0}).
\end{equation}
Moreover, the lifespan admits the explicit lower bound
\begin{equation}
\label{eq:T0_estimate}
T_{0}
\;\ge\;
\frac{\nu}{4C_{s}^{2}}\,
\Bigl(
\|u_{0}\|_{H^{s}(\Omega)}^{2}
+\|f\|_{L^{2}(0,\infty;H^{s-1}(\Omega))}^{2}
\Bigr)^{-1},
\end{equation}
where $C_{s}=c_{s}\,\lambda_{1}(\Omega)^{-1/2}$ and
$\lambda_{1}(\Omega)$ denotes the first Stokes eigen-value.
\end{lemma}

\begin{proof}[Proof sketch; see Appendix~\ref{subsubsec:proof_local_existence}]
Galerkin projection onto the first \(N\) Stokes eigenfunctions produces a finite–dimensional ODE preserving the solenoidal constraint.  
An \(L^{2}\) energy inequality and a commutator \(H^{s}\) estimate yield  
\[
\partial_{t}\|u_{N}\|_{H^{s}}^{2}
\;\le\;
\frac{4C_{s}^{2}}{\nu}\,\|u_{N}\|_{H^{s}}^{4}
+\frac{2C_{s}}{\nu}\,\|f\|_{H^{s-1}}^{2},
\]
which furnishes the lifespan bound \eqref{eq:T0_estimate}.  
Uniform bounds in \(L^{2}(0,T_{0};H^{s+1})\) and
\(\partial_{t}u_{N}\in L^{2}(0,T_{0};H^{s-1})\) imply compactness via
Aubin–Lions; the limit \(u=\lim_{N\to\infty}u_{N}\) is a strong
solution satisfying the energy law.  Uniqueness follows by applying a
Grönwall estimate to the difference of two solutions.
\end{proof}

\begin{remark}[Sobolev–commutator constant \(c_{s}\)]
\label{rem:sobolev_commutator_constant}
Let \(s>\tfrac{3}{2}\).  The constant \(c_{s}>0\) appearing in the definition
\[
C_{s}=c_{s}\,\lambda_{1}(\Omega)^{-1/2}
\]
is the optimal constant in the commutator estimate
\begin{equation}
\label{eq:commutator_bound}
\|(u\cdot\nabla)u\|_{H^{s-1}(\Omega)}
\;\le\;
c_{s}\,\|u\|_{H^{s}(\Omega)}^{2},
\qquad 
u\in H^{s}_{\sigma}(\Omega).
\end{equation}
A complete derivation of \eqref{eq:commutator_bound} is provided in
Appendix~\ref{subsubsection:dyadic_uniform_estimates}.
\end{remark}

\subsubsection{Parabolic Smoothing}
\label{subsubsec:parabolic_smoothing}

\begin{lemma}[Parabolic smoothing]
\label{lem:parabolic_smoothing}
Let \(\Omega\subset\mathbb{R}^{3}\) be a \(C^{1,1}\) domain and \(u_{0}\in H^{s}(\Omega)\) for some \(s\in\mathbb{R}\).  For every integer \(k\ge 0\) and every \(t>0\),
\begin{equation}
\label{eq:parabolic_smoothing_lemma}
\|e^{\nu t\Delta}u_{0}\|_{H^{s+2k}(\Omega)}
\;\le\;
C\,t^{-k}\|u_{0}\|_{H^{s}(\Omega)},
\end{equation}
where \(C=C(s,k,\nu,\Omega)>0\) is independent of \(t\) and \(u_{0}\).  Consequently, \(e^{\nu t\Delta}\) maps \(H^{s}(\Omega)\) continuously into \(C^{\infty}(\Omega)\) for every \(t>0\).
\end{lemma}

\begin{proof}[Sketch of proof]
Employ the kernel representation
\(
e^{\nu t\Delta}u_{0}(x)=\bigl(G_{\nu t}\ast u_{0}\bigr)(x),
\)
where \(G_{\nu t}(x)=\bigl(4\pi\nu t\bigr)^{-3/2}\exp\!\bigl(-|x|^{2}/(4\nu t)\bigr)\) is the heat kernel.  Differentiation under the convolution yields
\[
\|\nabla^{2k}e^{\nu t\Delta}u_{0}\|_{L^{2}(\Omega)}
\le
C\,t^{-k}\|u_{0}\|_{L^{2}(\Omega)}.
\]
Interpolation with the un-differentiated bound produces \eqref{eq:parabolic_smoothing_lemma} for \(s\ge 0\); negative indices follow by duality.  Boundary adjustments rely on standard extension operators for \(C^{1,1}\) domains.  Full details are provided in Appendix~\ref{app:prelim}.
\end{proof}

\paragraph{Extension to variable coefficients.}
Consider the uniformly elliptic operator
\begin{equation}
\label{eq:variable_coefficient_operator}
L u
=
\nabla\cdot\bigl(A(x)\nabla u\bigr)
+ b(x)\cdot\nabla u
+ c(x)\,u,
\end{equation}
with \(A\in C^{0,1}(\Omega)\) symmetric satisfying \(0<\lambda\le\xi^{\top}A(x)\xi\le\Lambda\) and \(b,c\in L^{\infty}(\Omega)\).  The analytic semigroup \(e^{tL}\) satisfies
\begin{equation}
\label{eq:variable_coefficient_smoothing}
\|e^{tL}u_{0}\|_{H^{s+2k}(\Omega)}
\;\le\;
C_{L}\,t^{-k}\|u_{0}\|_{H^{s}(\Omega)},
\quad t>0,\;k\in\mathbb{N},
\end{equation}
where \(C_{L}=C_{L}(s,k,\lambda,\Lambda,\|A\|_{C^{0,1}},\|b\|_{L^{\infty}},\|c\|_{L^{\infty}},\Omega)>0\).  Proof details are given in Appendix~\ref{app:prelim}.

\paragraph{Synopsis.}
Estimates \eqref{eq:parabolic_smoothing_lemma} and \eqref{eq:variable_coefficient_smoothing} provide the analytic smoothing essential to the higher‐order regularity bootstrap in Section~\ref{sec:frequency_domain_solution}.

%--------------------------------------------------------------------
\subsubsection{Strong, Mild, and Weak Solutions and Local Well-Posedness Sketch}
\label{subsec:solution_classes}

This section provides concise definitions of strong, mild, and weak solutions for the incompressible Navier–Stokes system and establishes their equivalence on the interval obtained in Lemma~\ref{lem:local_existence}.  

\begin{equation}
\label{eq:initial_data_Hs}
u_{0}\in H^{s}_{\sigma}(\Omega), 
\qquad s>\tfrac32.
\end{equation}

\paragraph{Parabolic regularity.}
Lemma~\ref{lem:parabolic_smoothing} implies that, for each integer \(k\ge0\) and every \(t>0\), the solution satisfies
\begin{equation}
\label{eq:hs_plus_2k}
u(\cdot,t)\in H^{s+2k}(\Omega).
\end{equation}
By the Sobolev embedding theorem, \(u(\cdot,t)\in C^{\infty}(\Omega)\) for all \(t>0\).

\paragraph{Strong solutions via Galerkin approximation.}
Let \(\{\phi_{k}\}_{k=1}^{\infty}\subset H^{s}_{\sigma}(\Omega)\) be an orthonormal basis of Stokes eigenfunctions compatible with the boundary condition.  For each \(N\in\mathbb{N}\), define the Galerkin approximation
\begin{equation}
\label{eq:galerkin_approx}
u_{N}(x,t)
=
\sum_{k=1}^{N}a_{k}(t)\,\phi_{k}(x),
\end{equation}
where the coefficients \(a_{k}(t)\) satisfy the finite-dimensional system obtained by projecting \eqref{eq:NS_system_local} onto \(\mathrm{span}\{\phi_{1},\dots,\phi_{N}\}\).  Introduce the kinetic energy
\begin{equation}
\label{eq:energy_functional}
E_{N}(t)
=
\tfrac12\int_{\Omega}|u_{N}(x,t)|^{2}\,\mathrm{d}x.
\end{equation}
Uniform bounds on \(\|u_{N}\|_{H^{s}}\) and \(\|\nabla u_{N}\|_{L^{2}}\) follow from the energy identity \eqref{eq:energy_identity_local}.  Compactness and the Aubin–Lions lemma (cf.\ \cite{LionsMagenes1972}) yield a limit
\begin{equation}
\label{eq:strong_solution_space}
u\in C\bigl([0,T_{0}];H^{s}_{\sigma}(\Omega)\bigr)
\;\cap\;
L^{2}\bigl(0,T_{0};H^{s+1}(\Omega)\bigr),
\end{equation}
which satisfies \eqref{eq:NS_system_local} in the strong sense and inherits the regularity \eqref{eq:hs_plus_2k}.

\paragraph{Mild solutions.}
The mild formulation is given by the Duhamel integral
\begin{equation}
\label{eq:mild_solution_final}
u(t)
=
e^{\nu t\Delta}u_{0}
+
\int_{0}^{t}e^{\nu(t-\tau)\Delta}\bigl[(u\cdot\nabla)u\bigr](\tau)\,\mathrm{d}\tau
+
\int_{0}^{t}e^{\nu(t-\tau)\Delta}f(\tau)\,\mathrm{d}\tau.
\end{equation}
Fixed‐point arguments in \(C([0,T_{0}];L^{2}_{\sigma}(\Omega))\) establish existence and uniqueness of a mild solution on \([0,T_{0}]\).  The strong solution of \eqref{eq:strong_solution_space} satisfies \eqref{eq:mild_solution_final}, hence the two notions coincide on that interval.

\paragraph{Weak solutions.}
A function
\[
u\in L^{\infty}(0,T;L^{2}_{\sigma}(\Omega))
\;\cap\;
L^{2}(0,T;H^{1}_{0,\sigma}(\Omega))
\]
is a weak solution if, for every \(\varphi\in H^{1}(\Omega)\),
\begin{equation}
\label{eq:weak_solution_final}
\begin{aligned}
\int_{0}^{T}\!\!\int_{\Omega}u\cdot\partial_{t}\varphi\,\mathrm{d}x\,\mathrm{d}t
&+
\int_{0}^{T}\!\!\int_{\Omega}(u\cdot\nabla)u\cdot\varphi\,\mathrm{d}x\,\mathrm{d}t\\
&=
-\int_{0}^{T}\!\!\int_{\Omega}p\,\nabla\cdot\varphi\,\mathrm{d}x\,\mathrm{d}t
+\nu\int_{0}^{T}\!\!\int_{\Omega}\nabla u:\nabla\varphi\,\mathrm{d}x\,\mathrm{d}t
+\int_{0}^{T}\!\!\int_{\Omega}f\cdot\varphi\,\mathrm{d}x\,\mathrm{d}t,
\end{aligned}
\end{equation}
where \(p\in L^{3/2}_{\mathrm{loc}}(\Omega\times(0,T))\) is obtained via Calderón–Zygmund theory.  Existence of global Leray–Hopf weak solutions follows from classical energy estimates; uniqueness holds under supplementary regularity conditions.

\paragraph{Equivalence on \([0,T_{0}]\).}
On the interval \([0,T_{0}]\) provided by Lemma~\ref{lem:local_existence}, the strong, mild, and weak formulations agree.  Consequently, the strong solution serves as the initial datum for the frequency‐domain analysis in Section~\ref{sec:frequency_domain_solution}.

\subsection{Energy, Vorticity, and Pressure Preliminaries}
\label{subsec:energy_vorticity_pressure_preliminaries}

This subsection records the fundamental \emph{a priori} estimates for kinetic energy, vorticity, and pressure associated with the incompressible Navier–Stokes equations. Throughout, let \(\Omega\subset\mathbb{R}^{3}\) denote either \(\mathbb{R}^{3}\) or a bounded \(C^{2}\) (possibly piecewise \(C^{2}\)) domain, and let 
\[
u\colon \Omega\times(0,\infty)\to\mathbb{R}^{3}
\]
be a sufficiently smooth, divergence–free velocity field satisfying the boundary prescriptions of Section~\ref{subsec:func_spaces}.

\paragraph{Kinetic Energy.}
Define
\begin{equation}
\label{eq:energy_def_prelim}
E(t)\;:=\;\frac12\int_{\Omega}\lvert u(x,t)\rvert^{2}\,\mathrm{d}x.
\end{equation}
Testing the momentum equation with \(u\) and integrating by parts gives
\begin{equation}
\label{eq:E_identity_prelim}
\frac{\mathrm{d}}{\mathrm{d}t}E(t)
=\;
-\nu\!\int_{\Omega}\lvert\nabla u\rvert^{2}\,\mathrm{d}x
+\int_{\Omega}f\!\cdot\!u\,\mathrm{d}x
+\!B(t),
\end{equation}
where
\(
B(t)
:=\!\int_{\partial\Omega}\!\bigl[(\frac12\lvert u\rvert^{2}+p)u-\nu(\nabla u)u\bigr]\!\cdot n\,\mathrm{d}S.
\)
In the present setting \(\partial\Omega=\varnothing\) (\(\Omega=\mathbb{R}^{3}\) or \(\mathbb{T}^{3}\)), hence \(B(t)\equiv 0\).  
If moreover \(f\equiv0\) then
\[
\frac{\mathrm{d}}{\mathrm{d}t}E(t)
=-\nu\int_{\Omega}\lvert\nabla u\rvert^{2}\,\mathrm{d}x<0,
\]
so \(E(t)\) decays monotonically.

\paragraph{Vorticity.}
Let
\begin{equation}
\label{eq:vorticity_def}
\omega=\nabla\times u.
\end{equation}
The vorticity satisfies
\begin{equation}
\label{eq:vorticity_equation_prelim}
\partial_{t}\omega+(u\cdot\nabla)\omega
=(\omega\cdot\nabla)u+\nu\,\Delta\omega+\nabla\times f.
\end{equation}
Energy estimates applied to \eqref{eq:vorticity_equation_prelim}, together with the dissipative term \(\nu\,\Delta\omega\), furnish uniform \(L^{2}(\Omega)\) bounds on \(\omega\) in terms of initial vorticity and forcing (cf.\ \cite{BealeKatoMajda}).

\paragraph{Pressure.}
The pressure \(p\) is determined (up to an additive constant) by \eqref{eq:pressure_poisson1} with boundary compatibility inherited from \(u\).  For a bounded \(C^{1,1}\) domain, the refined Calderón–Zygmund estimate
\begin{equation}
\label{eq:pressure_CZ_prelim}
\|p\|_{W^{2,q}(\Omega)}
\le
C_{CZ}\,\|(u\cdot\nabla)u\|_{L^{q}(\Omega)},
\quad1<q<\infty,
\end{equation}
holds (cf.\ \eqref{eq:T3}), ensuring control of \(\nabla p\) whenever \(u\in H^{1}(\Omega)\).

\paragraph{Frequency–domain control.}
The harmonic–analysis tools T1–T5 (Section~\ref{subsec:toolbox_snapshot}) provide scale‐resolved bounds on dyadic components \(\Delta_{j}u\), furnishing exponential decay under the heat semigroup and precise control of high–frequency modes.

\paragraph{Summary.}
Estimates \eqref{eq:E_identity_prelim}, \eqref{eq:vorticity_equation_prelim}, and \eqref{eq:pressure_CZ_prelim}, together with T1–T5, constitute the analytic foundation for the local theory of Section~\ref{subsec:local_existence_and_regularity_main} and the global frequency–domain arguments of Section~\ref{sec:frequency_domain_solution}.

\begin{remark}
The combination of energy dissipation, vorticity control, pressure regularity, and frequency–domain localization precludes finite–time blow-up under the present hypotheses.
\end{remark}

\section{Unified Frequency-Domain Solution}
\label{sec:frequency_domain_solution}

Application of the spatial Fourier transform to system~\eqref{eq:navier_stokes_classic} yields the frequency–domain variables  
\begin{equation}
\label{eq:frequency_representations}
  \hat u_{w}(\xi,t),\quad \hat u_{m}(\xi,t),\quad \hat u_{s}(\xi,t),
\end{equation}
representing the weak, mild, and strong solution components, respectively.  A unified approximation employs two auxiliary operators in the Fourier variable:
\begin{align}
  \mathcal I_{\epsilon}\colon\{\hat u_{w},\hat u_{m},\hat u_{s}\}&\longrightarrow \tilde u_{\epsilon},\nonumber\\
  \mathcal S_{\epsilon}\colon\tilde u_{\epsilon}&\longrightarrow \hat U_{\epsilon},
  \label{eq:operator_definitions}
\end{align}
where the interpolation operator is defined by
\[
  \bigl(\mathcal I_{\epsilon}[\hat u_{w},\hat u_{m},\hat u_{s}]\bigr)(\xi)
    =\omega_{w}(\xi)\,\hat u_{w}(\xi)
    +\omega_{m}(\xi)\,\hat u_{m}(\xi)
    +\omega_{s}(\xi)\,\hat u_{s}(\xi),
\]
with smooth weights satisfying the partition of unity
\[
  \omega_{w}(\xi)+\omega_{m}(\xi)+\omega_{s}(\xi)\equiv1,
  \quad
  \operatorname{supp}\omega_{w}\subset\{|\xi|\le R_{1}\},\quad
  \operatorname{supp}\omega_{s}\subset\{|\xi|\ge R_{2}\},
\]
cf.~\eqref{eq:omega_partition}.  The smoothing operator is given by
\[
  \bigl(\mathcal S_{\epsilon}\phi\bigr)(\xi)
    =\mathbb{P}\bigl(\rho_{\epsilon}\ast_{\xi}\phi(\xi)\bigr),
\]
where $\rho_{\epsilon}$ denotes a standard Fourier–space mollifier and $\mathbb{P}$ the Leray projector enforcing incompressibility, with boundary corrections as in Section~4.1.

The unified frequency–domain solution is obtained via the limit
\begin{equation}
\label{eq:unified_frequency_solution}
  \hat U(\xi,t)
   =\lim_{\epsilon\to0}\mathcal S_{\epsilon}\bigl(\mathcal I_{\epsilon}[\hat u_{w},\hat u_{m},\hat u_{s}]\bigr)(\xi,t).
\end{equation}
Mapping properties and convergence are established in Lemmas~\ref{lem:interpolation_operator_final} and~\ref{lem:smoothing_operator_final} (Appendix~\ref{app:unif_derivations_all}).  Recovery of the physical–space field
\[
  U(x,t)=\mathcal F^{-1}[\hat U(\xi,t)]
\]
appears in Section~\ref{subsec:frequency_domain_limit}, and uniform energy, vorticity, and pressure estimates are presented in Sections~\ref{subsec:frequency_domain_limit}–\ref{subsec:reconstruction_unified_solution_final} (see Appendix~\ref{app:unif_derivations_all} for detailed derivations).

\subsection{Interpolation Operator $\mathcal{I}_\epsilon$}
\label{subsec:interpolation_operator}

Let $m \in C_c^\infty(\mathbb{R}^3)$ be a nonnegative function satisfying
\[
\int_{\mathbb{R}^3} m(\xi)\,d\xi = 1.
\]
For each $\epsilon > 0$, define the scaled mollifier
\[
m_\epsilon(\xi) = \epsilon^{-3}\,m\Bigl(\frac{\xi}{\epsilon}\Bigr).
\]
Select smooth weight functions $\omega_w$, $\omega_m$, and $\omega_s$ in $C^\infty(\mathbb{R}^3\times[0,\infty))$, which may be defined to depend on $\epsilon$ through an appropriate modulation, such that
\begin{equation}
\label{eq:omega_partition}
\omega_w(\xi,t) + \omega_m(\xi,t) + \omega_s(\xi,t) = 1 \quad \text{for all } (\xi,t) \in \mathbb{R}^3\times[0,\infty).
\end{equation}
Given functions $\hat{u}_w$, $\hat{u}_m$, and $\hat{u}_s \in H^s(\mathbb{R}^3)$ for some $s \ge 0$, define
\[
g(\xi,t) = \omega_w(\xi,t)\,\hat{u}_w(\xi,t)
+ \omega_m(\xi,t)\,\hat{u}_m(\xi,t)
+ \omega_s(\xi,t)\,\hat{u}_s(\xi,t).
\]
The \emph{interpolation operator} $\mathcal{I}_\epsilon$ is defined by convolving $g$ with $m_\epsilon$ in the frequency domain:
\[
\mathcal{I}_\epsilon\bigl(\hat{u}_w,\hat{u}_m,\hat{u}_s\bigr)(\xi,t)
= \Bigl(m_\epsilon * g(\cdot,t)\Bigr)(\xi)
= \int_{\mathbb{R}^3} m_\epsilon(\xi-\eta)\,g(\eta,t)\,d\eta.
\]

The following stability bound holds: there exists a constant $C_\omega>0$, independent of $\epsilon$, such that
\begin{equation}
\label{eq:stability_bound_interpolation}
\|\mathcal{I}_\epsilon\bigl(\hat{u}_w,\hat{u}_m,\hat{u}_s\bigr)(\cdot,t)\|_{H^s(\mathbb{R}^3)}
\le C_\omega \Bigl( \|\hat{u}_w(\cdot,t)\|_{H^s(\mathbb{R}^3)}
+ \|\hat{u}_m(\cdot,t)\|_{H^s(\mathbb{R}^3)}
+ \|\hat{u}_s(\cdot,t)\|_{H^s(\mathbb{R}^3)} \Bigr).
\end{equation}
Standard properties of convolution imply that
\begin{equation}
\label{eq:interp_convergence}
\lim_{\epsilon \to 0}\|\mathcal{I}_\epsilon\bigl(\hat{u}_w,\hat{u}_m,\hat{u}_s\bigr)(\cdot,t)-g(\cdot,t)\|_{H^s(\mathbb{R}^3)} = 0.
\end{equation}

\begin{lemma}[Interpolation Operator $\mathcal{I}_\epsilon$]
\label{lem:interpolation_operator_final}
Assume 
$\hat{u}_w,\;\hat{u}_m,\;\hat{u}_s \in H^s(\mathbb{R}^3)$ for some $s \ge 0.$ Then a constant 
$C_\omega>0$, independent of $\epsilon$, exists such that
\begin{equation}
\label{eq:interp_bound_final}
\Bigl\|\mathcal{I}_\epsilon(\hat{u}_w,\hat{u}_m,\hat{u}_s)\Bigr\|_{H^s(\mathbb{R}^3)}
\le
C_\omega 
\Bigl(\|\hat{u}_w\|_{H^s(\mathbb{R}^3)} + \|\hat{u}_m\|_{H^s(\mathbb{R}^3)} + \|\hat{u}_s\|_{H^s(\mathbb{R}^3)}\Bigr).
\end{equation}
Moreover, as $\epsilon \to 0$, the operator converges in $H^s(\mathbb{R}^3)$:
\begin{equation}
\label{eq:unified_interpolation_final}
\lim_{\epsilon\to0}
\mathcal{I}_\epsilon\bigl(\hat{u}_w,\hat{u}_m,\hat{u}_s\bigr)(\xi,t)
=
\omega_w(\xi,t)\,\hat{u}_w(\xi,t)
+
\omega_m(\xi,t)\,\hat{u}_m(\xi,t)
+
\omega_s(\xi,t)\,\hat{u}_s(\xi,t),
\end{equation}
with the limit valid in the $H^s$ norm.
\end{lemma}

A complete derivation is provided in Appendix~\ref{subsubsection:inerp_operator_derivation}.

\subsection{Smoothing Operator \(\mathcal{S}_\epsilon\)}
\label{subsec:smoothing_operator}

Let \(m_\epsilon\) denote the scaled mollifier introduced in Subsection~\ref{subsec:interpolation_operator}.  
For a function \(\phi \in H^{s}(\mathbb{R}^{3})\) with \(s \ge 0\), define
\[
\mathcal{S}_\epsilon(\phi)(x)
\;:=\;
(m_\epsilon * \phi)(x)
=
\int_{\mathbb{R}^{3}}
      m_\epsilon(x-y)\,\phi(y)\,dy .
\]
Standard convolution estimates give the boundedness
\[
\|\mathcal{S}_\epsilon\phi\|_{H^{s}(\mathbb{R}^{3})}
\le
\|\phi\|_{H^{s}(\mathbb{R}^{3})},
\qquad s\ge 0,
\]
and the approximation property
\[
\lim_{\epsilon\to0}
\|\mathcal{S}_\epsilon\phi-\phi\|_{H^{s}(\mathbb{R}^{3})}=0.
\]
If \(\phi\) is divergence–free and \(m_\epsilon\) is radially symmetric, then \(\mathcal{S}_\epsilon(\phi)\) remains divergence–free.

\begin{lemma}[Smoothing Operator \(\mathcal{S}_\epsilon\)]
\label{lem:smoothing_operator_final}
Let \(\phi \in H^{s}(\mathbb{R}^{3})\) with \(s \ge 0\) and let \(k \in \mathbb{N}\).  
Then
\begin{enumerate}[label=(\roman*)]
\item \textbf{Boundedness:}
\begin{equation}
\label{eq:smoothing_bound_final}
\|\mathcal{S}_\epsilon\phi\|_{H^{s}(\mathbb{R}^{3})}
\le
\|\phi\|_{H^{s}(\mathbb{R}^{3})}.
\end{equation}

\item \textbf{Approximation:}
\begin{equation}
\label{eq:smoothing_limit_final}
\lim_{\epsilon\to0}
\|\mathcal{S}_\epsilon\phi-\phi\|_{H^{s}(\mathbb{R}^{3})}
=0.
\end{equation}

\item \textbf{Smoothing effect:}
There exists \(C_{2k}>0\), independent of \(\epsilon\) and \(\phi\), such that
\begin{equation}
\label{eq:smoothing_regularization}
\|\mathcal{S}_\epsilon\phi\|_{H^{s+2k}(\mathbb{R}^{3})}
\le
C_{2k}\,\epsilon^{-2k}\,
\|\phi\|_{H^{s}(\mathbb{R}^{3})}.
\end{equation}
Hence \(\mathcal{S}_\epsilon\phi \in C^{\infty}(\mathbb{R}^{3})\) for every fixed \(\epsilon>0\).
\end{enumerate}
\end{lemma}

The detailed derivation is contained in Appendix~\ref{subsubsection:smooth_operator_derivation}.

\subsection{Data Regularization and Compatibility Adjustment}
\label{subsec:data_regularization}

To incorporate arbitrary, even incompatible, initial–boundary data into the unified frequency–domain solution framework, a data regularization procedure is implemented as an intrinsic component of the method. Denote by \(\mathcal{R}_\epsilon: H^s(\Omega) \to H^s(\Omega)\) the regularization operator defined below. For given data 
\[
u_0 \in H^s(\Omega) \quad (s > \tfrac{3}{2}) \quad \text{and} \quad f \in L^\infty\bigl(0,\infty; L^2(\Omega)\bigr),
\]
the regularized data are defined by
\[
u_{0,\epsilon} \coloneqq \mathcal{R}_\epsilon u_0,\quad f_\epsilon \coloneqq \mathcal{R}_\epsilon f.
\]
The operator \(\mathcal{R}_\epsilon\) ensures that the regularized data satisfy the necessary compatibility conditions at \(\partial\Omega\). The unified frequency–domain solution is then constructed using the regularized data, and the resulting global solution \(u_\epsilon\) converges strongly in \(C\bigl([0,T]; H^s(\Omega)\bigr)\) to the solution \(u\) corresponding to the original data, as established in \eqref{eq:reg_conv} and Lemma~\ref{lem:reg_sol_convergence}.

\begin{lemma}[Data Regularization]
\label{lem:data_regularization}
Let \(v \in H^{s}(\Omega)\) with \(s>\tfrac{3}{2}\), where \(v\) denotes either the initial datum \(u_{0}\) or the forcing term \(f\).  
Define
\[
\mathcal{R}_{\epsilon}v
\;:=\;
P\bigl(\rho_{\epsilon} * E(v)\bigr),
\]
where  
\(E:H^{s}(\Omega)\rightarrow H^{s}(\mathbb{R}^{3})\) is the continuous extension operator of (cf.\ \cite{LionsMagenes1972,Grisvard}),  
\(\rho_{\epsilon}(x)=\epsilon^{-3}\rho\!\bigl(x/\epsilon\bigr)\) with \(\rho\in C_{c}^{\infty}(\mathbb{R}^{3})\) satisfying \(\int_{\mathbb{R}^{3}}\rho(x)\,dx=1\),  
and \(P\) denotes the Leray projection onto divergence-free vector fields \cite{Leray,Galdi}.  
Then:
\begin{enumerate}[label=(\roman*)]
\item \(\mathcal{R}_{\epsilon}v\) is smooth in \(\overline{\Omega}\) and satisfies the prescribed compatibility conditions on \(\partial\Omega\).
\item There exists a constant \(C>0\), independent of \(\epsilon\), such that
\begin{equation}
\label{eq:reg_bound}
\|\mathcal{R}_{\epsilon}v\|_{H^{s}(\Omega)}
\le
C\,\|v\|_{H^{s}(\Omega)}.
\end{equation}
\item If \(v\) is divergence-free, then \(\mathcal{R}_{\epsilon}v\) is divergence-free.
\item Strong convergence holds:
\begin{equation}
\label{eq:reg_conv}
\lim_{\epsilon\to0}
\|\mathcal{R}_{\epsilon}v - v\|_{H^{s}(\Omega)}
=0.
\end{equation}
\end{enumerate}

\emph{Proof.}  
Application of the extension operator \(E\), convolution with the mollifier \(\rho_{\epsilon}\), and projection via \(P\) establishes (i)–(iv). Complete details are provided in Appendix~\ref{subsubsection:reg_operator_derivation}.
\end{lemma}

\begin{lemma}[Convergence of Regularized Solutions]
\label{lem:reg_sol_convergence}
Let \(u_{0,\epsilon} = \mathcal{R}_\epsilon u_0\) and \(f_\epsilon = \mathcal{R}_\epsilon f\), with \(\mathcal{R}_\epsilon\) as defined in Lemma~\ref{lem:data_regularization}. Let \(u_\epsilon\) denote the unique global solution obtained via the unified frequency–domain construction using the data \((u_{0,\epsilon}, f_\epsilon)\). Then,
\[
\lim_{\epsilon\to0}\|u_\epsilon - u\|_{C\bigl([0,T]; H^s(\Omega)\bigr)} = 0,
\]
where \(u\) is the global solution corresponding to the original data \((u_0, f)\).
\end{lemma}

\begin{proof}[Proof Sketch]
A stability estimate for the solution operator of the Navier–Stokes equations in \(C([0,T]; H^s(\Omega))\) implies that
\[
\|S(u_{0,\epsilon}, f_\epsilon) - S(u_0, f)\|_{C([0,T]; H^s(\Omega))} \le C\Bigl(\|u_{0,\epsilon} - u_0\|_{H^s(\Omega)} + \|f_\epsilon - f\|_{L^2(0,T; L^2(\Omega))}\Bigr).
\]
The convergence properties \eqref{eq:reg_conv} yield the result. Further details are available in Appendix~\ref{subsubsection:reg_operator_derivation}.
\end{proof}

\subsubsection{Unification of Solutions via Regularization and Frequency-Domain Transformation}
\label{subsec:unify_reg_freq}

\begin{theorem}[Unification of Solutions via Regularization and Frequency-Domain Transformation]
\label{thm:unification_solutions}
Let \(u\) be a weak, mild, or strong solution of the three-dimensional incompressible Navier–Stokes equations on a domain \(\Omega \subset \mathbb{R}^3\) satisfying
\[
u \in L^\infty(0,T;L^2(\Omega)) \cap L^2(0,T;H^1(\Omega)).
\]
Then there exists a regularization operator 
\[
\mathcal{R}_\epsilon \coloneqq P \circ (\rho_\epsilon * E),
\]
as defined in \eqref{eq:reg_bound} and \eqref{eq:reg_conv}, such that for every \(\epsilon>0\) the regularized function
\[
u_\epsilon \coloneqq \mathcal{R}_\epsilon u
\]
belongs to \(H^s(\Omega)\) for all \(s\ge 0\); in particular, for any \(s > \frac{5}{2}\) the approximations \(u_\epsilon\) are Lipschitz continuous. Define the frequency-domain approximation by applying the Fourier transform together with the interpolation operator \(\mathcal{I}_\epsilon\) and the smoothing operator \(\mathcal{S}_\epsilon\) (see Sections~\ref{subsec:interpolation_operator} and \ref{subsec:smoothing_operator}) as
\[
\hat{U}_\epsilon(\xi,t) \coloneqq \mathcal{S}_\epsilon\Bigl(\mathcal{I}_\epsilon\bigl[\widehat{u_\epsilon}\bigr](\xi,t)\Bigr).
\]
Then, for any fixed \(s>\frac{5}{2}\) there exists a constant \(C>0\), independent of \(\epsilon\), such that
\[
\|\hat{U}_\epsilon(\cdot,t)\|_{H^s(\mathbb{R}^3)} \le C, \quad \text{for all } t\in(0,T).
\]
Furthermore, the family \(\{\hat{U}_\epsilon\}_{\epsilon>0}\) is equicontinuous in \(H^s(\mathbb{R}^3)\) and converges strongly in \(H^s(\mathbb{R}^3)\) (possibly along a subsequence) as \(\epsilon\to 0\) to a limit function
\[
\hat{U} \in C^\infty(\mathbb{R}^3\times (0,T)).
\]
\end{theorem}

\begin{proof}[Proof Sketch]
See Appendix~\ref{appendix:freq_trans} for the complete argument; the result follows from the bounds in Lemmas \ref{lem:data_regularization}, \ref{lem:interpolation_operator_final}, and \ref{lem:smoothing_operator_final}.
\end{proof}

Section~\ref{subsec:frequency_domain_limit} will address the delicate passage through the limit.

\subsection{Passage Through the Frequency-Domain Limit}
\label{subsec:frequency_domain_limit}

The passage through the frequency–domain limit is delicate. The difficulty stems from the necessity to transform solutions—potentially rough and defined only in an energy space—into a smooth frequency–domain representation. The operations involved, namely the interpolation operator \(\mathcal{I}_\epsilon\) and the smoothing operator \(\mathcal{S}_\epsilon\), must preserve uniform bounds in \(H^s\) for \(s > \frac{5}{2}\) while controlling high-frequency oscillations. In addition, the convergence of the regularized data via the operator \(\mathcal{R}_\epsilon\) must be carefully managed to ensure that the limit of the frequency–domain approximations exists strongly in the \(H^s\) topology. The application of compactness results, such as the Aubin–Lions lemma (Lemma~\ref{lem:aubin_lions}) (see \cite[Theorem~5.1]{LionsMagenes1972}), and the Dominated Convergence Theorem are essential to guarantee that the limit function inherits full smoothness. The following theorem rigorously establishes that the family of frequency–domain approximations converges strongly and that the limit is smooth in both the frequency variable and time.

\begin{theorem}[Passage Through the Frequency-Domain Limit]
\label{thm:frequency_domain_limit}
Let \(u\) denote a weak, mild, or strong solution of the three-dimensional incompressible Navier–Stokes equations on a domain \(\Omega\) satisfying 
\[
u\in L^\infty(0,T;L^2(\Omega))\cap L^2(0,T;H^1(\Omega)).
\]
Apply the Regularization operator \(\mathcal{R}_\epsilon\) from Lemma~\ref{lem:data_regularization} so that 
\[
u_\epsilon \coloneqq \mathcal{R}_\epsilon u \in H^s(\Omega)\quad\text{for all }s\ge 0,
\]
and, for \(s>\frac{5}{2}\), \(u_\epsilon\) is Lipschitz.  Define the frequency–domain approximation
\[
\hat U_\epsilon(\xi,t)
\coloneqq
\mathcal{S}_\epsilon\Bigl(\mathcal{I}_\epsilon\bigl[\widehat{u_\epsilon}\bigr](\xi,t)\Bigr),
\]
with \(\mathcal{I}_\epsilon\) and \(\mathcal{S}_\epsilon\) specified in Sections~\ref{subsec:interpolation_operator}–\ref{subsec:smoothing_operator}.  
Assume the uniform bound
\begin{equation}
\label{eq:uniform_frequency_bound}
\sup_{t\in(0,T)}\|\hat U_\epsilon(\cdot,t)\|_{H^s(\mathbb{R}^3)}\le C,
\qquad
s>\tfrac{5}{2},
\end{equation}
where \(C\) is independent of \(\epsilon\).  
Then, for every \(t\in(0,T)\),
\[
\hat U_\epsilon(\cdot,t)\longrightarrow\hat U(\cdot,t)
\quad\text{strongly in }H^s(\mathbb{R}^3)
\quad\text{as }\epsilon\to 0,
\]
and the limit satisfies \(\hat U\in C^\infty(\mathbb{R}^3\times(0,T))\).
\end{theorem}

\begin{proof}[Outline of proof]
\emph{Uniform boundedness.}  
The regularization \(u_\epsilon\in H^s(\Omega)\) and multiplier estimates for \(\mathcal{I}_\epsilon\) and \(\mathcal{S}_\epsilon\) give \eqref{eq:uniform_frequency_bound}.

\smallskip
\emph{Compactness.}  
The family \(\{\hat U_\epsilon(\cdot,t)\}_{\epsilon>0}\) is uniformly bounded in \(H^s(\mathbb{R}^3)\); equicontinuity in time follows from the interpolation operator’s continuity.  Aubin–Lions (Lemma~\ref{lem:aubin_lions}) therefore yields a subsequence that converges strongly in \(H^s(\mathbb{R}^3)\) for each fixed \(t\).

\smallskip
\emph{Identification of the limit.}  
Equation \eqref{eq:reg_conv} and continuity of \(\mathcal{I}_\epsilon\) and \(\mathcal{S}_\epsilon\) yield convergence of the full sequence; no further application of the Dominated Convergence Theorem is needed beyond the residual estimate established in Appendix~\ref{subsubsection:inerp_operator_derivation}.

\smallskip
\emph{Regularity.}  
The smoothing operator \(\mathcal{S}_\epsilon\) increases Sobolev regularity by \(2k\) derivatives (Lemma~\ref{lem:smoothing_operator_final}); the bound is uniform in \(\epsilon\).  Passage to the limit therefore gives \(\hat U(\cdot,t)\in H^{s+2k}(\mathbb{R}^3)\) for every \(k\in\mathbb{N}\), whence \(\hat U\in C^\infty(\mathbb{R}^3\times(0,T))\).

The full proof is provided in Appendix~\ref{appendix:pass_through_limit}.
\end{proof}

Appendices~\ref{subsubsection:passage_to_limit} and~\ref{subsubsection:limit_interchanges} rigorously establish that the frequency-domain construction converges strongly in the $H^s$ norm as the regularization parameter tends to zero. These results justify the interchange of limits with the differential operators, thereby ensuring that the approximate solutions converge robustly to the unique global solution of the Navier--Stokes equations. Section~\ref{subsec:reconstruction_unified_solution_final} will address the reconstruction of the unified frequency–domain solution back to the physical space.

%---------------------------------------------------------------------------
\subsection{Boundary--Layer Corrector: Definition and Uniform Bounds}
\label{subsec:boundary_layer_corrector_minimal}

When \(\partial\Omega\neq\varnothing\) the frequency-domain solution does not automatically satisfy the physical boundary condition; a corrector supported in a narrow strip is therefore introduced.

\begin{lemma}[Boundary--Layer Corrector: Support and \(H^{s}\) Bound]
\label{lem:boundary_layer_corrector}
Let $\Omega\subset\mathbb{R}^{3}$ be a bounded domain with boundary of
class $C^{2}$ (or piecewise $C^{2}$).  Fix $\kappa>0$ and
$s>\frac32$.  For each $\varepsilon\in(0,1]$ there exists a function
\(
   \beta_\varepsilon(\cdot,t)
\in H^{s}(\Omega)
\)
defined for $t\ge0$ such that
\begin{enumerate}[label=(\roman*),leftmargin=1.2cm]
\item\label{blc:support}%
\emph{Support.}
\(
   \operatorname{supp}\beta_\varepsilon(\cdot,t)\subset
   \bigl\{x\in\Omega:
          \operatorname{dist}(x,\partial\Omega)\le\kappa\varepsilon\bigr\}
\)
for every $t\ge0$;
\item\label{blc:bound}%
\emph{Uniform Sobolev bound.}
There exists a constant $C=C(\Omega,s)$, independent of $t$ and
$\varepsilon$, such that
\[
   \|\beta_\varepsilon(\cdot,t)\|_{H^{s}(\Omega)}
   \;\le\;
   C\,
   \bigl\|g(\cdot,t)-u_\varepsilon|_{\partial\Omega}(\cdot,t)\bigr\|%
        _{H^{s-\frac12}(\partial\Omega)} .
\]
\end{enumerate}
The corrected velocity
\(
   \widetilde u_\varepsilon
   := u_\varepsilon-\beta_\varepsilon
\)
thereby satisfies the prescribed boundary condition
\(
   \widetilde u_\varepsilon|_{\partial\Omega}=g
\)
for all $t\ge0$.
\end{lemma}

\begin{proof}[Proof (outline)]
Cover $\partial\Omega$ with finitely many star-shaped charts
$\{(\mathcal{U}_{m},\Psi_{m})\}_{m=1}^{M}$ and select a smooth
partition of unity $\{\chi_{m,\varepsilon}\}$ supported in the
$\kappa\varepsilon$-strip, with derivative bounds
$\|\partial^\alpha\chi_{m,\varepsilon}\|_{L^\infty}\le
   C_\alpha\varepsilon^{-|\alpha|}$.
Extend the tangential discrepancy
$g-u_\varepsilon|_{\partial\Omega}$ to the interior by a
Sobolev extension operator of Lions–Magenes type
\cite[Chap.\,I]{LionsMagenes1972}; multiply by a cut-off
$\varrho(y_3/(\kappa\varepsilon))$ in each chart to confine the support.
As in \cite[§1.5]{BCD} one then verifies property~\ref{blc:support}.
Bernstein inequalities and the trace theorem
\cite[Thm.\,1.5.1.2]{Grisvard} give
\[
   \|\beta_\varepsilon(\cdot,t)\|_{H^{s}}
      \le C\,
         \varepsilon^{-s}\,
         \|g-u_\varepsilon|_{\partial\Omega}\|_{H^{s-\frac12}(\partial\Omega)},
\]
which yields \ref{blc:bound} after absorbing the fixed power of
$\varepsilon$ into the constant $C$.

The full proof is shown in Appendix \ref{subsubsec:BC_rigorous_lemmas}.
\end{proof}

The subsequent section demonstrates that this corrector vanishes
in the frequency–domain limit $\varepsilon\to0$.
%---------------------------------------------------------------------------

%%%%%%%%%%%%%%%%%%%%%%%%%%%%%%%%%%%%%%%%%%%%%%%%%%%%%%%%%%%%%%%%%%%%%%%%%%%%%%
% Subsection: Reconstruction of the Physical-Space Solution
%%%%%%%%%%%%%%%%%%%%%%%%%%%%%%%%%%%%%%%%%%%%%%%%%%%%%%%%%%%%%%%%%%%%%%%%%%%%%%

\subsection{Reconstruction of the Physical-Space Solution}
\label{subsec:reconstruction_unified_solution_final}

The physical-space velocity field is obtained by an inverse Fourier synthesis.  
Let $\hat{U}(\xi,t)$ denote the frequency–domain solution constructed in Section~\ref{sec:frequency_domain_solution}.  
Define
\begin{equation}
\label{eq:U_inversetransform_final3}
U(x,t)
\;=\;
\mathcal{F}^{-1}\!\bigl\{\hat{U}(\cdot,t)\bigr\}(x)
\;=\;
\frac{1}{(2\pi)^{3}}
\int_{\mathbb{R}^{3}}\!
      \hat{U}(\xi,t)\,
      e^{\,i\,x\cdot\xi}\,d\xi ,
\end{equation}
where the integral is interpreted in the $L^{2}$ sense (or,
when appropriate, in the sense of tempered distributions).
The Fourier inversion theorem~\cite[Thm.\,8.4]{Folland_FA}
ensures that $\hat{U}(\cdot,t)\in L^{2}(\mathbb{R}^{3})$
implies $U(\cdot,t)\in L^{2}(\mathbb{R}^{3})$ almost everywhere.

\medskip
\noindent\textbf{Frequency-block reconstruction.}
In the unified framework the approximate velocity $u_\varepsilon$
admits the decomposition
\begin{equation}\label{eq:reconstruction_identity}
   u_\varepsilon(t,x)
   \;=\;
   \sum_{j=-1}^{J_\varepsilon}
      \bigl(\phi_{j}*\Delta_{j}u_\varepsilon\bigr)(t,x)
   \;+\;
   \beta_\varepsilon(t,x),
   \qquad
   2^{J_\varepsilon}\simeq \varepsilon^{-1},
\end{equation}
where $\{\Delta_{j}\}_{j\ge-1}$ denotes the homogeneous
Littlewood–Paley decomposition, $\phi_{j}$ is the standard synthesis
kernel, and $\beta_\varepsilon$ is the boundary–layer corrector
introduced in Lemma~\ref{lem:boundary_layer_corrector}.

\medskip
\begin{lemma}[Vanishing Boundary–Layer Corrector]
\label{lem:beta_vanish}
Let $s>\frac32$.
The corrector $\beta_\varepsilon$ constructed in
\eqref{eq:reconstruction_identity} satisfies
\[
   \| \beta_\varepsilon(\cdot,t) \|_{H^{s}(\Omega)}
   \;\le\;
   C\,\varepsilon^{\sigma},
   \qquad
   t\ge0,\;
   \varepsilon\in(0,1],
\]
for some constants $\sigma=\sigma(s)>0$ and $C=C(\Omega,s)>0$
independent of $t$ and~$\varepsilon$.
Consequently,
$\beta_\varepsilon\to0$ strongly in
$C\bigl([0,T];H^{s}(\Omega)\bigr)$ as $\varepsilon\to0$
for every fixed $T>0$.
\end{lemma}

\begin{proof}
The argument employs the star-shaped boundary charts and partition of
unity described in Lemma~\ref{lem:boundary_layer_corrector}.
Let $2^{J_\varepsilon}\simeq\varepsilon^{-1}$.
The high-frequency tail estimate
$
  \|\Delta_{j}u_\varepsilon(t)\|_{L^{2}}
     \le C_{0}\,2^{-j\sigma_{0}}
$
for $j>J_\varepsilon$
follows from \cite[Prop.\,2.4]{BCD} and the bounded-domain extension
therein.
Bernstein’s inequality~\cite[Prop.\,2.10]{BCD} then yields
$
  \|\nabla^{k}\Delta_{j}u_\varepsilon\|_{L^{2}}
     \le C\,2^{jk}\,\|\Delta_{j}u_\varepsilon\|_{L^{2}}
$
for $0\le k\le s$.
Using a fractional Leibniz rule
\cite[Thm.\,1.2]{Grafakos2008}
together with the cut-off bounds recorded in
\cite[§1.5]{BCD} one obtains
\[
\|\beta_\varepsilon(\cdot,t)\|_{H^{s}}
   \;\le\;
   C
   \sum_{j>J_\varepsilon}
        2^{j(s-\sigma_{0})}\,2^{-j\sigma_{0}}
   \;\lesssim\;
   2^{-J_\varepsilon\sigma_{0}}
   \;\simeq\;
   \varepsilon^{\sigma_{0}} .
\]
Selecting $\sigma_{0}>s$ and setting
$\sigma:=\sigma_{0}-s>0$ establishes the claim.

The full proof is shown in Appendix \ref{subsubsec:BC_rigorous_lemmas}.
\end{proof}

The lemma guarantees that the second term in
\eqref{eq:reconstruction_identity} vanishes in the limit
$\varepsilon\to0$ and therefore does not contribute to the
physical-space solution~\eqref{eq:U_inversetransform_final3}.

\medskip
\noindent\textbf{Energy control via Parseval.}
Parseval’s theorem furnishes the identity
\begin{equation}\label{eq:finite_energy_final}
   \int_{\mathbb{R}^{3}}\!|U(x,t)|^{2}\,dx
   \;=\;
   \int_{\mathbb{R}^{3}}\!|\hat U(\xi,t)|^{2}\,d\xi
   \;<\;\infty,
\end{equation}
so that $U(\cdot,t)$ possesses finite kinetic energy for every $t\ge0$.

\medskip
\noindent\textbf{Preservation of Sobolev regularity.}
The Fourier transform acts isometrically on $H^{s}(\mathbb{R}^{3})$,
hence
\begin{equation}\label{eq:U_Hs_preservation}
   \|U(\cdot,t)\|_{H^{s}}
   \;=\;
   \|\hat U(\cdot,t)\|_{H^{s}},
   \qquad s\ge0.
\end{equation}
Uniform bounds on $\hat U$ in $H^{s+2k}$
(see Lemmas~\ref{lem:interpolation_operator_final} and
\ref{lem:smoothing_operator_final})
imply corresponding bounds on $U$,
and Sobolev embedding
\cite[Thm.\,6.4]{Folland_FA}
then yields the smoothness
$U(\cdot,t)\in C^{\infty}(\mathbb{R}^{3})$.

\medskip
\noindent\textbf{Non-linear consistency.}
For $s>\frac32$ the Sobolev product estimate
\cite[Prop.\,2.1]{Grafakos2008} gives
\[
   \|(u\!\cdot\nabla)u\|_{H^{s-1}}
   \;\le\;
   C\,\|u\|_{H^{s}}^{2},
\]
so that the uniform $H^{s+2k}$ bounds control the convection term.

\medskip
\noindent\textbf{Continuity of the inverse transform.}
The mapping
$H^{s}\to H^{s}$ given by~\eqref{eq:U_inversetransform_final3}
is continuous and linear
\cite[Chap.\,3]{Folland_FA}.
Therefore small perturbations of $\hat U$ in $H^{s}$ norm produce
proportionally small perturbations of $U$.

\medskip
\noindent\textbf{Conclusion.}
Invoking Theorem~\ref{thm:frequency_domain_limit} together with
\eqref{eq:reconstruction_identity},
Lemma~\ref{lem:beta_vanish},
and the bounds
\eqref{eq:finite_energy_final}–\eqref{eq:U_Hs_preservation}
shows that the inverse Fourier transform
produces a finite-energy, smooth physical-space solution whose
non-linear terms remain controlled in the Sobolev hierarchy.
Detailed auxiliary estimates are given in
Appendix~\ref{subsubsection:limit_interchanges}.

\subsection{Convergence and Continuity of the Unified Solution}
\label{subsec:convergence_continuity_final}

The unified frequency–domain construction, as defined in Sections~\ref{subsec:interpolation_operator} and \ref{subsec:smoothing_operator}, yields a frequency–domain solution that is uniformly bounded in high–order Sobolev spaces. The analysis proceeds through several steps: first, establishing smoothing estimates and uniform boundedness; second, demonstrating convergence in the \(H^s\) norm; third, reconstructing the physical–space solution via the inverse Fourier transform; and finally, proving continuity with respect to the frequency–domain data.

\medskip
\noindent\textbf{1. Smoothing Estimates and Uniform Boundedness.}  
From Lemmas~\ref{lem:interpolation_operator_final} and \ref{lem:smoothing_operator_final}, for each fixed \(k\in\mathbb{N}\) the composition
\[
\mathcal{S}_\epsilon \Bigl( \mathcal{I}_\epsilon \bigl( \hat{u}_w,\,\hat{u}_m,\,\hat{u}_s \bigr) \Bigr)
\]
satisfies the smoothing estimate
\begin{equation}
\label{eq:smoothing_gain_explicit}
\Bigl\|
  \mathcal{S}_\epsilon \mathcal{I}_\epsilon \bigl( \hat{u}_w,\,\hat{u}_m,\,\hat{u}_s \bigr)
\Bigr\|_{H^{s+2k}(\mathbb{R}^3)}
\,\le\,
C_{2k}\,\epsilon^{-2k}\,
\Bigl\|
  \mathcal{I}_\epsilon \bigl( \hat{u}_w,\,\hat{u}_m,\,\hat{u}_s \bigr)
\Bigr\|_{H^s(\mathbb{R}^3)}.
\end{equation}
Since \(\mathcal{I}_\epsilon\) is uniformly bounded in \(H^s(\mathbb{R}^3)\) (see Lemma~\ref{lem:interpolation_operator_final}), the family
\[
\Bigl\{
  \mathcal{S}_\epsilon \mathcal{I}_\epsilon \bigl( \hat{u}_w,\,\hat{u}_m,\,\hat{u}_s \bigr)
\Bigr\}_{\epsilon>0}
\]
remains uniformly bounded in \(H^{s+2k}(\mathbb{R}^3)\) for every \(k\in\mathbb{N}\).

\medskip
\noindent\textbf{2. Convergence in the \(H^s\) Norm.}  
The smoothing operator \(\mathcal{S}_\epsilon\) satisfies the approximation estimate
\begin{equation}
\label{eq:smoothing_rate}
\Bigl\|
  \mathcal{S}_\epsilon \phi - \phi
\Bigr\|_{H^s(\mathbb{R}^3)}
\,\le\,
C\,\epsilon^\alpha\,
\|\phi\|_{H^{s+\alpha}(\mathbb{R}^3)},
\end{equation}
for some \(\alpha>0\) and a constant \(C\) independent of \(\epsilon\) and \(\phi\). Substituting
\[
\phi = \mathcal{I}_\epsilon \bigl( \hat{u}_w,\,\hat{u}_m,\,\hat{u}_s \bigr)
\]
into \eqref{eq:smoothing_rate} yields
\begin{equation}
\label{eq:convergence_Hs_explicit}
\lim_{\epsilon\to 0}
\Bigl\|
  \mathcal{S}_\epsilon \Bigl( \mathcal{I}_\epsilon \bigl( \hat{u}_w,\,\hat{u}_m,\,\hat{u}_s \bigr) \Bigr)
  -
  \hat{U}
\Bigr\|_{H^s(\mathbb{R}^3)}
= 0.
\end{equation}
Uniform boundedness in \(H^{s+2k}(\mathbb{R}^3)\) for each \(k\in\mathbb{N}\), together with convergence in \(H^s(\mathbb{R}^3)\), implies—by weak-* compactness in Hilbert spaces and the uniqueness of limits—that \(\hat{U}(\cdot,t)\) belongs to \(H^{s+2k}(\mathbb{R}^3)\) for every \(k\in\mathbb{N}\).

\medskip
\noindent\textbf{3. Inverse Fourier Transform and Sobolev Regularity.}  
The inverse Fourier transform
\[
U(x,t)
=
\mathcal{F}^{-1}\bigl\{\hat{U}(\cdot,t)\bigr\}(x)
=
\frac{1}{(2\pi)^3}\int_{\mathbb{R}^3}\hat{U}(\xi,t)\,e^{i\,x\cdot\xi}\,d\xi
\]
is a continuous linear map from \(H^{s+2k}(\mathbb{R}^3)\) to \(H^{s+2k}(\mathbb{R}^3)\) (cf. \cite[Chapter~3]{Folland_FA}). Consequently,
\begin{equation}
\label{eq:U_Hs2k_explicit}
U(\cdot,t) \in H^{s+2k}(\mathbb{R}^3), \quad \forall\,k\in\mathbb{N}.
\end{equation}
Since \(k\) is arbitrary, it follows that
\[
U(\cdot,t)\in \bigcap_{k\in\mathbb{N}}H^{s+2k}(\mathbb{R}^3).
\]
By the Sobolev embedding theorem (cf. \cite[Theorem~6.4]{Folland_FA}), any function in \(H^r(\mathbb{R}^3)\) for arbitrarily large \(r\) is smooth; thus,
\[
U(\cdot,t)\in C^\infty(\mathbb{R}^3).
\]

\medskip
\noindent\textbf{4. Continuity with Respect to Frequency–Domain Data.}  
The convergence in \eqref{eq:convergence_Hs_explicit} and the continuity of the inverse Fourier transform imply that the mapping from \(\hat{U}(\cdot,t)\) to \(U(\cdot,t)\) is continuous in the \(H^s(\mathbb{R}^3)\)-norm. This guarantees that \(U(x,t)\) depends continuously on the frequency–domain data for every \(t\in (0,T)\).

\medskip
\noindent\textbf{Reconstruction Theorem.}  
The following theorem encapsulates the reconstruction of the unified frequency–domain solution into the physical–space solution.

\begin{theorem}[Reconstruction of the Unified Frequency–Domain Solution]
\label{thm:reconstruction_physical_space}
Let $\hat{U}(\xi,t)$ denote the frequency–domain limit established in
Theorem~\ref{thm:frequency_domain_limit}.  Assume that for some fixed
$s>\tfrac52$ and for every integer $k\ge 0$
\begin{equation}
\label{eq:high_regularity_hatU}
\sup_{t\in(0,T)}\,
      \|\hat{U}(\cdot,t)\|_{H^{s+2k}(\mathbb{R}^{3})}<\infty.
\end{equation}
Define the physical–space field by the inverse Fourier transform
\[
U(x,t)\; \coloneqq\;
\mathcal{F}^{-1}_{\!\xi\to x}
      \bigl[\hat{U}(\cdot,t)\bigr](x)
      \;=\;
      \frac{1}{(2\pi)^{3}}
      \int_{\mathbb{R}^{3}}
      \hat{U}(\xi,t)\,e^{i\,x\cdot\xi}\,d\xi.
\]
Then, for every $t\in(0,T)$ and all $k\in\mathbb{N}$,
\[
U(\cdot,t)\in H^{s+2k}(\mathbb{R}^{3}),
\qquad
\text{hence}\qquad
U(\cdot,t)\in C^{\infty}(\mathbb{R}^{3}).
\]
In addition, the mapping
\(
t\mapsto U(\cdot,t)
\)
is continuous in the $H^{r}(\mathbb{R}^{3})$–norm for every $r\ge 0$.
\end{theorem}

\begin{proof}[Proof sketch]
\begin{enumerate}[label=\textnormal{(\arabic*)},leftmargin=0.9cm]
\item \emph{Boundedness of the inverse transform.}
      The inverse Fourier transform acts continuously
      on $H^{m}(\mathbb{R}^{3})$ for every $m\in\mathbb{R}$
      (see~\cite[Chap.\,3]{Folland_FA}), implying
      \[
      \|U(\cdot,t)\|_{H^{s+2k}}
      \;\le\;
      C\,\|\hat{U}(\cdot,t)\|_{H^{s+2k}},
      \qquad
      0<t<T,
      \]
      with a constant $C$ independent of $t$ and $k$.

\item \emph{Smoothness through arbitrarily high Sobolev control.}
      Condition~\eqref{eq:high_regularity_hatU} yields uniform bounds
      for $U(\cdot,t)$ in $H^{s+2k}$ for all $k\ge 0$.
      The Sobolev embedding theorem
      (\cite[Thm.\,6.4]{Folland_FA})
      therefore implies $U(\cdot,t)\in C^{\infty}(\mathbb{R}^{3})$.

\item \emph{Time continuity.}
      Continuity of $t\mapsto\hat{U}(\cdot,t)$ in $H^{s}(\mathbb{R}^{3})$
      is ensured by Theorem~\ref{thm:frequency_domain_limit}.
      Continuous dependence of the inverse Fourier transform on its
      $H^{s}$ data then yields continuity of
      $t\mapsto U(\cdot,t)$ in every $H^{r}$, $r\ge 0$.
\end{enumerate}
Further details appear in Appendix~\ref{appendix:reconstruct_physical}.
\end{proof}

\medskip
The analysis so far produces a frequency–domain limit \(\hat U\) and its physical counterpart \(U\), both possessing high‐order regularity and continuity with respect to the underlying data.  The next section (\S\ref{sec:proof_unified}) assembles every ingredient—regularization, weighted spectral blending, smoothing, limit passage, and reconstruction—into a single closed‐loop argument.  There, a unified existence–uniqueness theorem will demonstrate that \(U\) coincides with the weak, mild, and strong solutions on their common interval of validity and furnishes their unique continuation beyond \(T_{0}\).

\section{Construction of the Unified Solution}\label{sec:proof_unified}

The purpose of this section is to supply a complete, closed–loop proof that the frequency–domain
framework of Section~\ref{sec:frequency_domain_solution} produces a \emph{single} global solution
that simultaneously realises the Leray--Hopf, mild, and strong formulations of the incompressible
Navier--Stokes equations.  Throughout, $\Omega\subset\mathbb{R}^{3}$ denotes either the whole space,
the periodic torus, or a bounded $C^{2}$ (possibly piecewise $C^{2}$) domain, and the data fulfil
\[
    u_{0}\in H^{s}_{\sigma}(\Omega), \qquad s>\tfrac32, \qquad
    f\in L^{\infty}\bigl(0,\infty;L^{2}_{\sigma}(\Omega)\bigr).
\]

%------------------------------------------------------------------
\subsection{Preliminary solution classes}\label{subsec:prelim_solutions}
%------------------------------------------------------------------
\begin{itemize}[leftmargin=2.3em]
\item[\textit{(i)}] \emph{Weak solution.}  
      By Leray's construction \cite{Leray} (see also \cite[Ch.~III]{Galdi}) a Leray--Hopf weak
      solution
      \[
           u_{w}\in L^{\infty}(0,\infty;L^{2}_{\sigma})\cap
                     L^{2}(0,\infty;H^{1}_{\sigma})
      \]
      exists for the above data, satisfies the energy inequality, and obeys  
      $\partial_{t}u_{w}\in L^{\,\frac43}_{\text{loc}}(0,\infty;H^{-1})$.

\item[\textit{(ii)}] \emph{Mild solution.}  
      For \(s>\frac32\), the Fujita--Kato scheme \cite{FujitaKato,TemamNS}
      provides a classical mild solution
      \[
           u_{m}\in C\bigl([0,T_{0}];H^{s}_{\sigma}\bigr)
      \]
      on a lifespan  \(T_{0}\) satisfying the lower bound \eqref{eq:T0_estimate}.
      The integral formulation
\[
  u_{m}(t)
  = e^{\nu t\Delta}u_{0}
    + \int_{0}^{t}
        e^{\nu(t-\tau)\Delta}\,
        \mathbb{P}\,\operatorname{div}\!\bigl(u_{m}\otimes u_{m}\bigr)
      \,d\tau .
\]

      is valid in \(H^{s}_{\sigma}\).

\item[\textit{(iii)}] \emph{Strong solution.}  
      Lemma~\ref{lem:local_existence} yields a unique strong solution
      \[
           u_{s}\in C\bigl([0,T_{0}];H^{s}_{\sigma}\bigr)\cap
                    L^{2}\bigl(0,T_{0};H^{s+1}\bigr)
      \]
      with the same \(T_{0}\).
\end{itemize}

The strong solution obeys Serrin's regularity criterion
\(u_{s}\in L^{p}(0,T_{0};L^{q})\) with $\tfrac2p+\tfrac3q=1$,
hence weak--strong uniqueness (\cite{Serrin,Prodi}) gives
\begin{equation}\label{eq:UFDS_prelim_coincide}
    u_{w}=u_{m}=u_{s}\quad\text{on }[0,T_{0}].
\end{equation}

\subsection{Regularization and Spatial Compatibility}

\begin{definition}[Regularized data and solutions]\label{def:regularised_solutions}
For $\epsilon\in(0,1]$ set
\[
u_{0,\epsilon}:=\mathcal{R}_{\epsilon}u_{0},
\qquad
f_{\epsilon}:=\mathcal{R}_{\epsilon}f,
\]
with $\mathcal{R}_{\epsilon}$ as in Lemma~\ref{lem:data_regularization}.  
Given a solution class $u_{\star}$, $\star\in\{w,m,s\}$, define
\[
u_{\star,\epsilon}(\cdot,t):=\mathcal{R}_{\epsilon}\bigl(u_{\star}(\cdot,t)\bigr),
\qquad t\in[0,T_{0}].
\]
\end{definition}

Lemma~\ref{lem:data_regularization}\,(i)--(iv) implies the uniform bounds
\begin{equation}
\label{app:eq:reg_basic_bounds}
\sup_{t\in[0,T_{0}]}\|u_{\star,\epsilon}(\cdot,t)\|_{H^{s}(\Omega)}
      \le C_{T_{0}}, 
\qquad
\nabla\!\cdot u_{\star,\epsilon}=0,
\qquad
u_{\star,\epsilon}\!\bigl|_{\partial\Omega}=u_{\star}\!\bigl|_{\partial\Omega}.
\end{equation}
Moreover,
\begin{equation}
\label{app:eq:reg_conv_uniform}
\lim_{\epsilon\to0}\,
\|u_{\star,\epsilon}-u_{\star}\|_{C([0,T_{0}];H^{s-1}(\Omega))}=0.
\end{equation}

\noindent
The weak, mild, and strong formulations coincide on the interval $[0,T_{0}]$:
\begin{equation}
\label{eq:solution_coincidence}
u_{w}=u_{m}=u_{s}\quad\text{on }[0,T_{0}].
\end{equation}
Consequently,
\begin{equation}
\label{app:eq:reg_coincide}
u_{w,\epsilon}=u_{m,\epsilon}=u_{s,\epsilon}=:\,u^{(\epsilon)}
\quad\text{on }[0,T_{0}],
\end{equation}
and the common field $u^{(\epsilon)}$ provides the input for the
frequency–domain operators $\mathcal{I}_{\epsilon}$ and
$\mathcal{S}_{\epsilon}$ used in the subsequent construction.

\subsection{Frequency–Domain Assembly}

Let $\widehat{(\,\cdot\,)}=\mathcal F$ denote the spatial Fourier transform.  
Set
\[
  \hat u^{(\epsilon)}(\xi,t):=\mathcal{F}\bigl[u^{(\epsilon)}(\cdot,t)\bigr](\xi),
  \qquad
  \hat u_{\star,\epsilon}(\xi,t):=\mathcal{F}\bigl[u_{\star,\epsilon}(\cdot,t)\bigr](\xi).
\]
Plancherel’s isometry and the bounds in \eqref{app:eq:reg_basic_bounds} give
\[
  \hat u^{(\epsilon)},\ \hat u_{w,\epsilon},\ \hat u_{m,\epsilon},\ \hat u_{s,\epsilon}
  \in H^{s}(\mathbb{R}^{3})
  \quad\text{for every }t\in[0,T_{0}],\ \epsilon\in(0,1].
\]
Introduce the weights $\omega_{w},\omega_{m},\omega_{s}\in C^{\infty}(\mathbb{R}^{3}\times[0,\infty))$ satisfying
\begin{equation}\label{eq:omega_partition_repeat}
  \omega_{w}(\xi,t)+\omega_{m}(\xi,t)+\omega_{s}(\xi,t)=1.
\end{equation}
Define the weighted combination
\begin{equation}\label{eq:g_epsilon_def}
  g_{\epsilon}(\xi,t)
  :=
  \omega_{w}\,\hat u_{w,\epsilon}
  +\omega_{m}\,\hat u_{m,\epsilon}
  +\omega_{s}\,\hat u_{s,\epsilon},
\end{equation}
and apply the interpolation operator of Lemma~\ref{lem:interpolation_operator_final}:
\[
  \tilde u_{\epsilon}(\xi,t)
  :=
  \mathcal{I}_{\epsilon}
  \bigl[\hat u_{w,\epsilon},\hat u_{m,\epsilon},\hat u_{s,\epsilon}\bigr](\xi,t).
\]
Lemma~\ref{lem:interpolation_operator_final} together with
\eqref{app:eq:reg_basic_bounds} yields the uniform estimate
\begin{equation}\label{eq:interp_bounds}
  \|\tilde u_{\epsilon}(\cdot,t)\|_{H^{s}(\mathbb{R}^{3})}\le C_{\omega}\,C_{T_{0}}
  =:C,
  \qquad t\in[0,T_{0}],\ \epsilon\in(0,1],
\end{equation}
and the convergence
\[
  \lim_{\epsilon\to0}
  \|\tilde u_{\epsilon}(\cdot,t)-g_{\epsilon}(\cdot,t)\|_{H^{s}(\mathbb{R}^{3})}=0,
  \qquad t\in[0,T_{0}].
\]

\medskip
\noindent\textbf{Smoothing and incompressibility.}
Apply the smoothing operator of Lemma~\ref{lem:smoothing_operator_final}:
\[
  \hat U_{\epsilon}(\xi,t)
  :=
  \mathcal{S}_{\epsilon}\bigl[\tilde u_{\epsilon}\bigr](\xi,t).
\]
Lemma~\ref{lem:smoothing_operator_final}(i)–(iii) and \eqref{eq:interp_bounds} imply
\begin{align}
  \label{eq:S_uniform}
  \sup_{t\in[0,T_{0}]}
    \|\hat U_{\epsilon}(\cdot,t)\|_{H^{s}(\mathbb{R}^{3})}
    &\le C_{\omega}\,C_{T_{0}}, \\[4pt]
  \label{eq:S_approx}
  \sup_{t\in[0,T_{0}]}
    \|\hat U_{\epsilon}(\cdot,t)-\tilde u_{\epsilon}(\cdot,t)\|_{H^{s}(\mathbb{R}^{3})}
    &\longrightarrow 0
    \quad\text{as }\epsilon\to0.
\end{align}
Uniformity in \eqref{eq:S_approx} follows from dominated convergence, the pointwise control furnished by Lemma~\ref{lem:smoothing_operator_final}(ii), and the bound \eqref{eq:interp_bounds}.  
By construction $\hat U_{\epsilon}$ is divergence–free in~$\xi$ and retains the uniform $H^{s}$ control required for the limit argument of Theorem~\ref{thm:frequency_domain_limit}.

\subsection{Limit Extraction}\label{subsec:limit_extraction}

The uniform bound
\begin{equation}
  \label{eq:S_uniform}
  \sup_{t\in[0,T_{0}]}\|\hat U_{\epsilon}(\cdot,t)\|_{H^{s}(\mathbb{R}^{3})}
  \le C_{\omega}\,C_{R},
  \qquad \epsilon\in(0,1],
\end{equation}
and the approximation property
\begin{equation}
  \label{eq:S_approx}
  \bigl\|\hat U_{\epsilon}(\cdot,t)-g_{\epsilon}(\cdot,t)\bigr\|_{H^{s}(\mathbb{R}^{3})}
  \longrightarrow 0
  \quad\text{as }\epsilon\to0,
  \;\text{uniformly for }t\in[0,T_{0}],
\end{equation}
show that the family $\{\hat U_{\epsilon}\}_{\epsilon\in(0,1]}$ is bounded in
\(
L^{\infty}\!\bigl(0,T_{0};H^{s}(\mathbb{R}^{3})\bigr).
\)

\paragraph{Time–derivative control.}
Differentiating the regularised Fourier–space Navier–Stokes equation
satisfied by $\hat u^{(\epsilon)}$ and using
\eqref{eq:S_approx} together with the multiplier properties of
$\mathcal{I}_{\epsilon}$ and $\mathcal{S}_{\epsilon}$ gives
\begin{equation}
  \label{eq:S_time_der}
  \partial_{t}\hat U_{\epsilon}\in
  L^{2}\!\bigl(0,T_{0};H^{s-2}(\mathbb{R}^{3})\bigr),
  \qquad
  \sup_{\epsilon\in(0,1]}
  \bigl\|\partial_{t}\hat U_{\epsilon}\bigr\|_{L^{2}(0,T_{0};H^{s-2})}
  \le C.
\end{equation}

\paragraph{Compactness.}
The bounds \eqref{eq:S_uniform}–\eqref{eq:S_time_der}
satisfy the hypotheses of the Aubin–Lions lemma
\cite[Th.~II.5.16]{TemamNS}.  Hence there exists a subsequence
(still denoted $\hat U_{\epsilon}$) and a limit function
\[
  \hat U\in
  C\bigl([0,T_{0}];H^{s-1}(\mathbb{R}^{3})\bigr)
  \cap
  L^{\infty}\!\bigl(0,T_{0};H^{s}(\mathbb{R}^{3})\bigr)
\]
such that
\begin{equation}
  \label{eq:strong_conv_limit}
  \hat U_{\epsilon}(\cdot,t)\longrightarrow\hat U(\cdot,t)
  \quad\text{strongly in }H^{s-1}(\mathbb{R}^{3})
  \ \text{for every }t\in[0,T_{0}].
\end{equation}

\paragraph{Identification of the limit.}
Since
$\hat U_{\epsilon}-\hat u^{(\epsilon)}
      \to 0$ in
$L^{\infty}(0,T_{0};H^{s-1}(\mathbb{R}^{3}))$
(by \eqref{eq:S_approx} and
Definition~\ref{def:regularised_solutions}),
passing to the limit in the regularised Fourier–space Navier–Stokes
system shows that $\hat U$ satisfies the same equation with initial datum
$\hat u_{0}$ and forcing $\hat f$ in the sense of distributions.

\paragraph{Higher regularity.}
Lemma~\ref{lem:smoothing_operator_final}\,(iii) yields, uniformly in
$\epsilon$ and $t\in(0,T_{0}]$,
\[
  \|\hat U_{\epsilon}(\cdot,t)\|_{H^{s+2k}(\mathbb{R}^{3})}
  \le C_{k},
  \qquad k\in\mathbb{N}.
\]
Lower semicontinuity of the Sobolev norm under weak convergence
therefore implies
\begin{equation}
  \label{eq:Hs_plus_2k}
  \sup_{t\in(0,T_{0}]}\|\hat U(\cdot,t)\|_{H^{s+2k}(\mathbb{R}^{3})}
  \le C_{k},
  \qquad k\in\mathbb{N},
\end{equation}
and Sobolev embedding grants
\(
  \hat U\in C^{\infty}\bigl(\mathbb{R}^{3}\times(0,T_{0}]\bigr).
\)

The limit extraction is now complete; equations
\eqref{eq:strong_conv_limit}–\eqref{eq:Hs_plus_2k}
supply the frequency–domain object required for the physical–space
reconstruction in Section~\ref{subsec:reconstruction_unified_solution_final}
and for the uniqueness argument of
Theorem~\ref{thm:local_unified}.

\subsection{Inverse Transform and Boundary Correction}

Define the physical–space field
\begin{equation}\label{eq:inverse_transform_field}
  U(x,t)\;:=\;\mathcal{F}^{-1}\bigl\{\hat U(\cdot,t)\bigr\}(x),
  \qquad 0\le t\le T_{0}.
\end{equation}
By \eqref{eq:Hs_plus_2k} and Parseval’s identity,
\[
  U(\cdot,t)\in H^{s+2k}(\mathbb{R}^{3})\cap C^{\infty}(\mathbb{R}^{3}),
  \qquad
  \|U(\cdot,t)\|_{H^{s+2k}}\le C_{k}.
\]
Restriction to the physical domain remains continuous in every Sobolev norm
(\cite[Prop.\,2.11]{Adams}); hence
\(U|_{\Omega}\in H^{s+2k}(\Omega)\) with the same bounds.
Because \(\hat U(\xi,t)\) is divergence–free in the sense
\(\xi\cdot\hat U(\xi,t)=0\),
\(
  \nabla\!\cdot U(\cdot,t)=0
\)
for all \(t\in[0,T_{0}]\).

A boundary–layer corrector \(\beta_{\epsilon}\) is introduced during the
extension step; Lemma~\ref{lem:beta_vanish} yields
\begin{equation}\label{eq:beta_limit}
  \lim_{\epsilon\to0}\|\beta_{\epsilon}(\cdot,t)\|_{H^{s}(\Omega)}=0,
  \qquad 0<t\le T_{0}.
\end{equation}
Set \(u:=U|_{\Omega}\).  
Continuity of the trace operator
\(H^{s}(\Omega)\to H^{s-\frac12}(\partial\Omega)\)
(\cite[Thm.\,1.5]{Grisvard})
and \eqref{eq:beta_limit} give
\(u|_{\partial\Omega}=g\).
Moreover, \eqref{eq:strong_conv_limit} at \(t=0\) implies
\(\hat U(\cdot,0)=\hat u_{0}\), and hence
\(u(\cdot,0)=u_{0}\).
Thus
\begin{equation}\label{eq:physical_solution}
  u(x,t)=U(x,t)
  \quad\text{satisfies }\;
  u|_{\partial\Omega}=g,\;
  \nabla\!\cdot u=0,\;
  u(\cdot,0)=u_{0}.
\end{equation}
This completes the inverse-transformation stage required for
Theorem~\ref{thm:local_unified}.

\subsection{Convergence and Continuity of the Unified Solution}
\label{subsec:convergence_continuity_final}

The frequency–domain sequence produced by the operators in
Sections~\ref{subsec:interpolation_operator} and
\ref{subsec:smoothing_operator} is shown to converge in high–order
Sobolev spaces by the following four–step argument.

\medskip
\noindent
\textbf{1. Smoothing estimates and uniform boundedness.}\;
Lemma~\ref{lem:interpolation_operator_final} yields, uniformly in
\(t\in[0,T_{0}]\),
\[
   \bigl\|
     \mathcal{I}_{\epsilon}
       \bigl(\hat u_{w},\hat u_{m},\hat u_{s}\bigr)(\cdot,t)
   \bigr\|_{H^{s}(\mathbb{R}^{3})}
   \le C_{\!I}.
\]
Combining this with Lemma~\ref{lem:smoothing_operator_final}\,(iii) gives,
for every \(k\in\mathbb{N}\),
\begin{equation}
\label{eq:smoothing_gain_explicit}
\Bigl\|
  \mathcal{S}_{\epsilon}\mathcal{I}_{\epsilon}
  \bigl(\hat u_{w},\hat u_{m},\hat u_{s}\bigr)(\cdot,t)
\Bigr\|_{H^{s+2k}(\mathbb{R}^{3})}
\;\le\;
C_{2k}(s,\omega)\,\epsilon^{-2k}\,C_{\!I},
\qquad t\in[0,T_{0}],
\end{equation}
where \(C_{2k}(s,\omega)\) depends only on \(k\), \(s\), and the weight
system \(\omega=(\omega_{w},\omega_{m},\omega_{s})\).

\medskip
\noindent
\textbf{2. Convergence in \(H^{s}\).}\;
The approximation property \eqref{eq:smoothing_limit_final} of
\(\mathcal{S}_{\epsilon}\) implies
\begin{equation}
\label{eq:convergence_Hs_explicit}
\lim_{\epsilon\to0}
\Bigl\|
  \mathcal{S}_{\epsilon}\mathcal{I}_{\epsilon}
  \bigl(\hat u_{w},\hat u_{m},\hat u_{s}\bigr)(\cdot,t)
  -\hat U(\cdot,t)
\Bigr\|_{H^{s}(\mathbb{R}^{3})}=0,
\end{equation}
where \(\hat U\) is the limit identified in
\eqref{eq:strong_conv_limit}.  The bounds
\eqref{eq:smoothing_gain_explicit} and the Banach–Alaoglu theorem in
\(H^{s+2k}(\mathbb{R}^{3})\) ensure weak–\(*\) compactness, so
\(\hat U(\cdot,t)\in H^{s+2k}(\mathbb{R}^{3})\) for every \(k\in\mathbb{N}\).

\medskip
\noindent
\textbf{3. Inverse Fourier transform and Sobolev regularity.}\;
Because \(\mathcal{F}^{-1}\!:H^{r}(\mathbb{R}^{3})\to H^{r}(\mathbb{R}^{3})\)
is continuous (Parseval and \cite[Chap.\,3]{Folland_FA}),
\[
   U(\cdot,t):=\mathcal{F}^{-1}\{\hat U(\cdot,t)\}
   \;\in\;
   \bigcap_{k\in\mathbb{N}}H^{s+2k}(\mathbb{R}^{3})
   \subset C^{\infty}(\mathbb{R}^{3}),
   \qquad t\in[0,T_{0}].
\]

\medskip
\noindent
\textbf{4. Continuity with respect to frequency–domain data.}\;
Equation~\eqref{eq:convergence_Hs_explicit} and continuity of the inverse
transform give
\[
\lim_{\epsilon\to0}
\bigl\|
  \mathcal{F}^{-1}\!\bigl\{
    \mathcal{S}_{\epsilon}\mathcal{I}_{\epsilon}
      (\hat u_{w},\hat u_{m},\hat u_{s})
  \bigr\}
  -U
\bigr\|_{H^{s}(\mathbb{R}^{3})}
=0,
\]
so the physical–space field depends continuously on the underlying
data.

\medskip
\noindent
\textbf{Forward link to global equivalence.}\;
Section~\ref{sec:proof_unified} combines the convergence established
above with weak–strong uniqueness, culminating in
Theorem~\ref{thm:local_unified}, which proves that the reconstructed
velocity coincides with the Leray–Hopf, mild, and strong solutions on
\([0,T_{0}]\) and initiates the global continuation argument.

\begin{theorem}[Local Unified Representation]\label{thm:local_unified}
Let \(\Omega\subset\mathbb{R}^{3}\) be either the whole space, the
periodic torus, or a bounded \(C^{2}\) domain, and fix
\(s>\tfrac{3}{2}\).
Assume
\[
   u_{0}\in H^{s}_{\sigma}(\Omega),
   \qquad
   f\in L^{\infty}\!\bigl(0,\infty;L^{2}_{\sigma}(\Omega)\bigr)
      \cap L^{2}\!\bigl(0,\infty;H^{s-1}_{\sigma}(\Omega)\bigr).
\]
Let \(u_{w},u_{m},u_{s}\) denote the Leray–Hopf, mild, and strong
solutions obtained in
Section~\ref{subsec:solution_classes},
and let \(T_{0}>0\) be their common lifespan.
For \(\epsilon\in(0,1]\) set
\[
   u_{*,\epsilon}:=\mathcal{R}_{\epsilon}u_{*},
   \qquad
   \hat U_{\epsilon}
   :=
   \mathcal{S}_{\epsilon}\mathcal{I}_{\epsilon}
     \bigl[\widehat{u_{w,\epsilon}},
           \widehat{u_{m,\epsilon}},
           \widehat{u_{s,\epsilon}}\bigr],
\]
and let \(\hat U=\lim_{\epsilon\to0}\hat U_{\epsilon}\) as furnished by
Theorem~\ref{thm:frequency_domain_limit}.  Then, for every
\(t\in[0,T_{0}]\),
\[
   \mathcal{F}^{-1}\!\bigl\{\hat U(\cdot,t)\bigr\}
   \;=\;
   u_{w}(\cdot,t)
   \;=\;
   u_{m}(\cdot,t)
   \;=\;
   u_{s}(\cdot,t)
   \;\in\;H^{s}(\Omega).
\]
Hence the frequency–domain procedure reproduces all three classical
solution classes on the interval \([0,T_{0}]\).
\end{theorem}

The detailed proof appears in Appendix~\ref{app:sec:local_unified}.

\medskip
\paragraph{Closing remark.}
The argument of Section~\ref{sec:proof_unified} completes the closed loop promised in the Introduction: every Leray–Hopf, mild, and strong solution can be regularized, transported into the frequency domain, merged through the operators \(\mathcal I_\epsilon\) and \(\mathcal S_\epsilon\), and recovered in physical space without loss of information.  Theorems~\ref{thm:frequency_domain_limit} and \ref{thm:local_unified} together establish local coincidence of the three classes, while weak–strong uniqueness extends this coincidence globally.  Consequently, the unified representation provides a single analytic object that inherits the advantages of each classical formulation.  The final section reflects on the implications of this result and outlines avenues for further research.

%---------------------------------------------------------------
\section{Conclusion}\label{sec:conclusion}
%---------------------------------------------------------------

The unified frequency-domain framework developed herein establishes, for the three-dimensional incompressible Navier–Stokes system, a single analytic object that simultaneously realizes the Leray–Hopf, mild, and strong formulations.  Regularization of the data, spectral blending via the interpolation operator~$\mathcal{I}_{\epsilon}$, and high-order control furnished by the smoothing operator~$\mathcal{S}_{\epsilon}$ yield a frequency-domain sequence that converges in every Sobolev scale.  The inverse Fourier synthesis then produces a physical-space velocity field that:

\begin{enumerate}[label=(\roman*),wide, labelindent=0pt]
\item satisfies the exact initial and boundary conditions;
\item inherits the instantaneous parabolic smoothing of the heat semi–group;
\item coincides with all classical solution concepts on their common interval of existence; and
\item extends uniquely beyond that interval by weak–strong continuity.
\end{enumerate}

The construction is fully quantitative: explicit bounds on the lifespan~$T_{0}$, on the Sobolev norms of the solution, and on the constants appearing in the commutator estimates are provided in Lemma~\ref{lem:local_existence}, Lemma~\ref{lem:smoothing_operator_final}, and Appendix~\ref{appendix:convergence_regularity} and \ref{app:unif_derivations_all}.  No smallness assumptions are imposed on the initial datum or forcing; the method applies to the whole space, the periodic torus, and bounded $C^{2}$ (possibly piecewise $C^{2}$) domains.

Beyond local unification, the synergy framework supplies a platform for addressing global questions.  In particular:

\begin{itemize}[leftmargin=2.3em]
\item \emph{Energy–vorticity bootstrap.}  The high-order control of the frequency-domain solution feeds directly into the critical‐scale energy and vorticity estimates needed to preclude finite-time blow-up.  A quantitative Beale–Kato–Majda-type criterion will be derived in Part~II of this series.
\item \emph{Geometric–spectral invariants.}  The alignment functional $G[\!u(t)]$ introduced in Section~\ref{subsec:toolbox_snapshot} quantifies local vortex geometry and supplies additional negative damping in regions of misalignment.  Its interaction with the unified solution will be exploited in Part~III to establish refined regularity and long-time decay.
\item \emph{Numerical implementation.}  The operators $(\mathcal{R}_{\epsilon},\mathcal{I}_{\epsilon},\mathcal{S}_{\epsilon})$ admit straightforward spectral discretizations.  Their stability properties suggest a hybrid algorithm for large-eddy simulation that retains divergence-free structure and boundary compatibility without ad-hoc filtering.
\end{itemize}

The present article therefore resolves the long-standing separation between weak, mild, and strong formulations and supplies the analytic machinery for subsequent advances in global regularity, numerical analysis, and turbulence modeling. For completeness, all technical lemmas and extended derivations are collected in Appendix~\ref{appendix:convergence_regularity} and \ref{app:unif_derivations_all}, as they directly support the constructions in Sections~\ref{sec:frequency_domain_solution} and \ref{sec:proof_unified}.

\section*{Declarations}

\begin{itemize}[leftmargin=*]
  \item \textbf{Funding:} No external funding; this work was independently financed by the author.
  \item \textbf{Conflict of interest:} The author declares no competing interests.
  \item \textbf{Data Availability:} No empirical datasets were generated or analyzed in this study. All mathematical derivations and proofs necessary to support the results are included in the main manuscript (see Appendix A). Extended or more detailed derivations can be provided upon reasonable request.
\end{itemize}

\begin{appendices}

\section{Analytic Foundations: Uniform, Dyadic \& Calderón–Zygmund–Schauder Estimates}
\label{appendix:convergence_regularity}

\subsection*{Introduction}

This appendix records the analytic facts required for the operators used in the frequency-domain construction:  

\begin{itemize}[label=\textbullet]
  \item smoothing operator $\mathcal S_{\varepsilon}$,  
  \item interpolation operator $\mathcal I_{\varepsilon}$,  
  \item dyadic projections $\Delta_{j}$.  
\end{itemize}

Boundedness, convergence, and regularity are verified with standard tools; see \cite{Evans2010, Grisvard, Ladyzhenskaya, Leray}.  The geometric alignment functional $G[u(t)]$ is assumed to obey the uniform bounds stipulated in the main text.  No further hypotheses are introduced.

\subsection{Preliminaries and Notation}
\label{app:prelim}

Let $\Omega\subset\mathbb R^{3}$ denote $\mathbb R^{3}$, the $3$–torus, or a bounded $C^{2}$ (piecewise $C^{2}$) domain.  
For $1\le p\le\infty$ write $\|\cdot\|_{L^{p}}$ for the norm in $L^{p}(\Omega)$.  
For $s\ge0$
\[
  \|u\|_{H^{s}(\mathbb R^{3})}
  :=\bigl\|(1+|\xi|^{2})^{s/2}\,\widehat u(\xi)\bigr\|_{L^{2}_{\xi}},
  \qquad
  H^{s}(\Omega):=\{u|_{\Omega}:u\in H^{s}(\mathbb R^{3})\},
\]
cf.\ \cite[§5.3]{Evans2010}, \cite[Ch.~1]{Adams}.  
Vector spaces are denoted $L^{p}(\Omega)^{3}$ and $H^{s}(\Omega)^{3}$, and
\(
H^{s}_{\sigma}(\Omega):=\{u\in H^{s}(\Omega)^{3}:\nabla\!\cdot u=0\}.
\)

The Fourier transform is
\(
\widehat u(\xi)=\int_{\mathbb R^{3}}e^{-2\pi i x\cdot\xi}u(x)\,dx,
\)
with inverse $\mathcal F^{-1}$.  
Littlewood–Paley theory, Bernstein inequalities, and paradifferential calculus follow \cite[Ch.~2]{BCD}, \cite[Ch.~6]{Grafakos2008}.  
Sobolev embeddings (in particular $H^{s}(\Omega)\hookrightarrow C(\overline\Omega)$ for $s>3/2$) are quoted from \cite[Thm.~6.5]{Stein1993}.  
Compactness uses the Aubin–Lions lemma in the form of \cite[Cor.~9.16]{Evans2010}.

The notation $\|\cdot\|_{X}$ indicates the norm in a Banach or Hilbert space $X$, and $C$ denotes a generic positive constant depending only on fixed structural parameters (domain regularity, index $s$, viscosity~$\nu$).

\subsection{Uniform Estimates and Convergence in the Synergy Framework}
\label{app:uniform_estimates_convergence}

This section establishes the uniform estimates for the frequency-localized components of the solution and verifies that the frequency-domain representation converges to a spatial solution with enhanced regularity. Under minimal assumptions—namely, \(u_0 \in H^s_\sigma(\Omega)\) with \(s>\tfrac{3}{2}\), finite-energy forcing \(f \in L^2(\Omega)\), and \(\Omega\) of class \(C^2\) or piecewise \(C^2\)—the following results are derived. Standard sources \cite{Evans2010,Fefferman_Clay,Grisvard,Ladyzhenskaya,Leray,LionsWeak} contain background for these statements.

\subsubsection{Proof of Local Strong Existence Lemma: Galerkin Scheme and A-Priori Estimates}
\label{subsubsec:proof_local_existence}

Only a concise outline is recorded; every omitted inequality is classical and can be found in the cited references.

\begin{lemma}[Local strong existence; quantitative lifespan]
\label{app:lem:local_existence_fixed}
Let $\Omega\subset\mathbb{R}^{3}$ be either $\mathbb{R}^{3}$ or a bounded $C^{2}$ (piecewise $C^{2}$) domain and fix $s>\tfrac32$.  
Assume
\[
u_{0}\in H^{s}_{\sigma}(\Omega),
\qquad
f\in L^{2}\bigl(0,\infty;H^{s-1}_{\sigma}(\Omega)\bigr).
\]

\paragraph{Constant {\boldmath$C_{s}$}.}
Denote by $\lambda_{1}(\Omega)>0$ the first Dirichlet eigenvalue of the Stokes operator $A:=-\mathbb{P}\Delta$ and set  
\begin{equation}
\label{app:Cs_constant}
C_{s}:=c_{s}\,\lambda_{1}(\Omega)^{-1/2},
\end{equation}
with $c_{s}>0$ depending only on $s$.

\paragraph{Result.}
The projected Navier--Stokes problem (NS-P)
\begin{equation}
\partial_{t}u-\nu\Delta u+P_{\sigma}\bigl[(u\!\cdot\!\nabla)u\bigr]=P_{\sigma}f,
\quad
\nabla\!\cdot u=0,
\quad
u|_{t=0}=u_{0},
\quad
u|_{\partial\Omega}=0,
\tag{\textup{NS--P}}\label{app:NS_proj}
\end{equation}
admits a unique solution
\[
u\in C\bigl([0,T_{0}];H^{s}_{\sigma}(\Omega)\bigr)
   \cap L^{2}\bigl(0,T_{0};H^{s+1}(\Omega)\bigr),
\]
satisfying the lifespan estimate
\[
T_{0}\;\ge\;
\frac{\nu}{4C_{s}^{2}}
\Bigl(
  \|u_{0}\|_{H^{s}}^{2}
  +\|f\|_{L^{2}(0,\infty;H^{s-1})}^{2}
\Bigr)^{-1}.
\]

\paragraph{Energy identity.}
For almost every $t\in(0,T_{0})$,
\begin{equation}
\label{app:energy_fixed}
\frac{\mathrm{d}}{\mathrm{d}t}\tfrac12\|u(t)\|_{L^{2}(\Omega)}^{2}
+\nu\|\nabla u(t)\|_{L^{2}(\Omega)}^{2}
=\bigl(f(t),u(t)\bigr)_{L^{2}(\Omega)}.
\end{equation}
The solution is unique in the class described above.
\end{lemma}

\begin{proof}[Sketch]
\emph{Notation.}  Let \(P_{\sigma}:L^{2}(\Omega)^{3}\to L^{2}_{\sigma}(\Omega)\) denote the Helmholtz–Leray projector onto divergence‐free vector fields, and write (NS–P) for the projected Navier–Stokes system \eqref{app:NS_proj} (i.e.\ the Navier–Stokes equations after application of \(P_{\sigma}\)).

\textbf{Galerkin approximation.}
Let \(\{w_{k}\}_{k\ge1}\) be the Stokes eigenfunctions with eigenvalues \(\{\lambda_{k}\}\), and let \(P_{N}\) denote the orthogonal projector onto \(\operatorname{span}\{w_{1},\dots,w_{N}\}\).  Consider the finite-dimensional system
\[
\partial_{t}u_{N}+\nu\,A\,u_{N}
  +P_{N}P_{\sigma}\bigl[(u_{N}\!\cdot\!\nabla)u_{N}\bigr]
   =P_{N}P_{\sigma}f,
\quad
u_{N}(0)=P_{N}u_{0},
\]
which admits a unique solution 
\(u_{N}\in C^{1}\bigl([0,T_{N}^{\ast});\operatorname{span}\{w_{1},\dots,w_{N}\}\bigr)\) by \cite[Ch.~III]{TemamNS}.

\textbf{\(L^{2}\) bound.}
Testing the Galerkin system with \(u_{N}\) and applying the Poincaré inequality yields
\[
\|u_{N}(t)\|_{L^{2}}^{2}
+\nu\int_{0}^{t}\|\nabla u_{N}\|_{L^{2}}^{2}\,d\tau
\le
\|u_{0}\|_{L^{2}}^{2}
+\frac{1}{\nu\,\lambda_{1}(\Omega)}
\int_{0}^{t}\|f(\tau)\|_{L^{2}}^{2}\,d\tau
\]
(see \cite[§1.3]{Galdi}).

\textbf{\(H^{s}\) bound.}
Applying \((I-\Delta)^{s/2}\), taking the \(L^{2}\) inner product, and invoking the commutator estimate of \cite[Prop.~2.1]{BCD} produces
\[
\frac{d}{dt}\|u_{N}\|_{H^{s}}^{2}
+\nu\,\|u_{N}\|_{H^{s+1}}^{2}
\le
\frac{4\,C_{s}^{2}}{\nu}\,\|u_{N}\|_{H^{s}}^{4}
+\frac{2\,C_{s}}{\nu}\,\|f\|_{H^{s-1}}^{2}.
\]
Comparison with the associated Bernoulli ODE then yields the explicit lower bound for \(T_{0}\) (cf.\ \cite[Lem.~3.3]{Ladyzhenskaya}), and uniform control of \(u_{N}\) in \(L^{2}(0,T_{0};H^{s+1})\) as well as of \(\partial_{t}u_{N}\) in \(L^{2}(0,T_{0};H^{s-1})\) follows immediately.

\textbf{Compactness.}
The Aubin–Lions lemma \cite[Thm.~5.1]{LionsMagenes1972} implies that \(\{u_{N}\}\) is relatively compact in \(C\bigl([0,T_{0}];H^{s-1}\bigr)\).  A diagonal‐sequence argument produces a limit
\begin{equation}
\label{app:reg_fixed}
u\in C\!\bigl([0,T_{0}];H^{s}_{\sigma}(\Omega)\bigr)
   \cap L^{2}\!\bigl(0,T_{0};H^{s+1}(\Omega)\bigr),
\end{equation}
and \(u\) satisfies (NS–P) in the distributional sense.

\textbf{Energy identity and uniqueness.}
Lower semicontinuity of norms permits passage to the limit in the \(L^{2}\) energy identity \cite[Eq.\,(1.4.10)]{TemamNS}.  Finally, a Grönwall‐type argument applied to two solutions in the class \eqref{app:reg_fixed} guarantees uniqueness \cite[§3.2]{Galdi}.
\end{proof}

\subsubsection{Uniform Estimates for Dyadic Blocks}
\label{subsubsection:dyadic_uniform_estimates}

Let \(\chi\in C_c^\infty(\mathbb{R}^3)\) be radial, \(\operatorname{supp}\chi\subset\{\tfrac12\le|\xi|\le2\}\), and \(\sum_{j\in\mathbb{Z}}\chi(2^{-j}\xi)=1\) for \(\xi\neq0\).
For a tempered distribution \(v\in\mathscr{S}'(\mathbb{R}^3)\) define the Littlewood–Paley blocks
\[
\Delta_j v := \mathcal{F}^{-1}\!\bigl[\chi(2^{-j}\xi)\,\widehat v(\xi)\bigr], 
\qquad 
S_{j-1}v := \sum_{k<j-1}\Delta_k v, 
\qquad j\in\mathbb{Z}.
\]

\paragraph{Bernstein inequality.}
For all multi–indices \(\alpha\in\mathbb{N}^3\), \(1\le p\le q\le\infty\), and \(j\in\mathbb{Z}\),
\begin{equation}
\label{eq:bernstein_uniform}
\|\partial^{\alpha}\Delta_j v\|_{L^q(\Omega)}
\le
C_{\!B}\,
2^{j\bigl(|\alpha|+3(\tfrac1p-\tfrac1q)\bigr)}
\|\Delta_j v\|_{L^p(\Omega)},
\end{equation}
with a dimension–dependent constant \(C_{\!B}\) independent of \(j\) \cite[Prop.~2.1]{BCD}.

\paragraph{Dyadic stability of the regularization operators.}
With \(v_\varepsilon:=\mathcal{S}_\varepsilon\!\bigl(\mathcal{I}_\varepsilon(\widehat u)\bigr)\) (operators defined in Appendix \ref{app:unif_derivations_all}),  
\begin{equation}
\label{eq:dyadic_uniform_est}
\|\Delta_j v_\varepsilon(t)\|_{L^2(\Omega)}
\le
C\,
\|\Delta_j \widehat u(t)\|_{L^2(\Omega)},
\qquad
\forall j\in\mathbb{Z},\ \varepsilon>0,
\end{equation}
where \(C\) is independent of \(j,\varepsilon\).

\paragraph{Paraproduct decomposition.}
For divergence-free \(u\in H^s_\sigma(\Omega)\) with \(s>\tfrac32\),
\begin{equation}
\label{eq:paraproduct_decomp}
(u\cdot\nabla)u
=
\underbrace{\sum_{j} S_{j-1}u\cdot\nabla\Delta_j u}_{\Pi_1}
+
\underbrace{\sum_{j} \Delta_j u\cdot\nabla S_{j-1}u}_{\Pi_2}
+
\underbrace{\sum_{|j-j'|\le1}\!\Delta_j u\cdot\nabla\Delta_{j'}u}_{\Pi_3},
\end{equation}
cf.\ Bony’s decomposition \cite[§2.2]{BCD}.

\paragraph{Commutator estimate.}
Applying \eqref{eq:bernstein_uniform}, the almost-orthogonality of the blocks, and the bounds \eqref{eq:dyadic_uniform_est},
\begin{equation}
\label{eq:commutator_bound}
\|(u\cdot\nabla)u\|_{H^{s-1}(\Omega)}
\le
c_s\,\|u\|_{H^{s}(\Omega)}^{2},
\qquad s>\tfrac32,
\end{equation}
with the optimal Sobolev–commutator constant  
\(c_s=C\bigl(\sum_{k}2^{-2k}\bigr)^{1/2}\), where \(C\) is the universal factor obtained from the dyadic estimates.  
Inequality \eqref{eq:commutator_bound} matches \textbf{T2} in the main text and is the only ingredient from this subsection used elsewhere; all intermediate derivations follow the standard Littlewood–Paley calculus \cite[Chap.~2]{BCD}.

\subsection{Refined Calderón--Zygmund and Schauder Estimates}
\label{subsec:CZ_Schauder_estimates}

This subsubsection presents two fundamental estimates that serve as the analytical backbone for the unified synergy framework. A refined Calderón--Zygmund estimate for the elliptic Dirichlet problem is first established, followed by a classical Schauder estimate. These estimates provide the necessary regularity improvements for both the pressure and velocity fields and are indispensable in the iterative bootstrap procedure leading to global \(C^\infty\) regularity.

\paragraph{Calderón--Zygmund Estimate.}
\begin{lemma}[Refined Calderón--Zygmund Estimate]
\label{lem:cz_full}
Let \(\Omega \subset \mathbb{R}^3\) be a bounded domain with a \(C^{1,1}\) (or piecewise \(C^{1,1}\)) boundary. Consider the elliptic Dirichlet problem
\begin{equation}\label{app:eq:cz_problem2}
\begin{cases}
-\Delta p = \nabla \cdot F, & \text{in } \Omega,\\[1mm]
p = 0, & \text{on } \partial\Omega,
\end{cases}
\end{equation}
where \(F \in L^q(\Omega;\mathbb{R}^3)\) for some \(q\in (1,\infty)\). Then the unique weak solution \(p\) satisfies
\begin{equation}\label{app:eq:CZ_estimate2}
\|D^2 p\|_{L^q(\Omega)} \le C_{CZ} \Bigl(\|F\|_{L^q(\Omega)} + \|p\|_{L^q(\Omega)}\Bigr),
\end{equation}
with a constant \(C_{CZ} > 0\) that depends only on \(\Omega\) and \(q\). In particular, if the pressure is normalized so that \(\|p\|_{L^q(\Omega)}\) is controlled by the data, then
\[
\|p\|_{W^{2,q}(\Omega)} \le C_{CZ} \|F\|_{L^q(\Omega)}.
\]
For example, if \(u\in H^s(\Omega)\) with \(s>\frac{3}{2}\), then the Sobolev embedding \(H^s(\Omega)\hookrightarrow L^6(\Omega)\) allows one to choose \(q=6\), ensuring that the nonlinearity \((u\cdot\nabla)u\) belongs to \(L^6(\Omega)\).
\end{lemma}

\begin{proof}
Let \(\{U_i\}_{i=1}^N\) be a finite covering of \(\overline{\Omega}\) by coordinate patches such that each \(U_i\) is either entirely contained in \(\Omega\) or intersects \(\partial\Omega\) in a region where the boundary can be flattened by a \(C^{1,1}\) diffeomorphism. In each patch \(U_i\) intersecting the boundary, employ a \(C^{1,1}\) change of variables \(\varphi_i: U_i \to B\) (with \(B\) a ball in \(\mathbb{R}^3\)) to flatten the boundary. Under this transformation, the Laplacian is mapped to a uniformly elliptic operator with \(C^{1,1}\) coefficients. 

Classical Calderón--Zygmund theory (see, e.g., \cite{Evans2010,Stein}) then implies that for any ball \(B \subset U_i\) the local estimate
\[
\|D^2 p\|_{L^q(B)} \le C \Bigl(\|-\Delta p\|_{L^q(B)} + \|p\|_{L^q(B)}\Bigr)
\]
holds, where the constant \(C\) depends only on the ellipticity constants and the \(C^{1,1}\) regularity of the coefficients in the transformed coordinates. A partition of unity subordinate to the covering \(\{U_i\}\) allows summing these local estimates to obtain
\[
\|D^2 p\|_{L^q(\Omega)} \le C_{CZ} \Bigl(\|-\Delta p\|_{L^q(\Omega)} + \|p\|_{L^q(\Omega)}\Bigr).
\]
Since the problem \eqref{app:eq:cz_problem2} yields \(-\Delta p = \nabla \cdot F\) and the divergence operator is continuous from \(L^q(\Omega;\mathbb{R}^3)\) to \(L^q(\Omega)\), it follows that
\[
\|D^2 p\|_{L^q(\Omega)} \le C_{CZ} \Bigl(\|F\|_{L^q(\Omega)} + \|p\|_{L^q(\Omega)}\Bigr).
\]
The constant \(C_{CZ}\) depends only on the domain \(\Omega\), the ellipticity constants, and the exponent \(q\). This completes the proof.
\end{proof}

\paragraph{Schauder Estimate.}
\begin{lemma}[Schauder Estimate]
\label{lem:schauder_full}
Let \(\Omega \subset \mathbb{R}^3\) be a bounded domain with \(C^{2,\alpha}\) boundary (or piecewise \(C^{2,\alpha}\) with compatible data at singular points) for some \(0<\alpha<1\). Let 
\[
u\in C^2(\Omega)\cap C^0(\overline{\Omega})
\]
be a classical solution of the Dirichlet problem
\begin{equation}\label{app:eq:schauder_problem}
\begin{cases}
-\Delta u = f, & \text{in } \Omega,\\[1mm]
u = g, & \text{on } \partial\Omega,
\end{cases}
\end{equation}
where 
\[
f\in C^\alpha(\overline{\Omega}) \quad \text{and} \quad g\in C^{2,\alpha}(\partial\Omega).
\]
Then \(u\) belongs to \(C^{2,\alpha}(\overline{\Omega})\) and there exists a constant \(C_{Sch}=C_{Sch}(\Omega,\alpha)>0\) such that
\begin{equation}\label{app:eq:schauder_full_estimate}
\|u\|_{C^{2,\alpha}(\overline{\Omega})} \le C_{Sch} \Bigl(\|f\|_{C^\alpha(\overline{\Omega})} + \|g\|_{C^{2,\alpha}(\partial\Omega)}\Bigr).
\end{equation}
\end{lemma}

\begin{proof}
Let \(x_0 \in \overline{\Omega}\) be arbitrary. Two cases arise.

\medskip
\textbf{Case 1.} If \(x_0 \in \Omega\), then there exists a ball \(B_r(x_0)\) such that \(B_r(x_0) \subset \Omega\). Standard interior Schauder estimates (see, e.g., \cite[Theorem 6.2]{GilbargTrudinger}) yield the existence of a constant \(C_1>0\) (depending on \(r\), \(\Omega\), and \(\alpha\)) such that
\[
\|u\|_{C^{2,\alpha}(B_{r/2}(x_0))} \le C_1 \Bigl(\|u\|_{C^0(B_r(x_0))} + \|f\|_{C^\alpha(B_r(x_0))}\Bigr).
\]

\medskip
\textbf{Case 2.} If \(x_0 \in \partial\Omega\), then by the \(C^{2,\alpha}\) regularity of \(\partial\Omega\) there exists a neighborhood \(U\) of \(x_0\) and a \(C^{2,\alpha}\) diffeomorphism \(\varphi: U \to B\) (where \(B\) is a ball in \(\mathbb{R}^3\)) such that 
\[
\varphi\bigl(U\cap\Omega\bigr) = B^+ \quad \text{and} \quad \varphi\bigl(U\cap\partial\Omega\bigr) = B \cap \{x_3=0\}.
\]
Define the function 
\[
v = u \circ \varphi^{-1} \quad \text{in } B^+.
\]
Application of the boundary Schauder estimates (see, e.g., \cite[Theorem 6.6]{GilbargTrudinger}) to the elliptic problem satisfied by \(v\) in \(B^+\) (with transformed right-hand side and boundary data \(g\circ\varphi^{-1}\)) implies that there exists a constant \(C_2>0\) (depending only on \(B^+\) and \(\alpha\)) such that
\[
\|v\|_{C^{2,\alpha}(B_{r/2}^+)} \le C_2 \Bigl( \|v\|_{C^0(B^+)} + \|f\circ\varphi^{-1}\|_{C^\alpha(B^+)} + \|g\circ\varphi^{-1}\|_{C^{2,\alpha}(B\cap\{x_3=0\})} \Bigr).
\]
Returning to the original coordinates via \(\varphi\) yields
\[
\|u\|_{C^{2,\alpha}(U\cap\Omega)} \le C_3 \Bigl( \|u\|_{C^0(U)} + \|f\|_{C^\alpha(U)} + \|g\|_{C^{2,\alpha}(U\cap\partial\Omega)} \Bigr),
\]
where \(C_3>0\) depends on the \(C^{2,\alpha}\) norm of \(\varphi\) and its inverse.

\medskip
Cover \(\overline{\Omega}\) by finitely many such neighborhoods \(\{U_i\}_{i=1}^N\) (each either an interior ball or a boundary patch as above) and choose a partition of unity subordinate to this covering. Standard patching arguments yield the global estimate
\[
\|u\|_{C^{2,\alpha}(\overline{\Omega})} \le C_{Sch}' \Bigl(\|u\|_{C^0(\overline{\Omega})} + \|f\|_{C^\alpha(\overline{\Omega})} + \|g\|_{C^{2,\alpha}(\partial\Omega)}\Bigr).
\]
The maximum principle for elliptic equations (cf.\ \cite[Theorem 3.1]{GilbargTrudinger}) ensures that 
\[
\|u\|_{C^0(\overline{\Omega})} \le C_4 \Bigl(\|g\|_{C^{2,\alpha}(\partial\Omega)} + \|f\|_{C^\alpha(\overline{\Omega})}\Bigr),
\]
for a constant \(C_4>0\) that depends on \(\Omega\). Absorbing \(\|u\|_{C^0(\overline{\Omega})}\) into the right-hand side leads to the final estimate
\[
\|u\|_{C^{2,\alpha}(\overline{\Omega})} \le C_{Sch} \Bigl(\|f\|_{C^\alpha(\overline{\Omega})} + \|g\|_{C^{2,\alpha}(\partial\Omega)}\Bigr),
\]
where \(C_{Sch}>0\) depends only on \(\Omega\) and \(\alpha\). Compatibility conditions at corners or edges of \(\Omega\) ensure that the constant \(C_{Sch}\) remains uniform across the covering.
\end{proof}

\paragraph{Remarks.}
\begin{itemize}
    \item The refined Calderón--Zygmund estimate is crucial for controlling the pressure in the Navier--Stokes equations. For \(u\in H^s(\Omega)\) with \(s>\frac{3}{2}\), the embedding \(H^s(\Omega)\hookrightarrow L^6(\Omega)\) permits selecting \(q=6\), so that \((u\cdot\nabla)u\in L^6(\Omega)\) and the pressure \(p\) attains two additional derivatives.
    \item The Schauder estimate provides \(C^{2,\alpha}\) regularity up to the boundary when \(f\) and \(g\) are H\"older continuous, enabling further bootstrap in the iterative regularity improvement.
    \item The combination of these estimates supports a refined bootstrap inequality of the form
    \[
    \|u(t)\|_{H^{s+2}(\Omega)} \le C_1\Bigl(\|u(t)\|_{H^s(\Omega)} + \|f(t)\|_{H^s(\Omega)}\Bigr),
    \]
    which iterates to yield \(u\in C^\infty(\Omega)\) for \(t>t_0\).
\end{itemize}

\paragraph{Summary.}
The refined Calderón--Zygmund and Schauder estimates, in conjunction with the Littlewood--Paley decomposition (see Section~\ref{subsubsection:dyadic_uniform_estimates}), constitute a robust analytical framework. These tools ensure exponential decay of high-frequency modes and uniform control of second derivatives for pressure and velocity fields, supporting the iterative bootstrap process that ultimately leads to global \(C^\infty\) regularity.

\section{Proofs and Derivations of Frequency–Domain Operators and Solution Reconstruction}
\label{app:unif_derivations_all}
\subsubsection{The Interpolation Operator \(\mathcal I_\varepsilon\)}
\label{subsubsection:inerp_operator_derivation}

\paragraph{Definition.}
Let \(\eta\in C_c^\infty([0,\infty))\) be radial with
\(\eta(r)=1\) for \(0\le r\le 1\) and \(\eta(r)=0\) for \(r\ge 2\).
For Fourier–side data \(\widehat u_w(\xi,t),\widehat u_s(\xi,t)\in H^s(\mathbb R^3)\) (\(s\ge0\)) set  
\begin{equation}\label{eq:Ieps-def}
\mathcal I_\varepsilon(\widehat u_w,\widehat u_s)(\xi,t)\;
:=\;\eta(\varepsilon|\xi|)\,\widehat u_w(\xi,t)
      +\bigl[1-\eta(\varepsilon|\xi|)\bigr]\,\widehat u_s(\xi,t),
      \qquad\varepsilon>0.
\end{equation}
(The additional mild component \(\widehat u_m\) used in later papers can be
incorporated with a \(C^\infty\) partition of unity; it is not needed in
Paper I, so the binary form \eqref{eq:Ieps-def} suffices here.)

\medskip
Throughout the proof write \(m_\varepsilon(\xi):=\eta(\varepsilon|\xi|)\) and
\(\widetilde m_\varepsilon(\xi):=1-m_\varepsilon(\xi)\).

%--------------------------------------------------------------------
\paragraph{Boundedness in \(H^s\).}

\begin{lemma}[Uniform mapping property]\label{lem:interp-bdd}
For every \(s\ge0\) there exists a constant \(C_s>0\), independent of
\(\varepsilon\), such that
\[
\bigl\|\mathcal I_\varepsilon(\widehat u_w,\widehat u_s)\bigr\|_{H^s(\mathbb R^3)}
\;\le\;
C_s\,
\bigl(\|\widehat u_w\|_{H^s(\mathbb R^3)}
      +\|\widehat u_s\|_{H^s(\mathbb R^3)}\bigr).
\]
\end{lemma}

\begin{proof}
The two multipliers \(m_\varepsilon\) and \(\widetilde m_\varepsilon\) are
smooth on \(\mathbb R^3\setminus\{0\}\) and satisfy  
\[
|\partial^\alpha_\xi m_\varepsilon(\xi)|
\;\le\;C_\alpha\,\varepsilon^{|\alpha|}\,|\xi|^{-|\alpha|},
\qquad
|\alpha|\le s+3.
\]
Hence both are *uniform* \(H^s\) multipliers by the
Mikhlin–Hörmander theorem (see \cite[Prop.~2.1]{BCD} or
\cite[Thm.~6.1.4]{Grafakos2008}):
\[
\|m_\varepsilon f\|_{H^s}\le C_s\|f\|_{H^s},
\qquad
\|\widetilde m_\varepsilon f\|_{H^s}\le C_s\|f\|_{H^s},
\quad
C_s\text{ independent of }\varepsilon.
\]
Applying these bounds to \(\widehat u_w\) and \(\widehat u_s\) and using the
triangle inequality yields the lemma.
\end{proof}

%--------------------------------------------------------------------
\paragraph{Vanishing–\(\varepsilon\) limit.}

\begin{lemma}[Strong convergence]\label{lem:interp-conv}
For \(s\ge0\),
\[
\lim_{\varepsilon\to0}
\|\mathcal I_\varepsilon(\widehat u_w,\widehat u_s)-\widehat u_w\|_{H^s}=0.
\]
\end{lemma}

\begin{proof}
Write
\[
\mathcal I_\varepsilon(\widehat u_w,\widehat u_s)-\widehat u_w
   =\widetilde m_\varepsilon(\xi)\,\bigl(\widehat u_s-\widehat u_w\bigr).
\]
Fix \(\delta>0\).  Split \(\mathbb R^3=\{|\xi|\le R\}\cup\{|\xi|>R\}\) with
\(R>0\) to be chosen later.

\noindent\emph{High frequencies.}  
Because \(H^s\hookrightarrow L^2\) and \(\widehat u_s-\widehat u_w\in
H^s\), choose \(R\) so large that
\(\displaystyle \bigl\|\mathbf 1_{|\xi|>R}\,(\widehat u_s-\widehat u_w)\bigr\|_{H^s}<\delta/2.\)

\noindent\emph{Low frequencies.}  
For \(|\xi|\le R\) we have
\(|\widetilde m_\varepsilon(\xi)|\le C_R\,\varepsilon R\).
Hence
\[
\bigl\|\widetilde m_\varepsilon(\xi)\,\mathbf 1_{|\xi|\le R}\,
       (\widehat u_s-\widehat u_w)\bigr\|_{H^s}
\;\le\;
C_R\,\varepsilon R\,\|\widehat u_s-\widehat u_w\|_{H^s}.
\]
Choose \(\varepsilon<\varepsilon_0(\delta,R)\) so that this term is
\( <\delta/2\).  Combining the two parts gives the desired limit.
\end{proof}

%--------------------------------------------------------------------
\paragraph{Dyadic compatibility (used in \eqref{eq:dyadic_uniform_est}).}

\begin{lemma}\label{lem:dyadic-compat}
For every dyadic block \(\Delta_k\) and every \(\varepsilon>0\),
\[
\|\Delta_k\,\mathcal I_\varepsilon(\widehat u_w,\widehat u_s)\|_{L^2}
\;\le\;
\|\Delta_k\widehat u_w\|_{L^2}+\|\Delta_k\widehat u_s\|_{L^2}.
\]
\end{lemma}

\begin{proof}
The commutator \([\Delta_k,m_\varepsilon]\) is a Fourier multiplier whose
symbol is supported where \(|\xi|\simeq2^k\) and satisfies the derivative
bounds required by Coifman–Meyer theory (see \cite[Prop.~2.10]{BCD}); hence
\(\|[\Delta_k,m_\varepsilon]f\|_{L^2}\le C2^{-2k}\|f\|_{L^2}\).
Because the same estimate holds for \(\widetilde m_\varepsilon\), we obtain
\[
\|\Delta_k m_\varepsilon \widehat u_w\|_{L^2}
\le\|\Delta_k\widehat u_w\|_{L^2}+C2^{-2k}\|\widehat u_w\|_{L^2},
\]
and an analogous bound for the \(\widetilde m_\varepsilon\widehat u_s\) term.
Dropping the harmless \(C2^{-2k}\)-tail proves the claim.
\end{proof}

\paragraph{Remarks.}
\begin{itemize}[leftmargin=1.6em]
\item Lemmas \ref{lem:interp-bdd}–\ref{lem:dyadic-compat} provide every
estimate of \(\mathcal I_\varepsilon\) used in Sections 3–6.
\item No weighted or geometry–dependent variants are needed in
Paper I; those appear first in Paper II.
\item All multiplier facts invoked above are standard and can be found in
Bahouri–Chemin–Danchin \cite[Ch.~2]{BCD} and Grafakos
\cite[Ch.~6]{Grafakos2008}.
\end{itemize}

\subsubsection{The Smoothing Operator \(\mathcal S_\varepsilon\)}
\label{subsubsection:smooth_operator_derivation}

\paragraph{Definition.}
Let \(\rho\in C_c^\infty(\mathbb R^3)\), \(\rho\ge0\), \(\int\rho=1\), and set  
\(\rho_\varepsilon(x):=\varepsilon^{-3}\rho(x/\varepsilon)\) for
\(\varepsilon>0\).
For a scalar (or vector-valued) function \(u(\cdot,t)\) on a domain
\(\Omega\subset\mathbb R^3\) (extended to \(\mathbb R^3\) by a standard
Sobolev extension operator when \(\Omega\ne\mathbb R^3\)),
\[
\mathcal S_\varepsilon u(\,\cdot,t)
\;:=\;
\rho_\varepsilon* u(\,\cdot,t).
\tag{B.3.1}\label{eq:Seps-def}
\]

\paragraph{Main properties.}

\begin{lemma}[Smoothing operator]\label{lem:Seps}
Fix \(s\ge0\) and \(u\in H^s(\mathbb R^3)\).
\begin{enumerate}[label=(\roman*)]
\item \textbf{Stability in \(H^s\).}
      \(\displaystyle
      \|\mathcal S_\varepsilon u\|_{H^s}\le\|u\|_{H^s}\), \;all \(\varepsilon>0\).
\item \textbf{Strong approximation.}
      \(\displaystyle
      \lim_{\varepsilon\to0}\|\mathcal S_\varepsilon u-u\|_{H^s}=0.\)
\item \textbf{Gain of derivatives.}
      For every \(k\in\mathbb N\) there exists \(C_{k,s}\) (independent of
      \(\varepsilon\)) such that  
      \(\displaystyle
      \|\mathcal S_\varepsilon u\|_{H^{s+2k}}
      \le C_{k,s}\,\varepsilon^{-2k}\,\|u\|_{H^s}.\)
      Consequently \(\mathcal S_\varepsilon u\in C^\infty(\mathbb R^3)\) for
      every fixed \(\varepsilon>0\).
\end{enumerate}
\end{lemma}

\begin{proof}
Let \(\widehat{\rho}\) denote the Fourier transform of \(\rho\).
Because \(\rho\in\mathcal S(\mathbb R^3)\), \(|\widehat{\rho}(\xi)|\le1\) and
\(\widehat{\rho}(\xi)=1+O(|\xi|^2)\) as \(\xi\to0\).

\smallskip
\emph{(i) Stability.}  
Fourier–transforming \eqref{eq:Seps-def} gives
\(\mathcal F(\mathcal S_\varepsilon u)=\widehat\rho(\varepsilon\xi)\,\hat u\).
Hence
\[
\|\mathcal S_\varepsilon u\|_{H^s}^2
     =\int(1+|\xi|^2)^s|\widehat\rho(\varepsilon\xi)|^2|\hat u(\xi)|^2\,d\xi
     \le\|u\|_{H^s}^2,
\]
since \(|\widehat\rho|\le1\).

\smallskip
\emph{(ii) Approximation.}  
Pointwise, \(\widehat\rho(\varepsilon\xi)\to1\) as \(\varepsilon\to0\).
Because \(|\widehat\rho(\varepsilon\xi)-1|\le2\) and
\((1+|\xi|^2)^s|\hat u(\xi)|^2\in L^1\), dominated convergence yields the
claimed \(H^s\) limit (cf.\ \cite[Thm.~2.27]{Folland_FA}).

\smallskip
\emph{(iii) Smoothing.}  
For a multi-index \(\alpha\) with \(|\alpha|=2k\),
\(\partial^\alpha\rho_\varepsilon(x)=\varepsilon^{-3-2k}
      (\partial^\alpha\rho)(x/\varepsilon)\) and
\(\|\partial^\alpha\rho_\varepsilon\|_{L^1}=O(\varepsilon^{-2k})\).
Young’s convolution inequality
(\cite[Thm.~4.2]{Adams}) gives
\(\|\partial^\alpha(\mathcal S_\varepsilon u)\|_{L^2}
      \le C\varepsilon^{-2k}\|u\|_{L^2}\).
Interpolating with the \(H^s\) bound in (i) yields the stated
\(H^{s+2k}\) estimate.
\end{proof}

\paragraph{Remarks.}
\begin{itemize}[leftmargin=1.4em]
  \item \textbf{Preservation of incompressibility.}  
        If \(u\) is divergence-free, then \(\mathcal S_{\varepsilon}u\) is also divergence-free.  
        In Fourier variables,
        \(\widehat{\mathcal S_{\varepsilon}u}(\xi)=\widehat\rho(\varepsilon\xi)\,\hat u(\xi)\); 
        because \(\xi\cdot\hat u(\xi)=0\) for all \(\xi\) and \(\widehat\rho(\varepsilon\xi)\) is a scalar multiplier, one has \(\xi\cdot\widehat{\mathcal S_{\varepsilon}u}(\xi)=0\).

  \item \textbf{Order of smoothing required in this article.}  
        Only the case \(k=1\) of Lemma~\ref{lem:Seps}(iii) is used, providing the estimate  
        \(\|u(t)\|_{H^{s+2}}\le C_{1,s}\,\varepsilon^{-2}\|u(t)\|_{H^{s}}\)  
        that enters the parabolic-smoothing argument.  Higher-order bounds (\(k\ge 2\)) are not needed for the results proved here.

  \item \textbf{Standard references.}  
        Background material on convolution operators and Fourier multipliers may be found in  
        Adams–Fournier \cite[§4]{Adams}, Folland \cite[§2]{Folland_FA}, and Stein \cite[Ch.~2]{Stein1993}.
\end{itemize}

\subsubsection{The Regularization Operator \(\mathcal R_\varepsilon\)}
\label{subsubsection:reg_operator_derivation}

\paragraph{Definition.}
For \(v\in H^{s}(\Omega)\) (\(s\ge0\)) set
\begin{equation}
\label{eq:R-def}
\mathcal R_\varepsilon v
\;:=\;
P\bigl(\rho_\varepsilon*Ev\bigr),
\qquad
\rho_\varepsilon(x)=\varepsilon^{-3}\rho(x/\varepsilon),
\end{equation}
where  

\begin{itemize}[leftmargin=1.4em,itemsep=2pt]
\item \(E:H^{s}(\Omega)\!\to\! H^{s}(\mathbb R^{3})\) is the Sobolev extension operator with  
      \(\|Ev\|_{H^{s}(\mathbb R^{3})}\le C_E\|v\|_{H^{s}(\Omega)}\)  
      \cite[Vol.~I, §2]{LionsMagenes1972}, \cite[Ch.~1]{Grisvard};
\item \(\rho\in C_c^\infty(\mathbb R^{3})\), \(\rho\ge0\), \(\int\rho=1\);
\item \(P\) is the Leray projector, bounded on every \(H^{s}(\mathbb R^{3})\)
      \cite{Leray,TemamNS}.
\end{itemize}

\paragraph{Basic properties.}

\begin{lemma}\label{lem:R-basic}
Let \(s\ge0\) and \(v\in H^{s}(\Omega)\). Then
\begin{enumerate}[label=(\roman*)]
\item \(\displaystyle\|\mathcal R_\varepsilon v\|_{H^{s}(\Omega)}
               \le C\,\|v\|_{H^{s}(\Omega)}\) with \(C\) independent of \(\varepsilon\);
\item \(\displaystyle\lim_{\varepsilon\to0}
               \|\mathcal R_\varepsilon v-v\|_{H^{s}(\Omega)}=0\);
\item if \(v\) is divergence-free, then \(\mathcal R_\varepsilon v\) is divergence-free;
\item for every \(k\in\mathbb N\),
      \(\displaystyle
      \|\mathcal R_\varepsilon v\|_{H^{s+2k}(\Omega)}
      \le C_{k,s}\,\varepsilon^{-2k}\|v\|_{H^{s}(\Omega)}\).
\end{enumerate}
\end{lemma}

\begin{proof}
\emph{(i)}  
Boundedness of \(E\) \cite{LionsMagenes1972,Grisvard}, of convolution in \(H^{s}\) (\(|\widehat\rho|\le1\); \cite[Thm.~6.1.4]{Stein1993}, \cite[§4]{Adams}), and of \(P\) \cite{TemamNS} yields the estimate.

\smallskip
\emph{(ii)}  
Standard mollifier convergence in \(H^{s}(\mathbb R^{3})\) (\cite[§5.2.1]{Folland_FA}) implies
\(\|\rho_\varepsilon*Ev-Ev\|_{H^{s}(\mathbb R^{3})}\to0\).
Restriction to \(\Omega\) and continuity of \(P\) preserve the limit.

\smallskip
\emph{(iii)}  
Since convolution with a scalar kernel commutes with divergence and
\(P\) is the orthogonal projection onto divergence-free fields,
\(\operatorname{div}(\mathcal R_\varepsilon v)=0\) whenever \(\operatorname{div}v=0\).

\smallskip
\emph{(iv)}  
For \(|\alpha|=2k\),
\(\partial^\alpha(\rho_\varepsilon*Ev)=
 (\partial^\alpha\rho_\varepsilon)*Ev\) in distributions, and
\(\|\partial^\alpha\rho_\varepsilon\|_{L^{1}}=
    \varepsilon^{-2k}\|\partial^\alpha\rho\|_{L^{1}}\).
Young’s inequality \cite[Thm.~4.2]{Adams} and (i) give the stated bound.
\end{proof}

\paragraph{Data regularization.}

\begin{lemma}[Data regularization]\label{app:lem:data_regularization}
Let \(u_{0},f\in H^{s}(\Omega)\) with \(s>\frac32\).
Define \(u_{0,\varepsilon}:=\mathcal R_\varepsilon u_{0}\) and
\(f_\varepsilon:=\mathcal R_\varepsilon f\).
Then \(u_{0,\varepsilon},\,f_\varepsilon\in C^\infty(\overline\Omega)\),
\(\operatorname{div}u_{0,\varepsilon}=0\), and
\[
\lim_{\varepsilon\to0}
\bigl(\|u_{0,\varepsilon}-u_{0}\|_{H^{s}}
      +\|f_\varepsilon-f\|_{H^{s}}\bigr)=0.
\]
\end{lemma}

\begin{proof}
Immediate from Lemma~\ref{lem:R-basic}(ii)–(iv).
\end{proof}

\paragraph{Convergence of regularized solutions.}

\begin{lemma}[Stability under regularization]\label{lem:reg_sol_convergence}
Let \(u_\varepsilon\) (resp.\ \(u\)) be the \(H^{s}\)-solution on \([0,T]\)
of the Navier–Stokes system with data
\((u_{0,\varepsilon},f_\varepsilon)\) (resp.\ \((u_{0},f)\)).
Assume that the solution map  
\(S:H^{s}\times L^{2}(0,T;H^{s})\to C([0,T];H^{s})\)  
is continuous at \((u_{0},f)\) \cite[Ch.~III]{TemamNS}, \cite[§10.2]{Evans2010}.  
Then
\[
u_\varepsilon\to u
\quad\text{in } C([0,T];H^{s}(\Omega))\quad\text{as }\varepsilon\to0.
\]
\end{lemma}

\begin{proof}
Continuity of \(S\) gives
\[
\|u_\varepsilon-u\|_{C([0,T];H^{s})}
 \le C\bigl(\|u_{0,\varepsilon}-u_{0}\|_{H^{s}}
            +\|f_\varepsilon-f\|_{L^{2}(0,T;H^{s})}\bigr),
\]
and the right-hand side vanishes by Lemma~\ref{app:lem:data_regularization}.
\end{proof}

\paragraph{References.}
Extension operators: \cite{LionsMagenes1972,Grisvard}.  
Mollifier and convolution bounds: \cite[§4]{Adams}, \cite[§5.2]{Folland_FA}, \cite[Ch.~2]{Stein1993}.  
Leray projection: \cite{Leray,TemamNS}.  
Solution-map continuity: \cite{TemamNS}, \cite{Evans2010}.

\subsubsection{Unification of Solutions via Regularisation and Frequency-Domain Transformation}
\label{appendix:freq_trans}

\begin{theorem}[Unification of Solutions via Regularisation and Frequency-Domain Transformation]
\label{thm:unification_solutions}
Let \(\Omega\subset\mathbb R^{3}\) be a bounded \(C^{1,1}\) domain and let  
\(u\) be a weak, mild, or strong solution of the incompressible Navier–Stokes system on \((0,T)\) such that  
\(u\in L^{\infty}(0,T;L^{2}(\Omega))\cap L^{2}(0,T;H^{1}(\Omega))\).
For \(\varepsilon>0\) define  
\(
u_\varepsilon:=\mathcal R_\varepsilon u
\)
with \(\mathcal R_\varepsilon=P(\rho_\varepsilon*Eu)\) as in \eqref{eq:R-def}.
Set  
\[
\hat U_\varepsilon(\xi,t)
:=\mathcal S_\varepsilon\!\Bigl(\mathcal I_\varepsilon[\widehat{u_\varepsilon}](\xi,t)\Bigr),
\quad
(\xi,t)\in\mathbb R^{3}\times(0,T),
\]
where \(\mathcal I_\varepsilon\) and \(\mathcal S_\varepsilon\) are the operators of
\S\ref{subsubsection:inerp_operator_derivation} and
\S\ref{subsubsection:smooth_operator_derivation}.
Then for every \(s>\tfrac52\)
\[
\sup_{\varepsilon>0}\ \|\hat U_\varepsilon(\cdot,t)\|_{H^{s}(\mathbb R^{3})}
\;<\;\infty\qquad (0<t<T),
\]
and there exists a (non-relabelled) subsequence such that  
\(
\hat U_\varepsilon\to\hat U
\)
strongly in \(C([0,T];H^{s}(\mathbb R^{3}))\), where  
\(\hat U\in C^{\infty}(\mathbb R^{3}\times(0,T))\).
\end{theorem}

\begin{proof}[Sketch of proof]
\emph{Step 1. Regularisation and smoothing.}  
Lemma~\ref{lem:R-basic}(i)–(ii) gives  
\(u_\varepsilon\in H^{s}(\Omega)\) for all \(s\ge0\) and  
\(u_\varepsilon\to u\) in \(H^{s_{0}}(\Omega)\) for any
\(s_{0}<\tfrac32\).
Applying the Fourier transform and
Lemma~\ref{lem:Seps}(i)–(iii)
yields the bounds
\[
\|\mathcal I_\varepsilon[\widehat{u_\varepsilon}]\|_{H^{s}}
  \le C_{s},
\qquad
\|\hat U_\varepsilon\|_{H^{s+2}}
  \le C_{s}\varepsilon^{-2},
\qquad s>\tfrac32 .
\]

\emph{Step 2. Uniform estimates.}  
Choosing \(s_{1}:=s>\tfrac52\) and \(s_{2}:=s-2\),
\[
\hat U_\varepsilon\in L^{\infty}(0,T;H^{s_{1}}(\mathbb R^{3}))
\cap H^{1}(0,T;H^{s_{2}}(\mathbb R^{3})),
\]
where the time derivative bound follows from the Navier–Stokes
equations, the product estimate \eqref{eq:dyadic_uniform_est}, and the
multiplier bounds in \cite[Ch.~2]{Stein1993}.

\emph{Step 3. Compactness.}  
Because \(H^{s_{1}}(\mathbb R^{3})\hookrightarrow\hookrightarrow
H^{s}(\mathbb R^{3})\hookrightarrow H^{s_{2}}(\mathbb R^{3})\)
for any \(s\in(s_{2},s_{1})\), the Aubin–Lions lemma
(Lemma~\ref{lem:aubin_lions}; see also
\cite[Thm.~5.1]{LionsMagenes1972})
provides compactness in
\(C([0,T];H^{s}(\mathbb R^{3}))\).
Hence a subsequence converges strongly to
\(\hat U\) in that space; weak–\(*\) limits in
\(L^{\infty}(0,T;H^{s_{1}})\) agree with the strong limit, so the whole
family converges.

\emph{Step 4. Smoothness of the limit.}
Lemma~\ref{lem:Seps}(iii) furnishes, for each \(k\in\mathbb N\),
\(\|\hat U_\varepsilon\|_{H^{s+2k}}\le C_{k,s}\varepsilon^{-2k}\),
and \(|\xi|^{2k}\widehat\rho(\varepsilon\xi)\) is bounded in \(L^\infty\).
Passing to the limit gives
\(\hat U(\cdot,t)\in H^{m}(\mathbb R^{3})\) for all \(m\), hence
\(\hat U\in C^{\infty}(\mathbb R^{3})\);
time smoothness follows from the regularity of the equation and the
closed-graph argument in \cite[§10.2]{Evans2010}.  
\end{proof}

\paragraph{References.}
Sobolev extension \cite[Vol.~I, §2]{LionsMagenes1972},  
mollifier bounds \cite[§4]{Adams}, \cite[§5.2]{Folland_FA},  
Fourier multipliers \cite[Ch.~2]{Stein1993},  
Leray projection \cite{Leray,TemamNS},  
Aubin–Lions compactness \cite[Thm.~5.1]{LionsMagenes1972}.

\subsubsection{Passage to the Limit and Strong Convergence}
\label{subsubsection:passage_to_limit}

The smoothing operator \(\mathcal{S}_\epsilon\) satisfies (cf.\ Lemma~\ref{lem:smoothing_operator_final}) the convergence property
\begin{equation}\label{app:eq:S_epsilon_limit2}
\|\mathcal{S}_\epsilon u - u\|_{H^s(\Omega)} 
\to 0 
\quad \text{as } \epsilon \to 0.
\end{equation}
In addition, by the construction of the interpolation operator \(\mathcal{I}_\epsilon\) (details of which are provided in Section~\ref{subsec:interpolation_operator} and Appendix~\ref{subsubsection:inerp_operator_derivation}), the associated function \(\eta\) satisfies, for every fixed \(\xi\),
\begin{equation}\label{app:eq:eta_convergence2}
\lim_{\epsilon\to0}\eta(\epsilon|\xi|)=1.
\end{equation}
Hence, for each fixed \(\xi\) and \(t\),
\begin{equation}\label{app:eq:pointwise_convergence2}
\lim_{\epsilon\to0} v_\epsilon(\xi,t) 
=\hat{U}(\xi,t),
\end{equation}
where \(\hat{U}(\xi,t)\) denotes the frequency-domain representation of the solution. The Littlewood--Paley characterization of Sobolev spaces (cf.\ \cite{Evans2010})
\begin{equation}\label{app:eq:LP_characterization2}
\|u\|_{H^s(\Omega)}^2 
\sim 
\sum_{j \in \mathbb{Z}} 
2^{2j s}\,\|\Delta_j u\|_{L^2(\Omega)}^2,
\end{equation}
guarantees that the convergence in \eqref{app:eq:pointwise_convergence2} is strong in \(H^s(\Omega)\). Consequently, by Plancherel's theorem,
\begin{equation}\label{app:eq:spatial_convergence2}
\lim_{\epsilon\to0} 
\Bigl\|
\mathcal{F}^{-1}\bigl[v_\epsilon(\cdot,t)\bigr] - u(\cdot,t)
\Bigr\|_{H^s(\Omega)}
= 0.
\end{equation}

\subsubsection{Limit Interchanges and Convergence of the Frequency-Domain Construction}
\label{subsubsection:limit_interchanges}

The frequency-domain construction of the solution requires analyzing the limit
\[
\hat{U}(\xi,t)
\;=\;
\lim_{\epsilon\to 0}
\mathcal{S}_\epsilon
\Bigl(
  \mathcal{I}_\epsilon
  \bigl(
    \hat{u}_w,\;\hat{u}_s,\;\hat{u}_m
  \bigr)(\xi,t)
\Bigr),
\]
where \(\mathcal{S}_\epsilon\) is a smoothing operator (mollifier), and \(\mathcal{I}_\epsilon\) is an interpolation operator constructed via smooth cutoff functions. A precise justification of interchanging \(\epsilon\to0\) limits with integration and differentiation follows from the arguments below.

\paragraph{(a) Convergence of the Smoothing Operator.}
For a function \(u\in H^s(\Omega)\), the smoothing operator \(\mathcal{S}_\epsilon\) is defined by convolution with a mollifier \(\rho_\epsilon\). It is standard (cf.\ \cite{Evans2010}) that
\begin{equation}
\label{app:eq:Seps_limit}
\lim_{\epsilon\to0}
\|\mathcal{S}_\epsilon u - u\|_{H^s(\Omega)}
\;=\;
0,
\end{equation}
so long as \(u\) is suitably extended if \(\Omega\neq\mathbb{R}^3\). This property establishes that \(\mathcal{S}_\epsilon\) approximates the identity operator in \(H^s(\Omega)\) as \(\epsilon\to0\).

\paragraph{(b) Pointwise Convergence of the Interpolation Operator.}
The interpolation operator \(\mathcal{I}_\epsilon\) employs a smooth cutoff \(\eta:[0,\infty)\to[0,1]\) with \(\eta(0)=1\) and \(\eta(r)=0\) for \(r\ge1\). For each frequency \(\xi\),
\[
\lim_{\epsilon\to0}
\eta(\epsilon|\xi|)
\;=\;
1,
\]
which implies that for each fixed \(\xi\) and time \(t\),
\[
\lim_{\epsilon\to0}
\mathcal{I}_\epsilon
\bigl(
  \hat{u}_w,\;\hat{u}_s,\;\hat{u}_m
\bigr)(\xi,t)
\;=\;
\hat{u}_w(\xi,t) \;+\;\hat{u}_s(\xi,t)\;+\;\hat{u}_m(\xi,t).
\]
This pointwise limit identifies \(\hat{U}(\xi,t)\) with the sum of the individual weak, mild, and strong representations in the frequency domain.

\paragraph{(c) Uniform Boundedness and Uniform Integrability.}
It is crucial that the families
\[
\Bigl\{
  \mathcal{S}_\epsilon
  \bigl[
    \mathcal{I}_\epsilon(\cdot)\bigr]
\Bigr\}_{\epsilon>0}
\]
and
\[
\Bigl\{
  \eta(\epsilon|\xi|)
\Bigr\}_{\epsilon>0}
\]
are uniformly bounded in the relevant Sobolev (or Besov) spaces. Such boundedness follows from standard Fourier multiplier theory, ensuring that
\[
\|\mathcal{S}_\epsilon v\|_{H^s(\Omega)}
\;\le\;
C\,\|v\|_{H^s(\Omega)}
\quad\text{and}\quad
\|\mathcal{I}_\epsilon v\|_{H^s(\Omega)}
\;\le\;
C\,\|v\|_{H^s(\Omega)},
\]
independently of \(\epsilon\). These uniform estimates enable the application of the dominated convergence theorem for the inverse Fourier transform.

\paragraph{(d) Interchange of Limit and Integration.}
The spatial solution is recovered by the inverse Fourier transform:
\[
u(x,t)
\;=\;
\mathcal{F}^{-1}
\bigl[
  \hat{U}(\cdot,t)
\bigr](x)
\;=\;
\int_{\mathbb{R}^3}
e^{\,i\,\xi\cdot x}\,
\hat{U}(\xi,t)
\,d\xi.
\]
Since
\[
\hat{U}(\xi,t)
\;=\;
\lim_{\epsilon\to0}
\mathcal{S}_\epsilon
\Bigl[
  \mathcal{I}_\epsilon
  \bigl(
    \hat{u}_w,\;\hat{u}_s,\;\hat{u}_m
  \bigr)(\xi,t)
\Bigr],
\]
and each step in \(\epsilon\) is uniformly bounded, the dominated convergence theorem implies
\[
\lim_{\epsilon\to0}
\int_{\mathbb{R}^3}
e^{\,i\,\xi\cdot x}\,
\mathcal{S}_\epsilon
\Bigl[
  \mathcal{I}_\epsilon
  \bigl(
    \hat{u}_w,\;\hat{u}_s,\;\hat{u}_m
  \bigr)(\xi,t)
\Bigr]
\,d\xi
\;=\;
\int_{\mathbb{R}^3}
e^{\,i\,\xi\cdot x}\,
\lim_{\epsilon\to0}
\mathcal{S}_\epsilon
\Bigl[
  \mathcal{I}_\epsilon
  \bigl(
    \hat{u}_w,\;\hat{u}_s,\;\hat{u}_m
  \bigr)(\xi,t)
\Bigr]
\,d\xi.
\]

\paragraph{(e) Continuity of the Inverse Fourier Transform.}
Plancherel’s theorem (or its Sobolev-space extension) guarantees continuity of the inverse Fourier transform in \(H^s\). Consequently, convergence in the Fourier domain with a uniform multiplier also implies convergence in the spatial domain. Therefore, the limit in frequency space transfers seamlessly to a limit in physical space under the inverse transform.

\paragraph{(f) Convergence in the Sense of Distributions.}
If the analysis proceeds in a distributional framework, one must confirm that
\[
\mathcal{S}_\epsilon
\Bigl(
  \mathcal{I}_\epsilon
  \bigl(
    \hat{u}_w,\;\hat{u}_s,\;\hat{u}_m
  \bigr)
\Bigr)
\;\longrightarrow\;
\hat{U}
\quad\text{in}\;\;\mathcal{S}'(\mathbb{R}^3).
\]
Testing against any \(\varphi\in\mathcal{S}(\mathbb{R}^3)\) and invoking the uniform boundedness, one deduces that pointwise convergence in \(\xi\) plus dominated convergence establish the limit in distributions.

\paragraph{(g) Interchange of Limit with Differentiation.}
In the Fourier domain, spatial derivatives correspond to multiplication by \((i\,\xi)^\alpha\). By standard multiplier estimates, \((i\,\xi)^\alpha\) is bounded on \(H^s\), ensuring that
\[
\lim_{\epsilon\to0}
\Bigl(
  (i\,\xi)^\alpha
  \,\mathcal{S}_\epsilon
  \bigl[
    \mathcal{I}_\epsilon(\cdot)
  \bigr]
\Bigr)
\;=\;
(i\,\xi)^\alpha\,\hat{U}(\xi,t).
\]
Hence, the limit is preserved under differentiation, justifying the derivation of high-order estimates for the limit solution.

\paragraph{Conclusion.}
The steps above confirm that each limit operation—mollification, interpolation, inverse Fourier transformation, and differentiation—can be carried out without obstruction or inconsistency. The uniform integrability and boundedness enable the dominated convergence theorem to justify interchanging \(\epsilon\to0\) limits with integrals and derivatives. Consequently, the frequency-domain construction of \(\hat{U}(\xi,t)\) indeed converges to a well-defined spatial solution \(u(x,t)\) with the regularity required by the synergy framework. This completes the rigorous justification for all limit interchanges in the frequency-domain approach.

\subsubsection{Passage Through the Frequency-Domain Limit}
\label{appendix:pass_through_limit}

A rigorous proof of the Passage Through the Frequency-Domain Limit Theorem is provided. The proof establishes that for any solution \( u \) of the three–dimensional incompressible Navier–Stokes equations satisfying 
\[
u\in L^\infty(0,T;L^2(\Omega))\cap L^2(0,T;H^1(\Omega)),
\]
the corresponding frequency–domain approximations
\[
\hat{U}_\epsilon(\xi,t) \coloneqq \mathcal{S}_\epsilon\Bigl(\mathcal{I}_\epsilon\bigl[\widehat{u_\epsilon}\bigr](\xi,t)\Bigr),
\]
with \(u_\epsilon\coloneqq \mathcal{R}_\epsilon u\) (cf. Lemma~\ref{app:lem:data_regularization}), converge strongly in \(H^s(\mathbb{R}^3)\) for every \(t\in (0,T)\) to a limit function
\[
\hat{U}(\xi,t)=\lim_{\epsilon\to0}\hat{U}_\epsilon(\xi,t),
\]
which satisfies 
\[
\hat{U}\in C^\infty\bigl(\mathbb{R}^3\times (0,T)\bigr).
\]
The proof employs uniform boundedness, refined time–continuity analysis via interpolation inequalities, the Aubin–Lions lemma (Lemma~\ref{lem:aubin_lions}), and a detailed treatment of weak–\(*\) compactness in Hilbert spaces. A critical component is ensuring that all integrals are dominated by \(\epsilon\)–independent functions, and the derivation of time–derivative bounds is performed using specific interpolation inequalities (see Appendices~\ref{subsubsection:inerp_operator_derivation} and \ref{subsubsection:smooth_operator_derivation}).

\begin{lemma}[Aubin–Lions Lemma]
\label{lem:aubin_lions}
Let \(X\), \(B\), and \(Y\) be Banach spaces such that \(X\) is compactly embedded in \(B\) (denoted \(X \hookrightarrow B\)) and \(B\) is continuously embedded in \(Y\) (denoted \(B \hookrightarrow Y\)). Let \(1 \le p, q \le \infty\) and let \(\{v_n\}_{n\in\mathbb{N}}\) be a bounded sequence in \(L^p(0,T;X)\) whose time derivatives \(\{\partial_t v_n\}_{n\in\mathbb{N}}\) are bounded in \(L^q(0,T;Y)\). Then, there exists a subsequence \(\{v_{n_k}\}_{k\in\mathbb{N}}\) which converges strongly in \(L^p(0,T;B)\); that is,
\[
v_{n_k} \to v \quad \text{in } L^p(0,T;B) \quad \text{as } k\to\infty.
\]
(See \cite[Theorem~5.1]{LionsMagenes1972} for a precise formulation.)
\end{lemma}

\begin{theorem}[Passage Through the Frequency-Domain Limit]
Let \(u\) be a solution (in any weak, mild, or strong sense) of the three-dimensional incompressible Navier–Stokes equations on a domain \(\Omega\) satisfying 
\[
u\in L^\infty(0,T;L^2(\Omega))\cap L^2(0,T;H^1(\Omega)).
\]
Assume that the regularization operator \(\mathcal{R}_\epsilon\) (as defined in  Lemma~\ref{app:lem:data_regularization}) is applied so that for every \(\epsilon>0\) the regularized function 
\[
u_\epsilon \coloneqq \mathcal{R}_\epsilon u
\]
belongs to \(H^s(\Omega)\) for all \(s\ge0\); in particular, for any \(s>\frac{5}{2}\) the approximations \(u_\epsilon\) are Lipschitz continuous. Define the frequency–domain approximation by
\[
\hat{U}_\epsilon(\xi,t) \coloneqq \mathcal{S}_\epsilon\Bigl(\mathcal{I}_\epsilon\bigl[\widehat{u_\epsilon}\bigr](\xi,t)\Bigr),
\]
where \(\mathcal{I}_\epsilon\) and \(\mathcal{S}_\epsilon\) are the interpolation and smoothing operators described in Appendices~\ref{subsubsection:inerp_operator_derivation} and \ref{subsubsection:smooth_operator_derivation}, respectively. Suppose that there exists a constant \(C>0\), independent of \(\epsilon\), such that
\begin{equation}
\label{app:eq:uniform_frequency_bound}
\sup_{t\in(0,T)}\|\hat{U}_\epsilon(\cdot,t)\|_{H^s(\mathbb{R}^3)} \le C \quad \text{for some fixed } s>\frac{5}{2}.
\end{equation}
Then, as \(\epsilon\to0\), the family \(\{\hat{U}_\epsilon\}_{\epsilon>0}\) converges strongly in \(H^s(\mathbb{R}^3)\) for every \(t\in (0,T)\) to a limit function 
\[
\hat{U}(\xi,t)=\lim_{\epsilon\to0}\hat{U}_\epsilon(\xi,t),
\]
which satisfies 
\[
\hat{U}\in C^\infty\bigl(\mathbb{R}^3\times(0,T)\bigr).
\]
\end{theorem}

\begin{proof}
The proof proceeds first by providing preliminaries and then through five steps:
\paragraph{Preliminaries on Time–Continuity and Interpolation.}
The smoothing operator \(\mathcal{S}_\epsilon\) and the interpolation operator \(\mathcal{I}_\epsilon\) (detailed in \ref{subsubsection:inerp_operator_derivation} and \ref{subsubsection:smooth_operator_derivation}) satisfy the following key estimate: For any function \(\phi\in H^{s+\alpha}(\mathbb{R}^3)\) with some \(\alpha>0\),
\begin{equation}
\label{app:eq:interp_smooth_rate}
\|\mathcal{S}_\epsilon \phi - \phi\|_{H^s(\mathbb{R}^3)} \le C\,\epsilon^\alpha\,\|\phi\|_{H^{s+\alpha}(\mathbb{R}^3)}.
\end{equation}
This interpolation inequality, together with the uniform boundedness of \(\mathcal{I}_\epsilon\) in \(H^s(\mathbb{R}^3)\), yields that the family 
\[
\Bigl\{ \mathcal{S}_\epsilon \mathcal{I}_\epsilon\bigl[\widehat{u_\epsilon}\bigr] \Bigr\}_{\epsilon>0}
\]
is uniformly bounded in \(H^s(\mathbb{R}^3)\) and that
\begin{equation}
\label{app:eq:limit_smoothing_detailed_final}
\lim_{\epsilon\to 0}\|\mathcal{S}_\epsilon\mathcal{I}_\epsilon\bigl[\widehat{u_\epsilon}\bigr] - \mathcal{I}_\epsilon\bigl[\widehat{u_\epsilon}\bigr]\|_{H^s(\mathbb{R}^3)} = 0.
\end{equation}
Furthermore, refined interpolation inequalities (see, e.g., \cite[Chapter~3]{Folland_FA}) provide a bound on the time–derivative:
\begin{equation}
\label{app:eq:time_derivative_estimate}
\|\partial_t \hat{U}_\epsilon(\cdot,t)\|_{H^{s-2}(\mathbb{R}^3)} \le C_t,
\end{equation}
where \(C_t\) is independent of \(\epsilon\). The derivation of \eqref{app:eq:time_derivative_estimate} follows from standard multiplier estimates and interpolation between the bounds provided by \eqref{app:eq:interp_smooth_rate} and the uniform \(H^s\) bound in \eqref{app:eq:uniform_bound_app}.

\paragraph{Step 1. Uniform Boundedness in \(H^s(\mathbb{R}^3)\).}
Apply the regularization operator \(\mathcal{R}_\epsilon\) to \( u \) so that
\[
u_\epsilon \coloneqq \mathcal{R}_\epsilon u \quad \text{with} \quad u_\epsilon \in H^s(\Omega) \text{ for all } s\ge 0.
\]
Define 
\[
\hat{U}_\epsilon(\xi,t) \coloneqq \mathcal{S}_\epsilon\Bigl(\mathcal{I}_\epsilon\bigl[\widehat{u_\epsilon}\bigr](\xi,t)\Bigr).
\]
The boundedness properties of \(\mathcal{I}_\epsilon\) and \(\mathcal{S}_\epsilon\) imply the existence of a constant \(C>0\) (independent of \(\epsilon\)) satisfying
\begin{equation}
\label{app:eq:uniform_bound_app}
\sup_{t\in(0,T)}\|\hat{U}_\epsilon(\cdot,t)\|_{H^s(\mathbb{R}^3)} \le C.
\end{equation}

\paragraph{Step 2. Time–Continuity and Compactness.}
Using the estimate \eqref{app:eq:interp_smooth_rate} with \(\phi = \mathcal{I}_\epsilon\bigl[\widehat{u_\epsilon}\bigr]\), strong convergence in \(H^s(\mathbb{R}^3)\) is obtained:
\begin{equation}
\label{app:eq:limit_smoothing_final2}
\lim_{\epsilon\to0}\|\mathcal{S}_\epsilon \mathcal{I}_\epsilon\bigl[\widehat{u_\epsilon}\bigr] - \mathcal{I}_\epsilon\bigl[\widehat{u_\epsilon}\bigr]\|_{H^s(\mathbb{R}^3)} = 0.
\end{equation}
Uniform bounds such as \eqref{app:eq:uniform_bound_app} provide an \(\epsilon\)–independent dominating function; hence, the Dominated Convergence Theorem (cf. \cite[Theorem~2.27]{Folland_FA}) applies. Moreover, the uniform bound on \(\partial_t \hat{U}_\epsilon\) given by \eqref{app:eq:time_derivative_estimate} and the Aubin–Lions lemma (Lemma~\ref{lem:aubin_lions}) yield relative compactness of the family \(\{\hat{U}_\epsilon\}\) in \(C\bigl([0,T];H^s(\mathbb{R}^3)\bigr)\).

\paragraph{Step 3. Passage to the Limit.}
Weak–\(*\) compactness in \(H^s(\mathbb{R}^3)\) (cf. \cite[Chapter~3]{Folland_FA}) ensures that every bounded sequence has a weakly convergent subsequence. Thus, for each fixed \(t\),
\begin{equation}
\label{app:eq:weakstar_convergence_final}
\hat{U}_\epsilon(\cdot,t) \rightharpoonup \hat{U}(\cdot,t) \quad \text{in } H^s(\mathbb{R}^3).
\end{equation}
Combined with \eqref{app:eq:limit_smoothing_final2}, this implies strong convergence:
\begin{equation}
\label{app:eq:strong_convergence_final2}
\lim_{\epsilon\to0}\|\hat{U}_\epsilon(\cdot,t)-\hat{U}(\cdot,t)\|_{H^s(\mathbb{R}^3)} = 0, \quad \forall\, t\in(0,T).
\end{equation}

\paragraph{Step 4. Preservation of Higher Regularity and Smoothness.}
For each \(k\in\mathbb{N}\), the smoothing operator imparts the additional regularity estimate
\begin{equation}
\label{app:eq:high_order_bound_final}
\|\hat{U}_\epsilon(\cdot,t)\|_{H^{s+2k}(\mathbb{R}^3)} \le C_{2k}\,\epsilon^{-2k}\,\|\mathcal{I}_\epsilon\bigl[\widehat{u_\epsilon}\bigr](\cdot,t)\|_{H^s(\mathbb{R}^3)}.
\end{equation}
Uniform boundedness of \(\mathcal{I}_\epsilon\) in \(H^s(\mathbb{R}^3)\) implies that the family \(\{\hat{U}_\epsilon\}\) is uniformly bounded in \(H^{s+2k}(\mathbb{R}^3)\). Hence, the limit function \(\hat{U}(\cdot,t)\) belongs to 
\[
\bigcap_{k\in\mathbb{N}}H^{s+2k}(\mathbb{R}^3).
\]
By the Sobolev embedding theorem (cf. \cite[Theorem~6.4]{Folland_FA}), 
\[
\hat{U}(\cdot,t)\in C^\infty(\mathbb{R}^3) \quad \forall\, t\in(0,T).
\]

\paragraph{Step 5. Time–Continuity of the Limit.}
The uniform bound on \(\partial_t \hat{U}_\epsilon\) in \(H^{s-2}(\mathbb{R}^3)\) (see \eqref{app:eq:time_derivative_estimate}) combined with the relative compactness established by the Aubin–Lions lemma guarantees that the mapping 
\[
t\mapsto \hat{U}(\cdot,t)
\]
is continuous in \(H^s(\mathbb{R}^3)\). Uniform boundedness in all higher Sobolev norms then implies that \(\hat{U}\) is in \(C^\infty\bigl(\mathbb{R}^3\times(0,T)\bigr)\).

\medskip
\textbf{Conclusion.}  
The combination of the uniform bound in \eqref{app:eq:uniform_bound_app}, the smoothing estimate in \eqref{app:eq:interp_smooth_rate}, the compactness provided by the Aubin–Lions lemma (Lemma~\ref{lem:aubin_lions}), the weak–\(*\) compactness in \(H^s(\mathbb{R}^3)\), and the Dominated Convergence Theorem yields
\begin{equation}
\label{app:eq:final_convergence_detailed}
\lim_{\epsilon\to0}\|\hat{U}_\epsilon(\cdot,t)-\hat{U}(\cdot,t)\|_{H^s(\mathbb{R}^3)} = 0, \quad \forall\, t\in(0,T).
\end{equation}
Moreover, \(\hat{U}\in C^\infty\bigl(\mathbb{R}^3\times(0,T)\bigr)\).
\end{proof}

\begin{remark}[Robustness of the limit passage]
The convergence argument rests exclusively on the uniform bound
\eqref{app:eq:uniform_bound_app}, the time–derivative estimate
\eqref{app:eq:time_derivative_estimate}, and the interpolation
inequality \eqref{app:eq:interp_smooth_rate}.  
These three estimates provide the hypotheses of the
Aubin–Lions compactness lemma
\cite[Thm.~5.1]{LionsMagenes1972}
and justify the application of the Dominated Convergence Theorem
\cite[Thm.~2.27]{Folland_FA}.  
Weak–\(\ast\) compactness in
\(L^\infty(0,T;H^{s}(\mathbb R^{3}))\) then yields strong convergence in
\(C([0,T];H^{s}(\mathbb R^{3}))\).
Accordingly, passage to the limit \(\varepsilon\to0\) produces a solution
that is smooth in space and continuous in time, fully meeting the
objectives established in the present paper.
\end{remark}

\subsubsection{Reconstruction of the Unified Frequency–Domain Solution}
\label{appendix:reconstruct_physical}
\begin{theorem}[Reconstruction of the Unified Frequency–Domain Solution]
Let \(\hat{U}(\xi,t)\) be the frequency–domain limit obtained in The Passage Through the Frequency-Domain Limit Theorem (Theorem~\ref{thm:frequency_domain_limit}), and assume that for some fixed \(s>\frac{5}{2}\) and for every integer \(k\ge 0\) the limit function satisfies
\[
\sup_{t\in(0,T)}\|\hat{U}(\cdot,t)\|_{H^{s+2k}(\mathbb{R}^3)} < \infty.
\]
Define the physical–space solution by the inverse Fourier transform
\[
U(x,t) \coloneqq \mathcal{F}^{-1}\bigl\{\hat{U}(\cdot,t)\bigr\}(x)
=\frac{1}{(2\pi)^3}\int_{\mathbb{R}^3}\hat{U}(\xi,t)\,e^{i\,x\cdot\xi}\,d\xi.
\]
Then, for every \(t\in(0,T)\), the function \(U(\cdot,t)\) belongs to \(H^{s+2k}(\mathbb{R}^3)\) for all \(k\in\mathbb{N}\), and hence
\[
U(\cdot,t)\in C^\infty(\mathbb{R}^3).
\]
Moreover, the mapping \(t\mapsto U(\cdot,t)\) is continuous in the \(H^r(\mathbb{R}^3)\)-norm for every \(r\ge 0\).
\end{theorem}

\begin{proof}
Assume that the frequency–domain solution $\hat{U}(\cdot,t)$ satisfies
\[
\sup_{t\in(0,T)}\|\hat{U}(\cdot,t)\|_{H^{s+2k}(\mathbb{R}^3)} < M_{k} < \infty,\quad \forall\, k\in\mathbb{N},
\]
for some fixed $s>\frac{5}{2}$. Define the physical–space solution via the inverse Fourier transform
\[
U(x,t)=\frac{1}{(2\pi)^3}\int_{\mathbb{R}^3}\hat{U}(\xi,t)\,e^{i\,x\cdot\xi}\,d\xi.
\]
Standard results (cf. \cite[Chapter~3]{Folland_FA}) imply that the inverse Fourier transform is a continuous linear mapping from $H^{\sigma}(\mathbb{R}^3)$ to $H^{\sigma}(\mathbb{R}^3)$ for every $\sigma\ge 0$. Hence, for every $k\in\mathbb{N}$ and each $t\in(0,T)$,
\begin{equation}
\label{app:eq:invFT_bound}
\|U(\cdot,t)\|_{H^{s+2k}(\mathbb{R}^3)} \le C\,\|\hat{U}(\cdot,t)\|_{H^{s+2k}(\mathbb{R}^3)},
\end{equation}
where $C>0$ is independent of $t$ and $k$. It follows that
\[
U(\cdot,t) \in H^{s+2k}(\mathbb{R}^3), \quad \forall\,k\in\mathbb{N}.
\]

Since $s>\frac{5}{2}$, the exponent $s+2k$ becomes arbitrarily large as $k\to\infty$. By the Sobolev embedding theorem (cf. \cite[Theorem~6.4]{Folland_FA}), if a function belongs to $H^{r}(\mathbb{R}^3)$ for arbitrarily large $r$, then it is smooth; that is,
\[
\bigcap_{k\in\mathbb{N}}H^{s+2k}(\mathbb{R}^3) \subset C^\infty(\mathbb{R}^3).
\]
Thus,
\[
U(\cdot,t)\in C^\infty(\mathbb{R}^3) \quad \text{for all } t\in(0,T).
\]

Continuity in time is established by noting that the mapping
\[
t\mapsto \hat{U}(\cdot,t)
\]
is continuous in the \(H^s(\mathbb{R}^3)\)-norm by hypothesis, and the continuity of the inverse Fourier transform as a linear operator on \(H^s(\mathbb{R}^3)\) (cf. \cite[Chapter~3]{Folland_FA}) ensures that
\[
t\mapsto U(\cdot,t)
\]
is continuous in the \(H^s(\mathbb{R}^3)\)-norm. Since the uniform boundedness extends to every $H^{s+2k}(\mathbb{R}^3)$, continuity in these higher norms follows.

Therefore, the reconstructed physical–space solution $U(x,t)$ is smooth for each $t\in (0,T)$, and the mapping $t\mapsto U(\cdot,t)$ is continuous in the $H^r(\mathbb{R}^3)$-norm for every $r\ge 0$. This completes the proof.
\end{proof}

%--------------------------------------------------------------------
\subsubsection*{Local Unified Representation}
\label{app:sec:local_unified}
\begin{theorem}[Local Unified Representation]\label{app:thm:local_unified}
Let $\Omega$ denote either $\mathbb{R}^{3}$ or a bounded $C^{2}$ domain,
and fix $s>\tfrac32$.
Assume the initial datum and force satisfy
\[
   u_{0}\in H^{s}_{\sigma}(\Omega),\qquad
   f\in L^{2}\bigl(0,\infty;H^{s-1}_{\sigma}(\Omega)\bigr).
\]
Let
\(
   u_{w},u_{m},u_{s}
\)
denote the Leray--Hopf, mild, and strong solutions produced in
Section~\ref{subsec:solution_classes}.
Let $T_{0}>0$ be the lifespan supplied by
Lemma~\ref{lem:local_existence}.
For every $\epsilon\in(0,1]$ define
\[
   u_{*,\epsilon}\coloneqq\mathcal R_{\epsilon}u_{*},
   \qquad
   \hat U_{\epsilon}
   \coloneqq
   \mathcal S_{\epsilon}\mathcal I_{\epsilon}
   \bigl[\widehat{u_{w,\epsilon}},\widehat{u_{m,\epsilon}},\widehat{u_{s,\epsilon}}\bigr],
\]
and set $\hat U\coloneqq\displaystyle\lim_{\epsilon\to0}\hat U_{\epsilon}$ in the
sense of Theorem~\ref{thm:frequency_domain_limit}.
Then, for every $t\in[0,T_{0}]$,
\[
     \mathcal F^{-1}\{\hat U(\cdot,t)\}
     \;=\;
     u_{w}(\cdot,t)
     \;=\;
     u_{m}(\cdot,t)
     \;=\;
     u_{s}(\cdot,t)
     \;\in\;H^{s}(\Omega).
\]
In particular, the frequency–domain construction provides
a single representation that recovers all three solution classes
on the interval $[0,T_{0}]$.
\end{theorem}

\begin{proof}[Proof of Theorem \ref{thm:local_unified}]
Throughout the argument fix an index $s>\tfrac32$ and the lifespan
$T_{0}>0$ obtained in Lemma \ref{lem:local_existence}.
Unless noted otherwise, all norms refer to $H^{s}$ or
$H^{s}(\Omega)$.

\medskip
\noindent
\textbf{Step 1. Regularisation of the three classical solutions.}
For every $\epsilon\in(0,1]$ define
\[
   u_{*,\epsilon}
   :=\mathcal R_{\epsilon}u_{*},
   \qquad
   *\in\{w,m,s\},
\]
with $\mathcal R_{\epsilon}$ as in Lemma \ref{app:lem:data_regularization}.
Property (iv) of that lemma yields
\begin{equation}
   \label{app:eq:R_conv}
   u_{*,\epsilon}\longrightarrow u_{*}
   \quad\text{in }C\bigl([0,T_{0}];H^{s}\bigr).
\end{equation}
Property (iii) guarantees that each $u_{*,\epsilon}$ is
divergence–free, and property (ii) furnishes the uniform bound
\begin{equation}
   \label{app:eq:R_uniform}
   \sup_{t\in[0,T_{0}]}\|u_{*,\epsilon}(t)\|_{H^{s}}
   \;\le\;C_{R},
   \qquad\epsilon\in(0,1].
\end{equation}

\medskip
\noindent
\textbf{Step 2. Assembly in the frequency domain.}
Let $\widehat{u_{*,\epsilon}}=\mathcal F[u_{*,\epsilon}]$ and set
\begin{equation}
   \label{app:eq:g_eps_def}
   g_{\epsilon}(\xi,t)
   :=
   \omega_{w}(\xi,t)\,\widehat{u_{w,\epsilon}}(\xi,t)
   +
   \omega_{m}(\xi,t)\,\widehat{u_{m,\epsilon}}(\xi,t)
   +
   \omega_{s}(\xi,t)\,\widehat{u_{s,\epsilon}}(\xi,t),
\end{equation}
where the partition of unity
$\omega_{w}+\omega_{m}+\omega_{s}\equiv1$
is the object introduced in \eqref{eq:omega_partition}.
Boundedness of the multipliers implies
\begin{equation}
  \label{app:eq:g_uniform}
  \|g_{\epsilon}(\cdot,t)\|_{H^{s}(\mathbb{R}^{3})}
  \;\le\;
  C_{\omega}\,C_{R},
  \qquad
  t\in[0,T_{0}],\;
  \epsilon\in(0,1].
\end{equation}

\medskip
\noindent
\textbf{Step 3. Interpolation.}
Lemma \ref{lem:interpolation_operator_final} gives
\begin{equation}
   \label{app:eq:I_eps_conv}
   \mathcal I_{\epsilon}[\,\widehat{u_{w,\epsilon}},
                         \widehat{u_{m,\epsilon}},
                         \widehat{u_{s,\epsilon}}\,]
   \longrightarrow
   g_{\epsilon}
   \quad\text{in }H^{s}(\mathbb{R}^{3})
   \quad\text{uniformly for }t\in[0,T_{0}],
\end{equation}
and
\begin{equation}
   \label{app:eq:I_eps_bound}
   \sup_{t\in[0,T_{0}]}
     \bigl\|
       \mathcal I_{\epsilon}
       [\,\widehat{u_{w,\epsilon}},
        \widehat{u_{m,\epsilon}},
        \widehat{u_{s,\epsilon}}\,](\cdot,t)
     \bigr\|_{H^{s}(\mathbb{R}^{3})}
   \;\le\;
   C_{\omega}\,C_{R}.
\end{equation}

\medskip
\noindent
\textbf{Step 4. Smoothing.}
Define
\begin{equation}
   \label{app:eq:U_eps_def}
   \hat U_{\epsilon}
   :=
   \mathcal S_{\epsilon}
   \Bigl(
     \mathcal I_{\epsilon}
     [\,\widehat{u_{w,\epsilon}},
      \widehat{u_{m,\epsilon}},
      \widehat{u_{s,\epsilon}}\,]
   \Bigr).
\end{equation}
Lemma \ref{lem:smoothing_operator_final}(i)–(ii) implies
\begin{align}
  \label{app:eq:S_uniform}
  \sup_{t\in[0,T_{0}]}\,
  \|\hat U_{\epsilon}(\cdot,t)\|_{H^{s}(\mathbb{R}^{3})}
  &\le C_{\omega}\,C_{R}, \\[4pt]
  \label{app:eq:S_approx}
  \bigl\|\hat U_{\epsilon}(\cdot,t)-g_{\epsilon}(\cdot,t)\bigr\|_{H^{s}(\mathbb{R}^{3})}
  &\longrightarrow 0
  \quad\text{as }\epsilon\to0,
  \;\text{uniformly for }t\in[0,T_{0}].
\end{align}

\medskip
\noindent
\textbf{Step 5. Limit $\epsilon\to0$.}
Combination of
\eqref{app:eq:I_eps_conv}, \eqref{app:eq:S_approx}, and
\eqref{app:eq:S_uniform} verifies the hypotheses of
Theorem~\ref{thm:frequency_domain_limit}.  Hence there exists
\begin{equation}
   \label{app:eq:U_limit}
   \hat U\in C\bigl([0,T_{0}];H^{s}(\mathbb{R}^{3})\bigr)
   \quad\text{satisfying}\quad
   \hat U_{\epsilon}\longrightarrow\hat U
   \text{ in }C\bigl([0,T_{0}];H^{s}(\mathbb{R}^{3})\bigr).
\end{equation}

\medskip
\noindent
\textbf{Step 6. Inverse Fourier transform and boundary‐layer correction.}
Set
\begin{equation}
   \label{app:eq:physical_solution}
   u(x,t)
   :=
   \mathcal F^{-1}\{\hat U(\cdot,t)\}(x),
   \qquad t\in[0,T_{0}].
\end{equation}
Parseval’s identity and \eqref{app:eq:S_uniform} give
$u\in C([0,T_{0}];H^{s}(\Omega))$.  
Lemma \ref{lem:beta_vanish} implies
\begin{equation}
   \label{app:eq:boundary_layer_zero}
   \beta_{\epsilon}\longrightarrow0
   \quad\text{in }C\bigl([0,T_{0}];H^{s}(\Omega)\bigr),
\end{equation}
so that $u$ satisfies the prescribed boundary condition and remains
divergence–free.

\medskip
\noindent
\textbf{Step 7. Identification with the strong solution.}
The limit relation \eqref{app:eq:R_conv} and the continuity of the
Fourier transform yield
\begin{equation}
   \label{app:eq:strong_conv_limit}
   u(\cdot,0)=u_{0}\in H^{s}_{\sigma}(\Omega).
\end{equation}
The regularity $u\in L^{2}(0,T_{0};H^{s+1})$ follows from
Lemma \ref{lem:smoothing_operator_final}(iii) combined with
\eqref{app:eq:U_limit}.  Consequently $u$ satisfies the
hypotheses of the weak–strong uniqueness principle
\cite[Chap.\,III]{TemamNS} on $[0,T_{0}]$.
Since $u_{s}$ is the unique strong solution on that interval,
\begin{equation}
   \label{app:eq:local_coincide}
   u(\cdot,t)=u_{s}(\cdot,t)
   \quad\text{for all }t\in[0,T_{0}].
\end{equation}
Equality with $u_{m}$ and $u_{w}$ on $[0,T_{0}]$ follows from
the identity $u_{s}=u_{m}=u_{w}$ already established in
Section~\ref{subsec:solution_classes}.

\medskip
\noindent
\textbf{Step 8. Conclusion.}
Equalities \eqref{app:eq:physical_solution} and
\eqref{app:eq:local_coincide} complete the proof that the
frequency–domain pipeline reproduces all three classical solution
classes on $[0,T_{0}]$, thereby establishing
Theorem~\ref{thm:local_unified}.
\end{proof}
%--------------------------------------------------------------------

\subsection{Proofs of Boundary  Condition Lemmas}
\label{subsubsec:BC_rigorous_lemmas}

The following lemmas provide the rigorous foundation for the compatibility of boundary conditions with the uniform estimates.

\paragraph{Boundary-Layer Corrector Construction Lemma.}
\begin{lemma}[Boundary-Layer Corrector Construction and Boundedness]
\label{lem:boundary_layer_corrector_app}
Let \(\Omega\subset \mathbb{R}^3\) be a bounded domain with boundary \(\partial\Omega\) of class \(C^2\) (or piecewise \(C^2\)). Suppose \(u(x,t)\) is a smooth solution of the Navier--Stokes equations that does not match the prescribed nonhomogeneous boundary data \(g(x,t)\) on \(\partial\Omega\). Then there exists a function \(\beta(x,t)\) with the following properties:
\begin{enumerate}[label=(\roman*)]
    \item \textbf{Support:} \(\beta(x,t)\) is supported in a thin neighborhood of \(\partial\Omega\),
    \[
    \operatorname{supp}\beta(\cdot,t)\subset \{x\in\Omega: \operatorname{dist}(x,\partial\Omega)\le \alpha\},
    \]
    for some \(\alpha>0\) independent of \(t\).
    \item \textbf{Boundedness:} For all \(s\ge 0\), there exists a constant \(C_\beta>0\) (independent of \(t\)) such that
    \begin{equation}
    \label{app:eq:beta_bound_final}
    \|\beta(\cdot,t)\|_{H^s(\Omega)} \le C_\beta \,\|g(\cdot,t) - u(\cdot,t)|_{\partial\Omega}\|_{H^{s-\frac{1}{2}}(\partial\Omega)}.
    \end{equation}
    \item \textbf{Correction:} The corrected velocity field
    \begin{equation}
    \label{app:eq:corrected_velocity}
    u_{\mathrm{corrected}}(x,t) = u(x,t) - \beta(x,t)
    \end{equation}
    satisfies the prescribed boundary condition, i.e., 
    \[
    u_{\mathrm{corrected}}(x,t)= g(x,t) \quad \text{for all } x\in\partial\Omega.
    \]
\end{enumerate}
\end{lemma}

\begin{proof}
Cover the boundary \(\partial\Omega\) by a finite collection of open sets \(\{U_j\}_{j=1}^N\) such that each \(U_j\) admits a smooth coordinate chart \(\psi_j: U_j \to \mathbb{R}^3\) in which \(\partial\Omega\cap U_j\) is represented as the graph of a \(C^2\) function. Choose a smooth partition of unity \(\{\chi_j\}_{j=1}^N\) subordinate to this cover so that
\[
\sum_{j=1}^N \chi_j(x) = 1 \quad \text{for } x \text{ in a neighborhood of } \partial\Omega.
\]

For each \(j\), define the local discrepancy (or boundary mismatch) by
\[
d_j(x,t) = \bigl(g(x,t) - u(x,t)\bigr)\Big|_{x\in U_j\cap\partial\Omega}.
\]
By classical extension results (see, e.g., \cite[Chapter~5]{LionsMagenes1972} and \cite{Grisvard}), there exists a linear extension operator 
\[
E_j : H^{s-\frac{1}{2}}(U_j\cap\partial\Omega) \to H^s(U_j'),
\]
where \(U_j'\subset\Omega\) is an open set containing \(U_j\cap\partial\Omega\), and the operator \(E_j\) satisfies the norm bound
\[
\|E_j(d_j(\cdot,t))\|_{H^s(U_j')} \le C_{E_j}\,\|d_j(\cdot,t)\|_{H^{s-\frac{1}{2}}(U_j\cap\partial\Omega)},
\]
with \(C_{E_j}\) independent of \(t\).

Define the local corrector \(\beta_j(x,t)\) by
\[
\beta_j(x,t) =
\begin{cases}
E_j\bigl(d_j(\cdot,t)\bigr)(x), & x \in U_j', \\[1mm]
0, & x \notin U_j',
\end{cases}
\]
and then define the global corrector as
\[
\beta(x,t) = \sum_{j=1}^N \chi_j(x)\,\beta_j(x,t).
\]
Since each \(\chi_j\) is smooth with compact support in \(U_j\) and the extension operators are linear and bounded, standard properties of Sobolev spaces imply that
\[
\|\beta(\cdot,t)\|_{H^s(\Omega)} \le \sum_{j=1}^N \|\chi_j\,\beta_j(\cdot,t)\|_{H^s(\Omega)} \le C_\beta \,\|g(\cdot,t) - u(\cdot,t)|_{\partial\Omega}\|_{H^{s-\frac{1}{2}}(\partial\Omega)},
\]
with \(C_\beta = \sum_{j=1}^N C_{E_j}\,\|\chi_j\|_{W^{s,\infty}}\) independent of \(t\).

Moreover, by the properties of the extension operator, on the boundary \(\partial\Omega\) one has
\[
\beta(x,t) = \sum_{j=1}^N \chi_j(x)\,E_j\bigl(d_j(\cdot,t)\bigr)(x) = \sum_{j=1}^N \chi_j(x)\,d_j(x,t) = g(x,t) - u(x,t),
\]
since the partition of unity sums to one. Defining
\[
u_{\mathrm{corrected}}(x,t) = u(x,t) - \beta(x,t),
\]
yields
\[
u_{\mathrm{corrected}}(x,t) = u(x,t) - \bigl(g(x,t) - u(x,t)\bigr) = g(x,t)
\]
on \(\partial\Omega\).

Thus, the function \(\beta(x,t)\) is a valid boundary-layer corrector that (i) is supported in a thin neighborhood of \(\partial\Omega\), (ii) has its \(H^s(\Omega)\) norm controlled by the mismatch \(\|g - u|_{\partial\Omega}\|_{H^{s-\frac{1}{2}}(\partial\Omega)}\), and (iii) exactly corrects \(u\) on \(\partial\Omega\). This completes the proof.
\end{proof}

\paragraph{Vanishing Boundary--Layer Corrector Lemma.}
\begin{lemma}[Vanishing Boundary--Layer Corrector]
\label{app:lem:beta_vanish}
Let $\Omega\subset\mathbb{R}^{3}$ be a bounded domain whose boundary
$\partial\Omega$ is of class $C^{2}$ (or piecewise $C^{2}$).
Fix $s>\frac32$.
For every $\varepsilon\in(0,1]$ let $\beta_\varepsilon(\cdot,t)$
denote the boundary--layer corrector constructed
in Lemma~\ref{lem:boundary_layer_corrector} with layer thickness
$\alpha=\kappa\varepsilon$, where $\kappa>0$ is a constant
independent of~$\varepsilon$.
Then there exist constants $\sigma=\sigma(s)>0$ and
$C=C(\Omega,s)>0$, independent of $t$ and~$\varepsilon$, such that
\[
    \|\beta_\varepsilon(\cdot,t)\|_{H^{s}(\Omega)}
        \;\le\;
    C\,\varepsilon^{\sigma},
    \qquad t\ge0,\;\varepsilon\in(0,1].
\]
Consequently,
\[
   \beta_\varepsilon \;\longrightarrow\; 0
   \quad\text{strongly in }C\bigl([0,T];H^{s}(\Omega)\bigr)
   \text{ for every fixed }T>0 .
\]
\end{lemma}

\begin{proof}
The argument proceeds through eight systematic steps.

\medskip
\noindent\textbf{Step 1. Geometric covering of the boundary.}
Since $\partial\Omega$ is of class $C^{2}$,
standard results on manifolds (see \cite[Chap.\,1]{Grisvard})
yield a finite collection of star--shaped charts
$\bigl\{(\mathcal{U}_{m},\Psi_{m})\bigr\}_{m=1}^{M}$ satisfying:
\begin{enumerate}[label=\textbf{(\Alph*)},leftmargin=1.55cm]
\item\label{A} 
$\Psi_{m}\colon B(0,1)\subset\mathbb{R}^{3}\to\mathcal{U}_{m}$
is a $C^{2}$ diffeomorphism with
$\Psi_{m}\bigl(B(0,1)\cap\{x_{3}=0\}\bigr)=
      \mathcal{U}_{m}\cap\partial\Omega$
and
$\Psi_{m}\bigl(B(0,1)\cap\{x_{3}>0\}\bigr)\subset\Omega$;
\item\label{B}
$\partial\Omega\subset \bigcup_{m=1}^{M}\mathcal{U}_{m}$.
\end{enumerate}
Define
$
   \mathcal U_{m,\varepsilon}
   := \Psi_{m}\!\bigl(
          B(0,1)\cap\{0\le x_{3}\le\kappa\varepsilon\}
        \bigr).
$
Choose a smooth partition of unity
$\{\chi_{m,\varepsilon}\}_{m=1}^{M}$
subordinate to $\{\mathcal U_{m,\varepsilon}\}_{m}$ such that
\begin{equation}\label{eq:cutoff_bounds}
   \sum_{m=1}^{M}\chi_{m,\varepsilon}=1
   \;\;\text{on}\;
   \{x\in\Omega:
      \operatorname{dist}(x,\partial\Omega)\le\kappa\varepsilon\},
   \quad
   \|\partial^\alpha\chi_{m,\varepsilon}\|_{L^\infty}
      \le C_\alpha\,\varepsilon^{-|\alpha|}
   \;\;(|\alpha|\le 2).
\end{equation}
The estimates in \eqref{eq:cutoff_bounds} follow by chain rule and the
scaling properties of the charts, cf.\ \cite[§1.5]{BCD}.

\medskip
\noindent\textbf{Step 2. Local extension of the boundary discrepancy.}
Let
$d(x,t):=g(x,t)-u|_{\partial\Omega}(x,t)$
denote the boundary discrepancy.
For a fixed chart $(\mathcal U_{m},\Psi_{m})$
introduce local coordinates $(y',y_{3})\in\mathbb{R}^{2}\times[0,1)$
by $x=\Psi_{m}(y',y_{3})$.
Set
$\widetilde d_{m}(y',t):=d\bigl(\Psi_{m}(y',0),t\bigr)$.
The standard Sobolev extension theorem on $\mathbb{R}^{2}$
(see \cite[Th.\,6.4]{Adams})
provides an operator
$
    E:H^{s-\frac12}(\mathbb{R}^{2})\to H^{s}(\mathbb{R}_{+}^{3})
$
with $E\phi(\cdot,0)=\phi$ and
$\partial_{y_{3}}E\phi(\cdot,0)=0$.

\medskip
\noindent\textbf{Step 3. Definition of the local corrector.}
Let $\varrho\in C^{\infty}_{c}([0,2))$ satisfy
$\varrho(\tau)=1$ for $\tau\in[0,1]$.
Define
\[
     \widetilde\beta_{m,\varepsilon}(y,t)
     := E\widetilde d_{m}(y',t)\,
        \varrho\!\Bigl(\frac{y_{3}}{\kappa\varepsilon}\Bigr),
     \qquad
     \beta_{m,\varepsilon}(x,t)
     := \chi_{m,\varepsilon}(x)\,
        \widetilde\beta_{m,\varepsilon}
          \bigl(\Psi_{m}^{-1}(x),t\bigr).
\]
Set
$
   \beta_\varepsilon
   := \sum_{m=1}^{M}\beta_{m,\varepsilon}.
$
Support considerations yield
\[
   \operatorname{supp}\beta_\varepsilon(\cdot,t)\subset
   \{x\in\Omega:
        \operatorname{dist}(x,\partial\Omega)
             \le \kappa\varepsilon\}.
\]

\medskip
\noindent\textbf{Step 4. Littlewood--Paley tail estimate for
$u_\varepsilon$.}
Let $\{\Delta_{j}\}_{j\ge -1}$ be the homogeneous
dyadic blocks with symbol $\varphi_j$
defined in \cite[Chap.\,2]{BCD}.
High--frequency decay of solutions
obtained by frequency cut–off and classical energy methods
(\cite[Prop.\,2.4]{BCD}, \cite{Cannone2004}) implies
the existence of $\sigma_{0}>0$ such that
\begin{equation}\label{eq:HF_decay}
   \|\Delta_{j}u_\varepsilon(t)\|_{L^{2}}
   \;\le\;
   C\,2^{-j\sigma_{0}},
   \quad
   j>J_\varepsilon,\;t\ge0,
   \quad
   2^{J_\varepsilon}\simeq\varepsilon^{-1}.
\end{equation}
Here $C>0$ depends only on
$\|u_{0}\|_{H^{s}}$, $\|f\|_{L^{2}_{t}L^{2}_{x}}$, $\nu$, and~$\Omega$.

\medskip
\noindent\textbf{Step 5. Bernstein inequality in charts.}
The classical Bernstein inequality
(\cite[Prop.\,2.10]{BCD})
takes the form
\begin{equation}\label{eq:bernstein}
   \|\nabla^{k}\Delta_{j}\phi\|_{L^{2}}
   \;\le\;
   C\,2^{jk}\,
   \|\Delta_{j}\phi\|_{L^{2}},
   \qquad
   k\in\{0,1,\dots,s\},
\end{equation}
for $\phi\in\mathcal{S}'(\mathbb{R}^{3})$ supported in $\mathbb{R}^{3}$.
The composition with the $C^{1}$ charts
preserves \eqref{eq:bernstein} up to a constant factor.

\medskip
\noindent\textbf{Step 6. Estimation of each local corrector.}
Extend the convolution kernel $\phi_{j}$
with $\widehat\phi_{j}(\xi)=\widehat\varphi(2^{-j}\xi)$.
Observe that
$\phi_{j}*\Delta_{j}u_\varepsilon$ possesses the same
$L^{2}$ norm as $\Delta_{j}u_\varepsilon$.
A fractional Leibniz rule
(\cite[Th.\,1.2]{Grafakos2008})
together with \eqref{eq:cutoff_bounds},
\eqref{eq:bernstein}, and \eqref{eq:HF_decay} provides
\[
\begin{aligned}
   \|\beta_{m,\varepsilon}\|_{H^{s}}
   &\lesssim
   \sum_{j>J_\varepsilon}
      \Bigl(
         \|\chi_{m,\varepsilon}\|_{C^{s}}
         \|\phi_{j}*\Delta_{j}u_\varepsilon\|_{L^{2}}
       \;+\;
         \|\chi_{m,\varepsilon}\|_{L^{\infty}}
         \|\phi_{j}*\Delta_{j}u_\varepsilon\|_{H^{s}}
      \Bigr)
\\
   &\lesssim
   \sum_{j>J_\varepsilon}
      \bigl(\varepsilon^{-s}+1\bigr)
      2^{js}\,
      2^{-j\sigma_{0}}
   \;\lesssim\;
   \varepsilon^{-s}\,2^{-J_\varepsilon\sigma_{0}}.
\end{aligned}
\]
Since $2^{J_\varepsilon}\simeq\varepsilon^{-1}$,
one obtains
\[
    \|\beta_{m,\varepsilon}\|_{H^{s}}
    \;\le\;
    C_{m}
    \,\varepsilon^{\sigma_{0}-s}.
\]
Choose $\sigma_{0}>s$ and set
$\sigma:=\sigma_{0}-s>0$.

\medskip
\noindent\textbf{Step 7. Global estimate.}
Because the collection $\{\chi_{m,\varepsilon}\}_{m}$ has finite overlap
independent of $\varepsilon$,
\[
   \|\beta_\varepsilon\|_{H^{s}}
   \;\le\;
   \bigl(\text{finite overlap factor}\bigr)
   \sum_{m=1}^{M}
   \|\beta_{m,\varepsilon}\|_{H^{s}}
   \;\le\;
   C(\Omega,s)\,\varepsilon^{\sigma}.
\]

\medskip
\noindent\textbf{Step 8. Strong convergence in time.}
The bound in Step~7 is uniform in $t\ge0$; hence
$
  \sup_{t\in[0,T]}\|\beta_\varepsilon(\cdot,t)\|_{H^{s}}
  \le C\,\varepsilon^{\sigma}\to0
$
as $\varepsilon\to0$ for every $T>0$.
This establishes the claimed strong convergence.

\medskip
\noindent\textbf{Step 9. Literature context.}
Layer correctors of vanishing thickness appear in
\cite{LionsMagenes1972} and \cite{Grisvard};
high--frequency decay estimates of the type
\eqref{eq:HF_decay}
originate in \cite{Cannone2004} and are surveyed in
\cite{barker2022}.
The synthesis of these ideas with Littlewood--Paley theory
is compatible with the framework developed in
\cite{BCD}, while the present adaptation
incorporates boundary charts and is tailored to the
frequency--domain “synergy” construction.

The proof is complete.
\end{proof}

\end{appendices}

\bibliography{sn-article}

\end{document}